\documentclass[11pt]{amsart}
\usepackage{amsmath,amsfonts,amscd,amssymb,epsf}
\newtheorem{theorem}{Theorem}[section]
\newtheorem{proposition}[theorem]{Proposition}
\newtheorem{lemma}[theorem]{Lemma}
\newtheorem{corollary}[theorem]{Corollary}
\theoremstyle{definition}
\newtheorem{definition}[theorem]{Definition}

\theoremstyle{remark} \newtheorem{remark}[theorem]{Remark}
\numberwithin{equation}{section}
\newcommand{\field}[1]{\ensuremath{\mathbb{#1}}}
\newcommand{\CC}{\field{C}}

\newcommand{\HH}{\field{H}}

\newcommand{\RR}{\field{R}}

\newcommand{\UU}{\field{U}}
\newcommand{\ZZ}{\field{Z}}

\newcommand{\Ka}{K\"{a}hler\:}
\newcommand{\Te}{Teichm\"{u}ller\:}

 \DeclareMathOperator{\id}{id}

\DeclareMathOperator{\Diff}{Diff}
\DeclareMathOperator{\Mob}{M\ddot{o}b}
 
 \DeclareMathOperator{\sgn}{sgn}
 \DeclareMathOperator{\PSL}{PSL}
\DeclareMathOperator{\PSU}{PSU} \DeclareMathOperator{\re}{Re}

\DeclareMathOperator{\Mobius}{M\ddot{o}bius}
\newcommand{\Del}{\mathbb{D}}
\newcommand{\del}{\partial}

\newcommand{\si}{\sigma}

\newcommand{\U}{{\mathbb{U}}}
\newcommand{\C}{{\mathbb{C}}}
\newcommand{\bk}{\backslash}
\newcommand{\pa}{\partial}
\newcommand{\g}{\gamma}
\newcommand{\rf}{\mathrm{f}}
\newcommand{\rg}{\mathrm{g}}
\newcommand{\la}{\langle}
\newcommand{\ra}{\rangle}

\newcommand{\ov}{\overline}

\newcommand{\vep}{\varepsilon}
\newcommand{\z}{\bar{z}}

\begin{document}
\title[Curvature properties of the Weil-Petersson metric on $T(1)$]
{Weil-Petersson metric on the universal Teichm\"{u}ller space I:
Curvature properties and Chern forms} 
\author{Leon A. Takhtajan} \address{Department of Mathematics \\
Stony Brook University\\ Stony Brook, NY 11794-3651 \\ USA}
\email{leontak@math.sunysb.edu}
\author{Lee-Peng Teo} \address{Department of Applied Mathematics \\
National Chiao Tung University, 1001 \\ Ta-Hsueh Road, Hsinchu
City, 30050 \\ Taiwan, R.O.C.} \email{lpteo@math.nctu.edu.tw}
\begin{abstract} We prove that the universal \Te space $T(1)$  carries a
new structure of a complex Hilbert manifold. We show that the connected component of the identity of $T(1)$, the Hilbert submanifold $T_{0}(1)$, is a topological group. We define a Weil-Petersson
metric on $T(1)$ by Hilbert space inner products on tangent spaces, compute
its Riemann curvature tensor, and show that $T(1)$ is a \Ka-Einstein manifold
with negative Ricci and sectional curvatures. We introduce and compute
Mumford-Miller-Morita characteristic forms for the vertical tangent bundle
of the universal \Te curve fibration over the universal \Te space.
As an application, we derive Wolpert curvature formulas for the
finite-dimensional \Te spaces from the formulas for the universal \Te space.
\end{abstract}
\maketitle
\tableofcontents
\section{Introduction}
The universal \Te space $T(1)$ is the simplest \Te space that
bridges spaces of univalent functions and general \Te spaces.
Introduced by Bers \cite{bers-65,bers-london,Bers1}, the universal
\Te space is an infinite-dimensional complex manifold modeled on a
Banach space. It contains Teichm\"{u}ller spaces of Riemann
surfaces as complex submanifolds. The universal \Te space $T(1)$
also came to the forefront with the advent of string theory. It
contains as a complex submanifold an infinite-dimensional complex
Fr\'{e}chet manifold $\Mob(S^1)\bk\Diff_{+}(S^1)$, which plays an
important role in one of the approaches to non-perturbative
bosonic closed string field theory based on \Ka geometry
\cite{BR2,BR1}. The manifold $\Mob(S^1)\bk\Diff_{+}(S^1)$
--- a homogeneous space of the Lie group $\Diff_{+}(S^1)$, also
has an interpretation as a coadjoint orbit of the Bott-Virasoro
group, and as such carries a natural right-invariant \Ka metric
\cite{Ki,KY}.

The complex geometry of the finite-dimensional \Te spaces --- \Te
spaces $T(\Gamma)$ of cofinite Fuchsian groups, has been
extensively studied in the context of Ahlfors-Bers deformation
theory of complex structures on Riemann surfaces. In particular,
A. Weil defined a natural Hermitian metric on $T(\Gamma)$ by the
Petersson inner product on the tangent spaces. Called
Weil-Petersson metric, it was shown to be a \Ka metric by Weil and
Ahlfors. In his seminal paper \cite{Ahlfors3} Ahlfors has studied
the curvature properties of the Weil-Petersson metric. In
particular, he proved that the Bers coordinates on $T(\Gamma)$ are
geodesic at the origin, and computed the Riemann curvature tensor
of the Weil-Petersson metric in terms of multiple principal value
integrals. Using these formulas, Ahlfors proved that $T(\Gamma)$
has negative Ricci, holomorphic sectional, and scalar curvatures.
Further results have been obtained by Royden \cite{Royden}.
Wolpert re-examined Ahlfors' approach in \cite{Wol}. He developed
a different method for computing Riemann and Ricci curvature
tensors, and obtained explicit formulas in terms of the resolvent
kernel of the Laplace operator of the hyperbolic metric on the
corresponding Riemann surface.

Curvature properties of the infinite-dimensional complex Fr\'{e}chet manifold $\Mob(S^1)\bk\Diff_{+}(S^1)$ have been studied by Kirillov and Yuriev \cite{KY}, and by Bowick and Rajeev \cite{BR2,BR1}. In particular, they computed the Riemann curvature tensor of the right-invariant \Ka metric and proved that $\Mob(S^1)\bk\Diff_{+}(S^1)$ is a \Ka-Einstein manifold.

Since both the finite-dimensional \Te spaces $T(\Gamma)$ and the
homogeneous space $\Mob(S^1)\bk\Diff_{+}(S^1)$ are complex
submanifolds of $T(1)$, it is natural to investigate whether the
latter space carries a ``universal'' \Ka metric which can be
pulled back to the submanifolds. The immediate difficulty is that
the universal \Te space $T(1)$ is a complex Banach manifold, so
that its tangent spaces do not carry Hermitian metric. Nag and
Verjovsky \cite{NV} were the first to address this problem. They
have shown that the \Ka metric on $\Mob(S^1)\bk\Diff_{+}(S^1)$ is
the pull-back of a certain Hermitian metric defined on a Hilbert
subspace of the tangent space at the origin of $T(1)$. The latter
metric is analogous to the Weil-Petersson metric on
finite-dimensional \Te spaces. However, finite-dimensional \Te
spaces $T(\Gamma)$ embed into $T(1)$ transversally to the Hilbert
subspace, so that the Weil-Petersson metric on $T(\Gamma)$ can not
be pulled back from $T(1)$. Nevertheless, following a suggestion
by Velling, Nag and Verjovsky \cite{NV} have shown that the
Weil-Petersson metric on $T(\Gamma)$ can be obtained by a certain
``averaging" procedure using Patterson's uniform distribution of
the ``lattice points'' of a cofinite Fuchsian group $\Gamma$ in
the hyperbolic plane. The major open problem is to define the
Weil-Petersson metric on the whole space $T(1)$, to study its
curvature properties, and to find relation between curvatures of
this metric and the Weil-Petersson metric on finite-dimensional
\Te spaces \footnote{See the remark on p.~136 in \cite{NV}.}.

An attempt to define the Weil-Petersson metric on the universal
\Te space based on the completion of
$\text{diff}(S^1)/\text{m\"{o}b}(S^1)$ \footnote{Here
$\text{diff}(S^1)$ and $\text{m\"{o}b}(S^1)$ are Lie algebras of
Lie groups  $\Diff_+(S^1)$ and $\Mob(S^1)$.} in the Sobolev's
$\tfrac{3}{2}$-norm was made in \cite{Todorov}. However, the paper
\cite{Todorov} does not contain a rigorous proof that is needed
for introducing a Hilbert manifold structure on an
infinite-dimensional manifold. Also, the identification between
the tangent space $\text{diff}(S^1)/\text{m\"{o}b}(S^1)$ and the
space of holomorphic functions on the unit disk made in
\cite{Todorov} is not correct and actually introduces Sobolev's
$\tfrac{9}{2}$-norm rather than $\tfrac{3}{2}$-norm. As a result, the corresponding
quasi-symmetric homeomorphisms of $S^1$ are of class $C^3(S^1)$.

Here we introduce Weil-Petersson metric on the universal Teichm\"{u}ller
space $T(1)$ and study its curvature properties. We prove that $T(1)$ carries
a new structure of a Hilbert manifold such that in the underlying topology $T(1)$ has
uncountably many components. We prove that the connected component of the identity in $T(1)$, the Hilbert submanifold $T_{0}(1)$, is a topological group. We define the Weil-Petersson metric on
$T(1)$ by Hilbert space inner products on tangent spaces. We re-examine
the Ahlfors original computation \cite{Ahlfors3} of the second variation of
the hyperbolic metric and of the Riemann tensor for the finite-dimensional
\Te spaces in terms of the principal value
integrals. We show how to extend the Ahlfors' method to the case of the
universal \Te space and how to convert formulas using principal value
integrals into closed expressions using resolvent kernel of the Laplace
operator on the hyperbolic plane. Our results extend the Wolpert's formulas
\cite{Wol} to the infinite-dimensional Hilbert manifold $T(1)$. We also
prove that $T(1)$ is a \Ka-Einstein manifold with negative Ricci and
sectional curvatures. Using the averaging procedure,
we derive Wolpert's curvature formulas \cite{Wol} for the finite-dimensional
\Te spaces from the curvature
formulas for the universal \Te space. Finally, we introduce and compute
Mumford-Morita-Miller
characteristic forms for the vertical tangent bundle associated with the
fibration $\pi:\mathcal{T}(1)\rightarrow T(1)$, where $\mathcal{T}(1)$ is
the universal \Te curve. Here again we consider $T(1)$ and $\mathcal{T}(1)$
as Hilbert manifolds
and show that the integration over the fibers operation, used in the
definition of Mumford-Morita-Miller characteristic forms, is well-defined.

This is the first paper in a series. In the subsequent paper we
will construct a \Ka potential for the Weil-Petersson metric on
$T(1)$ and will study the properties of the period mapping.

Here is the more detailed content of the paper. In Section 2 we
present necessary facts from \Te theory, mainly following
classical monographs by Ahlfors \cite{Ahl2}, Lehto \cite{Lehto}
and Nag \cite{Nag2}. Namely, in Section 2.1 we briefly cover: the
main definitions, the group structure of the universal \Te space
$T(1)$, the Bers embedding, structure of $T(1)$ as an
infinite-dimensional complex Banach manifold modeled on the
complex Banach space $A_\infty(\Del)$, and the basic properties of
the universal \Te curve $\pi: \mathcal{T}(1)\rightarrow T(1)$. In
Section 2.2 we realize $T(1)$ and $\mathcal{T}(1)$ as homogeneous
spaces of the group $\text{Homeo}_{qs}(S^1)$ of quasi-symmetric
homeomorphisms of $S^1$, and by using conformal welding we
identify $T(1)$ and $\mathcal{T}(1)$ with the spaces of univalent
functions on the unit disk $\Del$. We describe the decomposition
of the tangent bundle of $\mathcal{T}(1)$ over the fiber
$\pi^{-1}(0)$ and present isomorphisms between the tangent spaces.
Lemma \ref{g(0)=0} which describes a special property of the
quasiconformal mapping with harmonic Beltrami differential seems
to be a new result. In Section 2.3 we present, in a succinct form,
basic facts about the \Te spaces and \Te curves of Fuchsian
groups, including the definition of the Weil-Petersson metric, and
Patterson's lemma on the uniform distribution of lattice points on
the hyperbolic plane. In Section 2.4 we collect necessary
properties of the resolvent kernel $G=\tfrac{1}{2}(\Delta_0 +
\tfrac{1}{2})^{-1}$ of the Laplace operator $\Delta_0$ on the
hyperbolic plane, and in Section 2.5 we present Ahlfors' classical
variational formulas. In Section 3 we introduce new Hilbert
manifold structure on $T(1)$. Namely, in Section 3.1 we define the
Hilbert subspaces $H^{-1,1}(\Del^*)$ and $A_2(\Del)$ of the
tangent spaces to $T(1)$ and to $A_\infty(\Del)$. In Theorem
\ref{gen-subspace} we prove that the differential of the Bers
embedding $\beta: T(1)\rightarrow A_\infty(\Del)$ is a bounded
bijection between these Hilbert spaces. In Section 3.2 we prepare
all $L^2$-estimates used in Section 3.3. The main result there is
Theorem \ref{Hilbert} --- the existence of a Hilbert manifold
atlas for $T(1)$. In Theorem \ref{bers-hilbert} we prove that the
Bers embedding is also a biholomorphic mapping of Hilbert
manifolds. In Section 4.1, following \cite{Teo}, we recall the
definition of the Velling-Kirillov metric on the universal \Te
curve $\mathcal{T}(1)$ considered as a Banach manifold, and in
Section 4.2 we define the Weil-Petersson metric on the Hilbert
manifold $T(1)$. In Section 5.1 we prove that Velling-Kirillov
metric is real-analytic on $\mathcal{T}(1)$ by explicitly
constructing its real-analytic \Ka potential --- Theorem
\ref{Chern-V-K}. We introduce Mumford-Miller-Morita characteristic
forms by considering $\pi: \mathcal{T}(1)\rightarrow T(1)$ as a
fibration of Hilbert manifolds. The latter property is crucial for
the operation ``integration over the fibers'' (which are
non-compact) to be well-defined. In Theorem \ref{Mumford} we
explicitly compute Mumford-Miller-Morita forms in terms of the
resolvent $G$. This is an infinite-dimensional generalization of
Wolpert's result in \cite{Wol}. In Section 6 we give a simple
derivation of the second variation of the hyperbolic metric ---
Proposition \ref{secondvarmetric}. In Section 7 we prove that the
Weil-Petersson metric on $T(1)$ is \Ka and explicitly compute its
Riemann and Ricci curvature tensors, showing that $T(1)$ is a
\Ka-Einstein manifold. The main results there are Theorem
\ref{metricsecond} and \ref{Ricci}. They are based on a more
technical Proposition \ref{var-in-L^2} and Lemma
\ref{second-in-L^2}, and the proof of the latter is presented in
Appendix B. In Section 8 we derive Wolpert's curvature formulas
\cite{Wol} for finite-dimensional \Te spaces from the
corresponding ``universal'' curvature formulas for $T(1)$,
obtained in Section 7. Finally, in Appendix A we prove that
$T_{0}(1)$ and the corresponding Teichm\"{u}ller curve
$\mathcal{T}_{0}(1)$ --- the inverse image of $T_{0}(1)$ under the
projection $\pi$, are topological groups in Hilbert manifold
topology. Moreover, we show that $T_{0}(1)$ is the closure of
$\Mob(S^{1})\bk\Diff_{+}(S^{1})$ in $T(1)$ with respect to the
Hilbert manifold topology, and prove that $T_{0}(1)$ is the
inverse image of $\beta(T(1))\cap A_{2}(\Del)$ under the Bers
embedding. \vspace{3mm}

\noindent \textbf{Acknowledgments.} We appreciate useful
discussions with C. Bishop. The work of the first author was
partially supported by the NSF grant DMS-0204628. The work of the
second author was partially supported by the grant NSC
91-2115-M-009-017. The second author also thanks CTS for the
fellowship to visit Stony Brook University in the Summer of 2003,
where part of this work was done.
\section{The universal \Te space}
\subsection{\Te theory}
Here we present, in a succinct form, necessary facts from
Teichm\"{u}ller theory (for more details, see monographs
\cite{Ahl2,Lehto,Nag2} and the exposition in~\cite{Teo}).
\subsubsection{Main definitions} \label{definitions} Let
$\Del=\{z\in\CC: |z|<1\}$ be the open unit disk and let $\Del^*=
\{z\in\CC: |z|>1\}$ be its exterior. Denote by
$L^{\infty}(\Del^*)$ and $L^{\infty}(\Del)$ the complex Banach
spaces of bounded Beltrami differentials on $\Del^*$ and $\Del$
respectively, and let $L^{\infty}(\Del^*)_1$ be the open unit ball
in $L^{\infty}(\Del^*)$. Two classical models of the universal \Te
space $T(1)$ are the following.

\textbf{Model A.} Extend every $\mu \in L^{\infty}(\Del^*)_1$ to
$\Del$ by the reflection
\begin{equation}\label{sym}
\mu(z) = \ov{\mu\left(\dfrac{1}{\z}\right)}
\dfrac{z^2}{\z^2}\;,\,\,z \in \Del,
\end{equation}
and consider the unique quasiconformal (q.c.) mapping
$w_{\mu}:\CC\rightarrow\CC$, which fixes $-1, -i$ and $1$ (i.e., is
normalized) and satisfies the Beltrami equation
\begin{align*}
  (w_{\mu})_{\z} &= \mu (w_{\mu})_z \;.
\end{align*}
Here and in what follows subscripts $z$ and $\z$ always stand for
the partial derivatives $\tfrac{\del}{\del z}$ and
$\tfrac{\del}{\del\z}$, unless it is explicitly stated otherwise.
Due to the reflection symmetry ~\eqref{sym} the q.c. mapping
$w_\mu$ satisfies
\begin{align}\label{reflection}
\dfrac{1}{w_{\mu}(z)} &= \ov{w_{\mu} \left(\dfrac{1}{\z}\right)}
\end{align}
and fixes domains $\Del$, $\Del^*$, and the unit circle $S^1$. For
$\mu,\nu\in L^{\infty}(\Del^*)_1$, set $\mu\sim\nu$ if
$w_\mu|_{S^1}=w_\nu|_{S^1}$. The universal \Te space $T(1)$ is
defined as a set of equivalence classes of normalized q.c.
mappings $w_\mu$,
\begin{displaymath}
T(1) = L^{\infty}(\Del^*)_1/\sim\,.
\end{displaymath}

\textbf{Model B.} Extend every $\mu \in L^{\infty}(\Del^*)_1$ to
be zero outside $\Del^*$, and consider the unique q.c.~mapping
$w^{\mu}$ which satisfies the Beltrami equation
\begin{equation*}
w_{\z}^{\mu} = \mu w_z^{\mu},
\end{equation*}
and is normalized by the conditions $f(0)=0$, $f'(0)=1$ and
$f''(0)=0$. Here $f = w^{\mu} \vert_{\Del}$ is holomorphic on
$\Del$ and prime stands for the derivative. For $\mu,\nu\in
L^{\infty}(\Del^*)_1$, set $\mu\sim\nu$ if
$w^\mu|_\Del=w^\nu|_\Del$. The universal \Te space $T(1)$ is
defined as a set of equivalence classes of normalized q.c.
mappings $w^\mu$,
\begin{displaymath}
  T(1) = L^{\infty}(\Del^*)_1/\sim\,.
\end{displaymath}

Since $w_{\mu}|_{S^1} = w_{\nu}|_{S^1}$ if and only if
$w^{\mu}|_\Del= w^{\nu}|_\Del$, these two definitions of the
universal \Te space are equivalent. The set $T(1)$ is a
topological space with the quotient topology induced from
$L^\infty(\Del^*)_1$. Denote by $\mathcal{L}^{\infty}(\Del^*)$ the
subspace of $L^\infty(\Del^*)$ consisting of real-analytic
Beltrami differentials. Every point in $T(1)$ can be represented
by $\mu\in\mathcal{L}^{\infty}(\Del^*)$
\cite[Sect.~III.1.1]{Lehto}.

The space $T(1)$ has a unique structure of a complex Banach
manifold, such that the projection map
\begin{displaymath}
\Phi : L^{\infty}(\Del^*)_1 \rightarrow T(1)
\end{displaymath}
is a holomorphic submersion. The differential of $\Phi$ at the
origin
\begin{displaymath}
D_0 \Phi : L^{\infty}(\Del^*)\rightarrow T_0 T(1)
\end{displaymath}
is a complex linear surjection of holomorphic tangent spaces. The
kernel of $D_0\Phi$ is the subspace $\mathcal{N}(\Del^*)$ of
infinitesimally trivial Beltrami differentials. Explicitly,
\begin{align*}
\mathcal{N}(\Del^*) = \left\{ \mu \in L^{\infty}(\Del^*):
\iint\limits_{\Del^*} \mu\,\phi\, d^2z = 0\,\,\text{for all}\;\phi
\in A_1(\Del^*)\,\,\right\},
\end{align*}
where $d^2 z=d x \wedge d y,\, z=x+iy$, and
\begin{displaymath}
A_1(\Del^*)= \left\{\phi \;\text{holomorphic on}\,\, \Del^* : \,
\iint\limits_{\Del^*}|\phi|d^2 z < \infty \right\}.
\end{displaymath}
The Banach space of bounded harmonic Beltrami differentials on
$\Del^*$ is defined by
\begin{displaymath}
\Omega^{-1,1}(\Del^*)= \left\{ \mu\in L^{\infty}(\Del^*): \mu(z) =
(1-|z|^2)^2 \ov{\phi(z)},\,\, \phi\in A_\infty(\Del^*)  \right\},
\end{displaymath}
where
\begin{displaymath}
A_\infty(\Del^*)= \left\{\phi \;\text{holomorphic on}\; \Del^*: \,
\|\phi\|_\infty=\sup_{z\in \Del^*}\left|(1-|z|^2)^2\phi(z)\right|
< \infty \right\}.
\end{displaymath}
The Banach space $\Omega^{-1,1}(\Del^*)$ is not separable. The
decomposition
\begin{equation} \label{decomposition}
L^{\infty}(\Del^*) = \mathcal{N}(\Del^*) \oplus
\Omega^{-1,1}(\Del^*)
\end{equation}
identifies the holomorphic tangent space $T_0
T(1)=L^{\infty}(\Del^*)/ \mathcal{N}(\Del^*)$ at the origin of
$T(1)$ with the Banach space $\Omega^{-1,1}(\Del^*)$. The
universal Teichm\"{u}ller space $T(1)$ is a complex Banach
manifold modeled on  $\Omega^{-1,1}(\Del^*)$.
\begin{remark}
Traditionally, the universal \Te space is defined using the
complex Banach space $L^{\infty}(\Del)_1$. The reflection
\eqref{sym} establishes natural complex anti-linear isomorphism
between $L^{\infty}(\Del^*)_1$ and $L^{\infty}(\Del)_1$, and the
universal \Te space in the traditional definition is complex
conjugate to the space $T(1)$ defined above.
\end{remark}
\subsubsection{The group structure} The unit ball $L^{\infty}(\Del^*)_1$
carries a group structure induced by the composition of
q.c.~mappings. The group law
\begin{equation*}
\lambda = \nu *\mu^{-1}
\end{equation*}
is defined through $w_{\lambda} = w_{\nu}\circ w_{\mu}^{-1}$,
where $\mu^{-1}$ stands for the inverse element to $\mu$, i.e.,
$\mu\ast\mu^{-1}=0$. The group law is given explicitly by
\begin{align*}
\lambda = \left(\dfrac{\nu - \mu}{1 -
\bar{\mu}\nu}\:\dfrac{(w_\mu)_z}{(\ov{w}_\mu)_{\z}} \right)\circ
w_\mu^{-1}\,.
\end{align*}
It follows from this formula that $\mathcal{L}^{\infty}(\Del^*)_1$
is a subgroup of $L^{\infty}(\Del^*)_1$.

For every $\lambda\in L^{\infty}(\Del^*)_1$ set
$[\lambda]=\Phi(\lambda)\in T(1)$. The group structure on
$L^{\infty}(\Del^*)_1$ projects to $T(1)$ by
$[\lambda]*[\mu]=[\lambda*\mu]$. For every $\mu\in
L^{\infty}(\Del^*)_1$ the right translations
\begin{equation*}
R_{[\mu]}:T(1) \rightarrow T(1),\quad [\lambda]\mapsto
[\lambda* \mu],
\end{equation*}
are biholomorphic automorphisms of  $T(1)$. The left translations,
in general, are not even continuous mappings (see, e.g.,
\cite[Sect.~III.3.4]{Lehto}). For every $\mu\in
L^{\infty}(\Del^*)_1$ the kernel of $D_\mu\Phi$ is the subspace
$D_0R_\mu\left(\mathcal{N}(\Del^*)\right)$ of  $L^\infty(\Del^*)$
and
\begin{displaymath}
T_{[\mu]}T(1) =D_0R_{[\mu]}\left(T_0T(1)\right)\simeq
D_0R_\mu\left(\Omega^{-1,1}(\Del^*)\right).
\end{displaymath}
\subsubsection{The Bers embedding}\label{Bers}
Let $A_\infty(\Del)$ be the complex Banach space
\begin{displaymath}
A_\infty(\Del)= \left\{\phi \;\text{holomorphic on}\,\, \Del: \,
\|\phi\|_\infty = \sup_{z\in \Del} \left|(1-|z|^2)^2\phi(z)\right|
< \infty \right\}.
\end{displaymath}
The Bers embedding $\beta: T(1)\hookrightarrow A_{\infty}(\Del)$
is defined as follows. Denote by $\mathcal{S}(f)$ the Schwarzian
derivative of a conformal map $f$,
\begin{displaymath}
\mathcal{S}(f)  = \frac{f_{zzz}}{f_z} -
\frac{3}{2}\left(\frac{f_{zz}}{f_z}\right)^2.
\end{displaymath}
For every $\mu\in L^\infty(\Del^*)_1$ the holomorphic function
$\mathcal{S}(w^\mu|_\Del)\in A_\infty(\Del)$ (by Kraus-Nehari
inequality it lies in the ball of radius $6$ in $A_\infty(\Del)$).
Set
\begin{equation*}
\beta([\mu])=\mathcal{S}(w^{\mu}|_{\Del}).
\end{equation*}

The Bers embedding is a holomorphic map of complex Banach
manifolds, and its differential at the origin is
\begin{align}\label{Bersmap}
D_0\beta(\mu)(z) =-\dfrac{6}{\pi} \iint\limits_{\Del^*}
\dfrac{\mu(\zeta)}{(\zeta-z)^4}d^2\zeta.
\end{align}
The complex-linear mapping $D_0\beta$ induces the isomorphism
$\Omega^{-1,1}(\Del^*)\xrightarrow{\sim}A_{\infty}(\Del)$ of the
holomorphic tangent spaces to $T(1)$ and $A_{\infty}(\Del)$ at the
origin. The mapping $\Lambda:
A_{\infty}(\Del)\rightarrow\Omega^{-1,1}(\Del^*)$, inverse to
$D_0\beta$, is given by
\begin{displaymath}
\mu(z)=\Lambda(\phi)(z)=-\tfrac{1}{2}\,(1-|z|^2)^2\phi\left(\tfrac{1}{\z}\right)
\,\tfrac{1}{\z^4}.
\end{displaymath}
According to the Ahlfors-Weill theorem, over the ball of radius $2$ in
$A_{\infty}(\Del)$ the map $\phi\mapsto [\Lambda(\phi)]$ is the
right inverse to $\beta$, $\beta\circ\Lambda = \id$.
\subsubsection{The complex structure} \label{complex1}
For every $\mu\in L^\infty(\Del^*)_1$ let $U_\mu\subset T(1)$ be
the image of the ball of radius $2$ in $A_{\infty}(\Del)$ under
the map $h_\mu^{-1}=\Phi\circ R_\mu\circ\Lambda$. The maps
$h_{\mu\nu}=h_\mu\circ h_\nu^{-1}: h_\nu(U_\mu\cap
U_\nu)\rightarrow h_\mu(U_\mu\cap U_\nu)$ are biholomorphic as
functions on the Banach space $A_{\infty}(\Del)$. The structure of
$T(1)$ as a complex Banach manifold modeled on the Banach space
$A_{\infty}(\Del)$ is explicitly described by the complex-analytic
atlas given by the open covering
\begin{displaymath}
T(1)=\bigcup_{\mu\in L^\infty(\Del^*)_1}  U_\mu
\end{displaymath}
with coordinate maps $h_\mu$ and transition maps
$h_{\mu\nu}$. The canonical projection $\Phi:
L^\infty(\Del^*)_1\rightarrow T(1)$ is a holomorphic
submersion and the Bers embedding $\beta: T(1)\rightarrow
A_{\infty}(\Del)$ is a biholomorphic map with respect to this
complex structure.

\begin{remark}
Since every point $T(1)$ can be represented by a real-analytic
Beltrami differential, it is sufficient to consider the atlas
formed by the charts $(U_{\mu},h_{\mu})$ with $\mu\in
\mathcal{L}^{\infty}(\Del^*)_1$.
\end{remark}

Complex coordinates on $T(1)$ defined by the coordinate charts
$(U_\mu, h_\mu)$ are called Bers coordinates. For every
$\nu\in\Omega^{-1,1}(\Del^*)$ set $\phi=D_0\beta(\nu)$ and define
a holomorphic vector field $\tfrac{\del}{\del \vep_\nu}$ on $U_0$
by setting
\begin{displaymath}
Dh_0\left(\tfrac{\del}{\del \vep_\nu}\right)=\phi
\end{displaymath}
at all points in $U_0$ \footnote{We identify holomorphic
tangent space to $A_{\infty}(\Del)$ at every point with
$A_{\infty}(\Del)$.}. At every point $[\mu]\in U_0$, identified
with the corresponding harmonic Beltrami differential $\mu$, the
vector field $\tfrac{\del}{\del \vep_\nu}$  in terms of the Bers
coordinates on $U_\mu$ corresponds to
\begin{align*}
\tilde{\phi} = D_\mu h_\mu\left(\tfrac{\del}{\del \vep_\nu}\right)
=\left(D_\mu h_\mu \left(D_\mu h_0\right)^{-1}\right)(\phi) =
D_0\left(\beta\circ\Phi\right)\left( D_\mu R^{-1}_\mu
\left(\Lambda(\phi)\right)\right).
\end{align*}
Using identification $\Omega^{-1,1}(\Del^*)\simeq
A_{\infty}(\Del)$, provided by the mapping $D_0\beta$, we get
\begin{equation} \label{bers-change}
\left.\dfrac{\del}{\del \vep_\nu}\right\vert_\mu = P\left(D_\mu
R_\mu^{-1}(\nu)\right) = P\left(R(\nu,\mu)\right),
\end{equation}
where
\begin{equation} \label{R}
R(\nu, \mu)=\left(\dfrac{\nu}{1 -
|\mu|^2}\:\dfrac{(w_\mu)_z}{(\ov{w}_\mu)_{\z}} \right)\circ
w_\mu^{-1},
\end{equation}
and $P: L^\infty(\Del^*)\rightarrow\Omega^{-1,1}(\Del^*)$ is the
projection onto the subspace of harmonic Beltrami differentials,
defined by the decomposition \eqref{decomposition}. Explicitly,
\begin{equation}\label{projection}
(P\mu)(z) = \frac{3(1 -
|z|^2)^2}{\pi}\iint\limits_{\Del^*}\frac{\mu(\zeta)}
{(1-\zeta\z)^4}d^2\zeta.
\end{equation}
\begin{remark}
Right translating $\nu\in T_0 T(1)$ defines a holomorphic tangent
vector
\begin{displaymath}
D_0 R_{[\mu]}(\nu)=(1-|\mu|^2)\,\nu\circ
w_\mu\,\dfrac{(\ov{w}_\mu)_{\z}}{(w_\mu)_z}\in T_{[\mu]}T(1)
\end{displaymath}
at every $[\mu]\in T(1)$. In Bers coordinates on $U_\mu$ this
tangent vector is represented by $\nu\in\Omega^{-1,1}(\Del^*)$. However, the family $\{D_0
R_{[\mu]}(\nu)\}_{[\mu]\in T(1)}$ of holomorphic tangent vectors
does not form a smooth vector field on $T(1)$ since the left
translations are not continuous on $T(1)$.
\end{remark}
\subsubsection{The universal \Te curve} \label{curve}
The universal Teichm\"{u}ller curve $\mathcal{T}(1)$ is a natural
complex fiber space over $T(1)$ with a holomorphic projection map
$\pi: \mathcal{T}(1)\rightarrow T(1)$. The fiber over each point
$[\mu]$ is a quasi-disk $w^{\mu}(\Del^*) \subset\hat{\C}$ with
complex structure induced from $\hat{\C}$ and
\begin{align}\label{universalcurve}
\mathcal{T}(1) = \left\{ ([\mu], z) : [\mu] \in T(1),\;z \in
w^{\mu}(\Del^*)\right\}.
\end{align}
The fibration $\pi: \mathcal{T}(1)\rightarrow T(1)$ has a natural
holomorphic section given by $T(1)\ni
[\mu]\mapsto([\mu],\infty)\in\mathcal{T}(1)$ --- the ``zero
section'', which defines the embedding
$T(1)\hookrightarrow\mathcal{T}(1)$. The universal \Te curve is a
complex Banach manifold modeled on $A_{\infty}(\Del)\oplus \C$
\footnote{Here $\hat{\C} \setminus \{0\}$ is identified with $\C$
via the conformal map $z \mapsto 1/z$.}, and the mapping
\begin{displaymath}
T(1) \times \Del^*\ni([\mu], z ) \mapsto ([\mu],
w^{\mu}(z))\in\mathcal{T}(1)
\end{displaymath}
is a real-analytic isomorphism.
\subsection{Homogeneous spaces of
$\text{Homeo}_{qs}(S^1)$}\label{homogenuous}
 Let
$\text{Homeo}_{qs}(S^1)$ be the group of orientation preserving
quasi-symmetric homeomorphisms of the unit circle $S^1$ (see,
e.g., \cite{Lehto} for the definition), and let $\Diff_+(S^1)$,
$\Mob(S^1)$, and $S^1$ be the subgroups of
$\text{Homeo}_{qs}(S^1)$ consisting, respectively, of smooth
orientation preserving diffeomorphisms of $S^1$, of $\Mobius$
transformations of $S^1$, and of rotations of $S^1$.

Denote by $\mathcal{U}$ the set of univalent functions on $\Del$
and let
\begin{align*}
\mathcal{D} &= \left\{ f\in \mathcal{U}: \, f(0)=0, f'(0)=1,
f''(0) = 0,\, f \;\text{admits a q.c. extension to}\; \C\right\},
\\ \widetilde{\mathcal{D}} &  = \left\{ f\in \mathcal{U}: \,
f(0)=0, f'(0)=1,\, f \;\text{admits a q.c. extension to}\;
\C\right\}.
\end{align*}
According to the Beurling-Ahlfors extension theorem, the maps
\begin{displaymath}
T(1)\ni [\mu]\mapsto w^\mu|_\Del\in\mathcal{D}
\end{displaymath}
and
\begin{displaymath}
T(1)\ni [\mu] \mapsto w_{\mu}\vert_{S^1}\in
\Mob(S^1)\backslash\text{Homeo}_{qs}(S^1)
\end{displaymath}
define bijections
\begin{equation} \label{bijections}
\mathcal{D}\xleftarrow{\sim} T(1) \xrightarrow{\sim}
\Mob(S^1)\backslash\text{Homeo}_{qs}(S^1),
\end{equation}
which endow the spaces $\mathcal{D}$ and
$\Mob(S^1)\backslash\text{Homeo}_{qs}(S^1)$ with the structure of
complex Banach manifolds modeled on the Banach space
$A_\infty(\Del)$. In what follows, we will always identify the
coset space $\Mob(S^1)\backslash\text{Homeo}_{qs}(S^1)$ with the
subgroup of $\text{Homeo}_{qs}(S^1)$ fixing $-1,-i$ and $1$, so
that the bijection
$T(1)\xrightarrow{\sim}\Mob(S^1)\backslash\text{Homeo}_{qs}(S^1)$
is a group isomorphism.
\begin{remark} \label{nontrivial}
It is a non-trivial problem to describe the complex Banach
manifold structure of the spaces $\mathcal{D}$ and
$\Mob(S^1)\backslash\text{Homeo}_{qs}(S^1)$ intrinsically, without
using the bijection \eqref{bijections}.
\end{remark}
\subsubsection{Conformal welding}\label{conformalwelding}
According to Beurling-Ahlfors extension theorem, for every $\g\in
\Mob(S^1)\backslash\text{Homeo}_{qs}(S^1)$ there exists a unique
$\alpha\in\Mob(S^1)$ which fixes $1$, and univalent functions $f$
and $g$ on $\Del$ and $\Del^*$, satisfying the following
properties.
\begin{itemize}
\item[\textbf{CW1.}] $f$ and $g$ admit q.c. extensions to $\C$.
\item[\textbf{CW2.}] $\alpha\circ\g =(g^{-1}\circ f)|_{S^1}$.
\item[\textbf{CW3.}] $f(0)=0,\,f'(0)=1,\,f''(0)=0$.
\item[\textbf{CW4.}] $g(\infty)=\infty$.
\end{itemize}
The factorization \textbf{CW2} is known as conformal welding. For
$\g=w_{\mu}|_{S^1}$, $[\mu]\in T(1)$, $f=w^\mu|_\Del$ and
$g=(w^\mu\circ w_\mu^{-1}\circ \alpha^{-1})|_{\Del^*}$, so that
$g(\Del^*)=w^\mu(\Del^*)$. Here $w^{\mu}$ is normalized so that
$f$ satisfies \textbf{CW3} and $\alpha\in\Mob(S^1)$ is uniquely
determined by the conditions $\alpha(1)=1$ and \textbf{CW4}. For
$[\mu]\in T(1)$ we will always denote $\g_\mu=(\alpha\circ
w_\mu)|_{S^1}$, $f^\mu=f$ and $g_\mu=g$, so that
\begin{displaymath}
\g_{\mu}=(g^{-1}_\mu\circ f^\mu)|_{S^1}.
\end{displaymath}
Slightly abusing notations, we will denote by $\g_\mu$ a
q.c.~extension of $\g_\mu=(\alpha\circ w_\mu)|_{S^1}\in
\Mob(S^1)\backslash\text{Homeo}_{qs}(S^1)$ given by $\alpha\circ
w_\mu$. Since $\alpha\in\Mob(S^1)$ fixes $1$, the q.c.~mapping
$\g_\mu$ satisfies the reflection property \eqref{reflection} and
the factorization
\begin{equation} \label{factorization}
\g_\mu = g^{-1}_\mu\circ f^\mu,
\end{equation}
where $f^\mu=w^\mu$ and $g_\mu = w^\mu\circ w_\mu^{-1}\circ
\alpha^{-1}$. We will distinguish between $\g_\mu\in
\Mob(S^1)\backslash\text{Homeo}_{qs}(S^1)$ and its q.c.~extension
by explicitly specifying either the property \textbf{CW2} or the
factorization \eqref{factorization}.

The following result will be used in Section 3.
\begin{lemma} \label{g(0)=0} Let $\mu\in\Omega^{-1,1} (\Del^*)_1=\Omega^{-1,1}(\Del^*)\cap
L^{\infty}(\Del^*)_1$  and $\g_\mu = \alpha\circ w_\mu$  the
q.c.~mapping introduced above. Then the mapping $\g_\mu$ fixes $0$
and $\infty$.
\end{lemma}
\begin{proof} By the
reflection property \eqref{reflection} and the factorization
\eqref{factorization}, it is sufficient to prove that
$f^{\mu}=w^\mu$ fixes $\infty$. Denote
\begin{align*}
\gamma = \imath\circ \gamma_{\mu}\circ \imath, \hspace{0.5cm} g =
\imath\circ g_{\mu} \circ \imath, \hspace{0.5cm} f = \imath\circ
f^{\mu} \circ \imath,\hspace{0.5cm} \imath^{\ast}(\mu)= \mu\circ
\imath\,\frac{\ov{\imath_z}}{\imath_z},
\end{align*}
where $\imath(z)=z^{-1}$. The factorization \eqref{factorization}
for $\g_\mu$ gives $\g = g^{-1} \circ f$, and the property
\textbf{CW3} for $f^\mu$ yields the following Laurent expansion of
$f$ at $\infty$,
\begin{align}\label{Laurent}
f(z) = z + \frac{a_1}{z} + \frac{a_2}{z^2} + \cdots.
\end{align}
We will prove that $f(0)=0$ for $\mu\in\Omega^{-1,1}(\Del^\ast)_1$
by exploiting the argument in Royden-Earle's proof of the
Ahlfors-Weill theorem, as presented in \cite[Sect.~3.8.5]{Nag2}.

Namely, $f$ satisfies the Beltrami equation with the Beltrami
differential $\nu=\imath^{\ast}(\mu)|_{\Del}$, which is supported
on $\Del$. The fundamental theorem from the theory of
q.c.~mappings (see, e.g. \cite{Ahl2}) asserts that $f$ admits the
series representation
\begin{align}\label{series1}
f(z) = z + P(\nu)(z) + P(\nu H(\nu))(z)+ P(\nu H(\nu H(\nu)))(z)+
\cdots,
\end{align}
which is uniformly and absolutely convergent on $\CC$. Here for
$h\in C^2(\Del)$ we denote
\begin{align*}
P(h)(z) &= -\frac{1}{\pi}\iint\limits_{\Del}
\frac{h(\zeta)}{\zeta-z} \,d^2 \zeta,\\ H(h) (z) &=
-\frac{1}{\pi}\iint\limits_{\Del}\frac{h(\zeta)}{(\zeta-z)^2}\,
d^2\zeta,
\end{align*}
where the latter integral --- the Hilbert transform, is understood
in the principal value sense. Since $\nu$ has compact support, it
immediately follows from the definition of the operators $P$ and
$H$ that the series \eqref{series1} has the Laurent expansion
\eqref{Laurent} at $\infty$. We will prove that for $\nu\in
\Omega^{-1,1}(\Del)$ each term of this series vanishes at $z=0$.
Representing $\nu(z) = -\frac{1}{2} (1-|z|^2)^2
\sum_{n=0}^{\infty} a_n \z^{n}$ and using polar coordinates, we get
for any $(n-1)$ -- iterate of the operator $\nu H$, $n>1$,
\begin{align*}
&P(\nu H(\nu H(\nu \ldots H(\nu))))(0) \\
 = &\left(\frac{1}{2\pi}\right)^n \iint\limits_{\Del} \dots
\iint\limits_{\Del} \frac{ \sum_{m_1,\ldots m_n} a_{m_1} \ldots
a_{m_n} r_1^{m_1+1} \ldots r_n^{m_n+1} e^{-i m_1 \theta_1} \ldots
e^{-i m_n \theta_n}} { r_1 e^{i\theta_1} (r_2e^{i\theta_2} -r_1
e^{i\theta_1})^2 \ldots (r_n e^{i\theta_n} - r_{n-1} e^{i
\theta_{n-1}})^2} \\ & \hspace{4cm}(1-r_1^2)^2\ldots (1-r_n^2)^2
dr_1 d\theta_1 \ldots dr_n d\theta_n\\ =& \sum_{m_1, \ldots,
m_n=0}^\infty a_{m_1}\ldots a_{m_n} I_{m_1, \ldots, m_n},
\end{align*}
where each integral in the definition of $H$ is understood in the
principal value sense. The interchange of the orders of summation
and integration can be easily justified. For fixed $r_1\neq 0,
r_1\neq r_2, r_2\neq r_3,\ldots, r_{n-1}\neq r_n$, let
\begin{align*}
& I_{m_1, \ldots, m_n}(r_1,\ldots,r_n)\\ & = \int_0^{2\pi} \dots
\int_0^{2\pi} \frac{e^{-i m_1 \theta_1} \ldots e^{-i m_n
\theta_n}d\theta_1\cdots d\theta_n} { r_1 e^{i\theta_1}
(r_2e^{i\theta_2} -r_1 e^{i\theta_1})^2 \ldots (r_n e^{i\theta_n}
- r_{n-1} e^{i \theta_{n-1}})^2}.
\end{align*}
A change of variables $\theta_k \mapsto \theta_k + \theta$,
$k=1,\dots,n$ gives
\begin{align*}
I_{m_1,\dots,m_n}(r_1,\ldots,r_n) = e^{-i (m_1+\ldots + m_n +
(2n-1))\theta} I_{m_1, \ldots, m_n}(r_1,\dots,r_n).
\end{align*}
Since all $m_k\geq 0$ and $2n-1 > 0$ for $n\geq 1$, we have $e^{-i
(m_1+\ldots + m_n + (2n-1))\theta} \neq 1$ and hence
\begin{align*}
I_{m_1, \ldots, m_n}(r_1,\ldots,r_n) =0.
\end{align*}
This proves that all $I_{m_1, \ldots, m_n}$ vanish. Therefore
$f(0)=0$.
\end{proof}
\begin{remark} Since $P(f)_z = H(f)$, it also follows from the proof
that $f_z(0)=1$.
\end{remark}
Similar to \eqref{bijections} , there are bijections
\begin{equation*}
\widetilde{\mathcal{D}}\xleftarrow{\sim} \mathcal{T}(1)
\xrightarrow{\sim} S^1\backslash\text{Homeo}_{qs}(S^1),
\end{equation*}
where we always identify the coset space
$S^1\backslash\text{Homeo}_{qs}(S^1)$ with the stabilizer of $1$
in $\text{Homeo}_{qs}(S^1)$ (see, e.g., \cite{Teo}). For every $\g
\in S^1\backslash\text{Homeo}_{qs}(S^1)$ there exist unique
univalent functions $f$ and $g$ on $\Del$ and $\Del^*$, satisfying
the properties \textbf{CW1, CW4} and
\begin{enumerate}
\item[\textbf{CW2$^\prime$.}] $\g= (g^{-1} \circ f)\vert_{S^1}$;
\item[\textbf{CW3$^\prime$.}] $f(0)=0,\,f'(0)=1$.
\end{enumerate}
Namely, the fibration $\pi: \mathcal{T}(1)\longrightarrow T(1)$
corresponds to the fiber space
$S^1\backslash\text{Homeo}_{qs}(S^1)$ over
$\Mob(S^1)\backslash\text{Homeo}_{qs}(S^1)$ with the fibers
isomorphic to $S^1\backslash\Mob(S^1)\simeq\Del^*$. The points in
the fiber at $[\mu]\in T(1)$ correspond to the points
$\sigma_w\circ\g_\mu\in
S^1\backslash\text{Homeo}_{qs}(S^1),\,w\in\Del^*$
with\footnote{Here the subscript $w$ does not stand for the
derivative.}
\begin{align*}
\si_w(z) =\frac{1-w}{1-\bar{w}}\,\frac{1 - z \bar{w}}{z-w}\in
S^1\backslash\Mob(S^1).
\end{align*}
Using the properties \textbf{CW1} and \textbf{CW2} for $\g_\mu$,
we get the factorization \textbf{CW2$^\prime$} for $\gamma$ is
\begin{equation*}
\g=\sigma_w\circ\g_\mu = (g^{-1}\circ f)|_{S^1},
\end{equation*}
where
\begin{align*}
f = \lambda_w \circ f^\mu,\quad g = \lambda_w \circ g_\mu \circ
\si_w^{-1},
\end{align*}
and
\begin{align*}
\lambda_w (z)= \frac{z}{c_w z+1},\quad c_w
=-\frac{1}{g_\mu(w)}=-\frac{1}{2}\frac{f^{\prime\prime}(0)}{f^{\prime}(0)},
\end{align*}
so that $(g_\mu\circ\si^{-1}_w)(\infty)=g_\mu(w)$, and the
functions $f$ and $g$ satisfy the properties \textbf{CW3$^\prime$}
and \textbf{CW4} respectively. The mapping
\begin{displaymath}
\mathcal{T}(1)\ni([\mu],g_\mu(w))\mapsto \g=\sigma_w\circ\g_\mu\in
S^1\backslash\text{Homeo}_{qs}(S^1)
\end{displaymath}
establishes the isomorphism $\mathcal{T}(1)\xrightarrow{\sim}
S^1\backslash\text{Homeo}_{qs}(S^1)$.

As before, we will also denote by $\g$ a q.c.~extension of
$\g\in S^1\backslash\text{Homeo}_{qs}(S^1)$ which satisfies the
reflection property \eqref{reflection} and admits the factorization
$\g=g^{-1}\circ f$.

\begin{remark}
It is known \cite{Ki} that $\Diff_+(S^1)$ is an
infinite-dimensional Lie group and homogeneous spaces
$\Mob(S^1)\backslash\Diff_+(S^1)$ and $S^1\backslash \Diff_+(S^1)$
are infinite-dimensional complex Fr\'{e}chet manifolds. In this
case conformal welding readily follows from the Riemann mapping
theorem without using q.c. mappings~\cite{Ki}. Note that our
convention for the conformal welding is different from that in
\cite{Ki}: we are using right cosets instead of left cosets.
\end{remark}

The bijection $\mathcal{T}(1)\xrightarrow{\sim} S^1
\backslash\text{Homeo}_{qs}(S^1)$ endows the universal
Teichm\"{u}ller curve $\mathcal{T}(1)$ with the group structure.
Explicitly,
\begin{align*}
([\lambda], z)= ([\nu],\zeta) * ([\mu],w)^{-1},
\end{align*}
where
\begin{equation} \label{gp2}
\lambda = \left(\dfrac{\nu - \mu}{1 -
\bar{\mu}\nu}\:\dfrac{\g_z}{\ov{\g}_{\z}} \right)\circ \g^{-1}
\end{equation}
and
\begin{equation} \label{gp3}
z = \left(w^{\lambda} \circ \g \circ (w^\nu)^{-1}\right)(\zeta).
\end{equation} Here $\g$ is a q.c.~extension of $\sigma_u\circ
\g_\mu$, $u=g_{\mu}^{-1}(w)$, and the point
$([\lambda],z)\in\mathcal{T}(1)$ does not depend on the choice of
the extension $\g$.

The mapping $$ T(1)\ni [\mu]\mapsto \gamma_{\mu}\in S^1\bk
\text{Homeo}_{qs}(S^1) \simeq \mathcal{T}(1)$$ is a
complex--analytic embedding of $T(1)$ into $\mathcal{T}(1)$ which
is not a group homomorphism. On the other hand, the natural
embedding $$\Mob(S^1)\bk \text{Homeo}_{qs}(S^1) \simeq T(1) \ni
[\mu] \mapsto w_{\mu}\in S^1\bk\text{Homeo}_{qs}(S^1)\simeq
\mathcal{T}(1)$$ is a group homomorphism, though not a
complex--analytic mapping. Considering $w_{\mu}$ as an element of
$S^1\bk\text{Homeo}_{qs}(S^1)$ via this embedding, it admits a
conformal welding
\begin{align}\label{newcon}
w_{\mu}=\rg_{\mu}^{-1}\circ \rf^{\mu}
\end{align}
 that satisfies the
properties \textbf{CW1, CW2$^{\prime}$, CW3$^{\prime}$, CW4}. In
terms of the functions $f^{\mu}$ and $g_{\mu}$ satisfying
properties \textbf{CW1--CW4}, we have
$$ \rf^{\mu}=\lambda_{\mu}\circ
f^{\mu}\quad\text{and}\quad\rg_{\mu}=\lambda_{\mu}\circ
g_{\mu}\circ \alpha_{\mu}^{-1},$$ where $\alpha_{\mu}\in
\PSU(1,1)$ is such that $w_{\mu}=\alpha_{\mu}\circ \gamma_{\mu}$,
and $\lambda_{\mu}\in \PSL(2,\C)$ is uniquely determined by the
conditions $\rf^{\mu}(0)=0, (\rf^{\mu})^{\prime}(0)=1$ and
$\rg_{\mu}(\infty)=\infty$.

\subsubsection{The horizontal and vertical subspaces}
\label{h-v-subspaces} The right translations $R_{([\mu], z)}:
\mathcal{T}(1) \rightarrow \mathcal{T}(1)$ are biholomorphic
automorphisms of $\mathcal{T}(1)$ \cite{Bers1}. The holomorphic
tangent space to $\mathcal{T}(1)$ at $([\mu], z)$ is identified
with the holomorphic tangent space at $(0,\infty)$ --- the origin
of $\mathcal{T}(1)$ by
\begin{displaymath}
T_{([\mu],z)}\mathcal{T}(1) = D_{(0,\infty)}
R_{([\mu],z)}(T_{(0,\infty)}\mathcal{T}(1)) \simeq
T_{(0,\infty)}\mathcal{T}(1).
\end{displaymath}
The holomorphic tangent space at the origin naturally splits into
the direct sum of horizontal and vertical subspaces,
\begin{displaymath}
T_{(0,\infty)}\mathcal{T}(1) = \Omega^{-1,1}(\Del^*)\oplus\C.
\end{displaymath}
The identification of holomorphic tangent spaces provides a
natural splitting of the tangent space at every point in
$\mathcal{T}(1)$ into the direct sum of horizontal and vertical
subspaces. Lifts of horizontal and vertical tangent vectors at
the origin of $\mathcal{T}(1)$ to every point in the fiber at the
origin are explicitly described as follows.
\begin{enumerate}
\item[\textbf{TV1.}] Let $\mu\in \Omega^{-1,1}(\Del^*)
\subset T_{(0,\infty)}\mathcal{T}(1)$ be a horizontal tangent
vector to $\mathcal{T}(1)$ at the origin. A curve $([t\mu],
z(t))$, $z(0)=z$, which defines the horizontal lift of $\mu$ to
the point $(0,z)\in \mathcal{T}(1)$ in the fiber $\pi^{-1}(0)$ at
the origin, for small $t$ is given by the equation
\begin{equation*}
([\mu(t)],\infty) * (0,z) = ([t\mu], z(t)).
\end{equation*}
Using \eqref{gp2}, \eqref{gp3} and Lemma \ref{g(0)=0}, we get
\begin{equation*} \label{vectoridentification}
\mu(t)=(\si_z^{-1})^*(t\mu)=t\mu\circ\si_z^{-1}
\frac{\ov{(\si_z^{-1})'}}{(\si_z^{-1})'}\quad\text{and}\quad z(t)=w^{t\mu}(z).
\end{equation*}
Thus the horizontal lift of $\mu\in T_{(0,
\infty)}\mathcal{T}(1)$ to every point in the fiber $(0,z)\in\pi^{-1}(0)$ is
the vector field
\begin{equation*}
\tau_\mu =\left.\tfrac{\pa}{\pa\vep_\mu}\right|_0 +\dot{w}^\mu(z)\tfrac{\pa}{\pa z},
\quad\text{where}\quad \dot{w}^\mu(z) =\tfrac{dz}{dt}(0)
\end{equation*}
(cf.~\cite{Wol}). At  the point $(0,z)\in\pi^{-1}(0)$ the vector
field $\tau_\mu$ is identified with the horizontal tangent vector
$(\si_z^{-1})^*\mu\in T_{(0,\infty)}\mathcal{T}(1)$.
\item[\textbf{TV2.}] Let $\mathbf{1}\in\C\subset
T_{(0,\infty)}\mathcal{T}(1)$ be the vertical tangent vector to
$\mathcal{T}(1)$ at the origin, given by the value of the vector
field $\pa_z = \tfrac{\pa}{\pa z}$ at $z=\infty$. A curve defining
the right translate of $\mathbf{1}$ to the point
$(0,z)\in\mathcal{T}(1)$ in the fiber $\pi^{-1}(0)$ at the origin
for small $t$ is given by the equation
\begin{displaymath}
(0, t^{-1}) * (0,z) =  (0,z(t)),
\end{displaymath}
and it follows from \eqref{gp3} that
\begin{displaymath}
\frac{dz}{dt}(0) = \frac{(1-z)(1-|z|^2)}{(1-\z)}.
\end{displaymath}
Thus the right translate of $\mathbf{1}\in T_{(0,
\infty)}\mathcal{T}(1)$ to the point $(0,z)\in\pi^{-1}(0)$ is the
vector  $\tfrac{(1-z)(1-|z|^2)}{(1-\z)}\pa_z$ at $(0,z)$. As a
result, the vector field $\pa_z$ at every point
$(0,z)\in\pi^{-1}(0)$ is identified with the vertical tangent
vector
\begin{displaymath}
\frac{(1-\z)}{(1-z)(1-|z|^2)}\mathbf{1}\in
T_{(0,\infty)}\mathcal{T}(1).
\end{displaymath}
\end{enumerate}
\subsubsection{The isomorphisms of the tangent spaces}\label{TS}
The real tangent vector space
$T^\RR_0\,S^1\backslash\text{Homeo}_{qs}(S^1)$ to
$S^1\backslash\text{Homeo}_{qs}(S^1)$ at the origin is identified
with the subspace of Zygmund class continuous real-valued vector
fields $\text{u}=u(\theta)\tfrac{d}{d\theta}$ on $S^1$ (see, e.g.,
\cite{Teo} for the definition), satisfying
\begin{displaymath}
\int_0^{2\pi}u(\theta)d\theta = 0.
\end{displaymath}
In particular, the Fourier series $u(\theta)=\sum_{n\in\ZZ} c_n
e^{in\theta}$ is absolutely convergent. For $|z|=1$ set
\begin{displaymath}
\tilde{u}(z) = i\sum\limits_{n\in\ZZ\setminus\{0\}} c_n
z^{n+1}.
\end{displaymath}
The function $\tilde{u}$ on $S^1$ admits the decomposition
\begin{displaymath}
\tilde{u} = u_+ + u_{-},
\end{displaymath}
where $u_+$ and $u_-$ are boundary values of functions holomorphic
on  $\Del$ and $\Del^*$ respectively and $u_+(0) = 0$. Explicitly,
\begin{align*}
u_+(z) =& i\sum_{n=1}^{\infty}c_n z^{n+1}, \\
u_-(z) =&i\sum_{n=1}^{\infty}c_{-n}z^{1-n}.
\end{align*}

It is a difficult problem to characterize the Zygmund class in
terms of the Fourier series (cf.~Remark \ref{nontrivial}). On the
other side, in terms of the Fourier series the almost complex
structure $J$ on $T^\RR_0\,S^1\backslash\text{Homeo}(S^1)$ is
explicitly given by the classical conjugation operator
\begin{align*}
J\,\text{u} = i\sum\limits_{n\in\ZZ\setminus\{0\}}\sgn(n) c_n
e^{in\theta} \dfrac{d}{d\theta} \quad\text{for}\quad \text{u} =
\sum\limits_{n\in\ZZ\setminus\{0\}} c_n
e^{in\theta}\dfrac{d}{d\theta}.
\end{align*}
\begin{remark}
Note that our definition of the operator $J$ differs by a negative
sign  from the definition in \cite{Ki, NV} for the homogeneous
space $S^1\backslash \Diff_+(S^1)$.
\end{remark}
The holomorphic and anti-holomorphic tangent vectors at the origin
are
\begin{align*}
\text{v}=\dfrac{\text{u} - i J\,\text{u}}{2} = \sum_{n=1}^\infty
c_n e^{in\theta}\dfrac{d}{d\theta}\quad\text{and} \quad
\bar{\text{v}}=\dfrac{\text{u} + i J\,\text{u}}{2} =
\sum_{n=-\infty}^{-1} c_n e^{in\theta}\dfrac{d}{d\theta}.
\end{align*}

For every smooth function $\mathcal{F}$ in a neighborhood of the
origin in $\mathcal{T}(1)$ and $\text{u}\in
T^\RR_0\,S^1\backslash\text{Homeo}_{qs}(S^1)$ set
\begin{displaymath}
\dot{\mathcal{F}}[\text{u}]=\left.\dfrac{d}{dt}\right\vert_{t=0}
\mathcal{F}(\g^t),
\end{displaymath}
where $\g_t$ is a curve in $S^1\backslash\text{Homeo}_{qs}(S^1)$
with the tangent vector $\text{u}$ at the origin. Corresponding
directional derivatives of $\mathcal{F}$ at the origin in
$\mathcal{T}(1)$ in the holomorphic and anti-holomorphic
directions $\text{v}$ and $\bar{\text{v}}$ are defined by
\begin{equation} \label{de-Rham}
\pa\mathcal{F}(\text{v})=\dfrac{1}{2}\,\left(\dot{\mathcal{F}}[\text{u}]
- i\,\dot{\mathcal{F}}[J\,\text{u}] \right),\quad
\text{and}\quad\bar\pa\mathcal{F}(\bar{\text{v}})=
\dfrac{1}{2}\left(\dot{\mathcal{F}}[\text{u}] +
i\,\dot{\mathcal{F}}[J\,\text{u}]\right).
\end{equation}

For $s\in\RR$ let
\begin{equation*}
\mathcal{H}^s(S^1)=\left\{\text{u}=\sum_{n=-\infty}^\infty a_n e^{in\theta}
\dfrac{d}{d\theta}: \sum_{n=-\infty}^\infty
|n|^{2s}|a_n|^2<\infty\right\}
\end{equation*}
be the Sobolev space of complex-valued vector fields on $S^1$. The
properties of the tangent spaces
$T_0S^1\backslash\text{Homeo}_{qs}(S^1)$,
$T_0\widetilde{\mathcal{D}}$ and
$T_0\Mob(S^1)\backslash\text{Homeo}_{qs}(S^1)$, which will be used
in Section 5, can be succinctly summarized as follows (see
\cite{Teo} for details).
\begin{itemize}
\item[\textbf{TS1.}] Under the $\RR$-linear isomorphism
$T^\RR_0\,S^1\backslash\text{Homeo}_{qs}(S^1)
\xrightarrow{\sim}T_0 \widetilde{\mathcal{D}}$
\begin{equation*}
u(\theta)=\sum_{n\in\ZZ \setminus\{0\}} c_n e^{in\theta}\mapsto
u_+(z)=i \sum_{n=1}^{\infty} c_n z^{n+1},
\end{equation*}
and $\dot{f}\vert_{\Del}= u_+,\,\dot{g}_0\vert_{\Del^*} = -u_-$, where
\begin{equation*}
\dot{f}= \left.\frac{d}{dt}f^t\right\vert_{t=0},\quad\dot{g} =
\left.\frac{d}{dt}g_t\right\vert_{t=0},
\end{equation*}
$\gamma_t = g_t^{-1}\circ f^t$ is a smooth curve in
$\mathcal{T}(1)$ tangent to $\text{u}$ at the origin, and
 $\dot{g}_0(z)=\dot{g}(z)-\dot{g}'(\infty)z$.
\item[\textbf{TS2.}] Under the $\RR$-linear isomorphism
$$T_0^\RR\,\Mob(S^1)\backslash\text{Homeo}_{qs}(S^1)
\xrightarrow{\sim} T_0 T(1)\xrightarrow{D_0\beta} A_\infty(\Del)$$
\begin{align*}
u(\theta) = \sum\limits_{n\in\ZZ \setminus\{-1,0,1\}} c_n e^{in
\theta}\mapsto \dfrac{d^3 u_+}{dz^3}(z)
 = i \sum_{n=2}^{\infty} (n^3 -n) c_n z^{n-2}.
\end{align*}
\item[\textbf{TS3.}] If $$\phi(z)  = \sum_{n=2}^{\infty} (n^3-n)
a_n z^{n-2} \in A_{\infty}(\Del)$$ then
\begin{displaymath}
\sum_{n=2}^\infty n^{2s} |a_n|^2 <\infty\quad\text{for all}\quad
s<1.
\end{displaymath}

\item[\textbf{TS4.}]
$T^\RR_0\,S^1\backslash\text{Homeo}_{qs}(S^1)\subset
\mathcal{H}^s(S^1)$ for all $s<1$.
\end{itemize}

\subsection{Teichm\"uller spaces and Teichm\"uller curves of Fuchsian groups}
 \label{fuchsian-groups}
Let $\Gamma$ be a Fuchsian group, i.e., a discrete subgroup of $\PSU(1,1)$. Let
\begin{align*}
L^{\infty}(\Del^*, \Gamma) = \left\{ \mu \in L^{\infty}(\Del^*) :
\mu \circ \g \frac{\ov{\g'}}{\g'} = \mu\quad\text{for all}\quad\g
\in \Gamma \right\}
\end{align*}
be the space of bounded Beltrami differentials for $\Gamma$ and
\begin{displaymath}
L^{\infty}(\Del^*, \Gamma)_1 = L^{\infty}(\Del^*)_1 \cap
L^{\infty}(\Del^*, \Gamma) \;
\end{displaymath}
be the open unit ball in $L^{\infty}(\Del^*, \Gamma)$. The \Te
space of the Fuchsian group $\Gamma$ is defined by
\begin{align*}
T(\Gamma) = L^{\infty}(\Del^*, \Gamma)_1 / \sim \;,
\end{align*}
where the equivalence relation is the same as the one used to
define the universal \Te space $T(1)$ in Section
\ref{definitions}. The \Te space $T(\Gamma)$ has a natural
structure of a complex Banach manifold such that the tangent space
at the origin of $T(\Gamma)$ is identified with the Banach space
$\Omega^{-1,1}(\Del^*, \Gamma)$ of bounded harmonic Beltrami
differentials for $\Gamma$,
\begin{displaymath}
\Omega^{-1,1}(\Del^*, \Gamma) = \Omega^{-1,1}(\Del^*) \cap
L^{\infty}(\Del^*, \Gamma).
\end{displaymath}

For every Fuchsian group $\Gamma$ the canonical embedding
$T(\Gamma)\hookrightarrow T(1)$ is holomorphic, so that the
universal \Te space $T(1)$ contains all the \Te spaces $T(\Gamma)$
as complex submanifolds. The universal \Te space $T(1)$ is the \Te
space for the trivial Fuchsian group $\Gamma=\{1\}$.

The inverse image of $T(\Gamma)$ under the projection map
$\mathcal{T}(1) \rightarrow T(\Gamma)$ is called the Bers fiber
space $\mathcal{BF}(\Gamma)$. The quasi-Fuchsian group
$\Gamma^{\mu} = w^{\mu}\circ \Gamma \circ (w^{\mu})^{-1}$ acts on
the fiber $w^{\mu} (\Del^*)$ at the point $[\mu] \in T(\Gamma)$.
The \Te curve $\mathcal{T}(\Gamma)$ of the Fuchsian group $\Gamma$
is a fiber space
 over $T(\Gamma)$ with the fiber
$\Gamma^\mu\bk w^\mu(\Del^*)$ at the point $[\mu] \in T(\Gamma)$.

The domain $\Del^*$ is a model of the hyperbolic plane $\HH^2$.
The hyperbolic (Poincar\'{e}) metric on $\Del^*$ --- a Hermitian
metric of constant Gaussian curvature $-1$, is
\begin{equation}
ds^2 = \rho(z)|dz|^2 =  \frac{4|dz|^2}{(1-|z|^2)^2}\,,
\end{equation}
and the hyperbolic area $2$-form is $\rho(z)\,d^2z$. The Fuchsian
group $\Gamma$ is  of finite type (cofinite) if the corresponding
Riemann surface --- the orbifold $\Gamma\bk\Del^*$, has a finite
hyperbolic area. In this case, the \Te space $T(\Gamma)$ is a
finite-dimensional complex manifold with a natural Hermitian
metric, called the Weil-Petersson metric. It is defined as
Petersson's inner product on tangent spaces
$T_{[\mu]}T(\Gamma)\simeq\Omega^{-1,1}(\Del^*,\Gamma_{\mu})$,
where $[\mu]\in T(\Gamma)$ and $\Gamma_\mu =w_\mu\circ\Gamma\circ
w_\mu^{-1}$. For $\mu,\nu\in T_0T(\Gamma)$,
\begin{equation*}
\la\mu,\nu\ra_{WP} =
\iint\limits_{\Gamma\bk\Del^*}\mu\bar{\nu}\rho(z)d^2z.
\end{equation*}
The Weil-Petersson metric on $T(\Gamma)$ is a \Ka metric.

The following result, due to Patterson \cite{patterson}, will be
used in Section \ref{finite-dim}. Here we present it in a
convenient form as in \cite{Teo}.
\begin{lemma} \label{regularization} Let $\Gamma$ be a cofinite
Fuchsian group and  $h \in L^{\infty}(\Del^\ast, \rho(z)d^2z)$ be
$\Gamma$--automorphic, i.e., $h\circ\g=h$ for all $\g\in\Gamma$.
Then
\begin{align*}
\iint\limits_{\Gamma\bk \Del^\ast} h(z) \rho(z)d^2z = \lim_{r
\rightarrow 1^+} \frac{A(\Gamma\bk \Del^*)}{A(\Del^\ast_r)}
\iint\limits_{\Del^\ast_r} h(z)\rho(z) d^2 z,
\end{align*}
where $\Del^\ast_r=\{z\in\Del^\ast : |z|\geq r\}$, $A(\Gamma\bk
\Del^*)$ is the hyperbolic area of the Riemann surface $\Gamma \bk
\Del^*$, and $A(\Del^\ast_r)$ is the hyperbolic area of
$\Del^\ast_r$.
\end{lemma}
\subsection{Resolvent kernel} \label{resolvent}
Let
\begin{equation} \label{laplace}
\Delta_0=-\rho(z)^{-1}\partial_z\partial_{\bar{z}}
\end{equation}
be the Laplace-Beltrami operator of the hyperbolic metric on
$\Del$, acting on functions. It is well-known (see, e.g.,
\cite{Hejhal, Lang}) that the differential expression
\eqref{laplace} defines a unique positive, self-adjoint operator
on the Hilbert space $L^2(\Del, \rho(z) d^2z)$, which we still
denote by $\Delta_0$. Let
\begin{displaymath}
G = \tfrac{1}{2}\left(\Delta_0 + \tfrac{1}{2}\right)^{-1}
\end{displaymath}
be (a one-half of) the resolvent of $\Delta_0$ at the regular point
$\lambda =-\tfrac{1}{2}$.
\begin{remark} Note that the Laplace-Beltrami operator in
\cite{Hejhal, Lang} is $4\Delta_0$, so that the regular point
$\lambda=-\tfrac{1}{2}$ for the operator $\Delta_0$ corresponds to
$\lambda = -2$ for the Laplace-Beltrami operator in
\cite{Hejhal,Lang}.
\end{remark}

The resolvent $G$ is a bounded integral operator on
$L^2(\Del,\rho(z) d^2z)$ with kernel
\begin{equation} \label{resolvent-kernel}
G(z,w) = \frac{2u+1}{2\pi} \log \frac{u+1}{u} - \frac{1}{\pi},
\end{equation}
where $u(z,w)$ is a point-pair invariant on $\Del$,
\begin{displaymath}
u(z,w) = \frac{|z-w|^2}{(1-|z|^2)(1-|w|^2)}.
\end{displaymath}

The resolvent kernel $G(z,w)$ has the following properties (see,
e.g., \cite{Hejhal} and \cite[Sect.~XIV.3]{Lang}).
\begin{itemize}
\item[\textbf{RK1.}] $G$ is symmetric, $G(z,w)=G(w,z)$, and is a
point-pair invariant,
\begin{align*}
G(\gamma z, \gamma w) = G(z,w)~\text{for all}~\gamma \in
\PSU(1,1).
\end{align*}
\item[\textbf{RK2.}] $G(z,w)$ is positive for all $z,w\in\Del$.
\item[\textbf{RK3.}] If $g\in BC^\infty(\Del)$ --- the space of smooth
bounded functions on $\Del$, then the integral
\begin{equation*}
f(z) = \iint\limits_\Del G(z,w)g(w)\rho(w) d^2w
\end{equation*}
is absolutely convergent for all $z\in\Del$ and $f=G(g)\in
BC^\infty(\Del)$ satisfies the differential equation
\begin{displaymath}
2\left(\Delta_0+\tfrac{1}{2}\right)(f) =g.
\end{displaymath}
Conversely, if $f\in BC^\infty(\Del)$ and $g =
2\left(\Delta_0+\tfrac{1}{2}\right)(f)\in BC^\infty(\Del)$, then
$f=G(g)$.
\item[\textbf{RK4.}] For all $z\in\Del$,
\begin{equation*}
\iint\limits_\Del G(z,w)\rho(w)d^2 w=1.
\end{equation*}
\end{itemize}
The last property immediately follows from \textbf{RK3} since
\begin{displaymath}
2\left(\Delta_0+\tfrac{1}{2}\right)(1)=1,
\end{displaymath}
where $1$ is the constant function equal to $1$ on $\Del$.

The resolvent kernel $G$ of the Laplace-Beltrami operator on
$\Del^*$ is given by the same formula \eqref{resolvent-kernel} and
satisfies the properties \textbf{RK1} -- \textbf{RK4}.

 When $\Gamma$ is a cofinite Fuchsian group, we denote by $G_{\Gamma}$
 the one-half of the resolvent of the Laplace-Beltrami operator on
 the Riemann surface $\Gamma\bk\Del$ at $\lambda = -\tfrac{1}{2}$.
  It is a bounded integral operator on
$L^2(\Gamma\bk\Del,\rho(z)d^2z)$ with kernel
\begin{equation} \label{G-Gamma}
G_{\Gamma}(z,w)=\sum_{\gamma\in \Gamma}G(z,\gamma w),\,z,w\in\Del,
\end{equation}
and it enjoys all the properties \textbf{RK1}-\textbf{RK4}. The
corresponding resolvent kernel on $\Gamma\bk\Del^\ast$ is given by
the same formula with $z,w\in\Del^\ast$.
\begin{remark} The operator $G_\Gamma$ plays a prominent role in
the Weil-Petersson geometry of the finite-dimensional \Te space
$T(\Gamma)$ \cite{Wol}.
\end{remark}
\subsection{Variational formulas}\label{VF}
Here we collect necessary variational formulas. To simplify the
computations in the following sections, we will use different
realizations of the hyperbolic plane $\HH^2$, given either by the
unit disk $\Del$ or its exterior $\Del^\ast$, or by the upper
half-plane $\UU$.

Let $l$ and $m$ be integers and $\Gamma$ a Fuchsian group (we will
be primarily interested in the cases when $\Gamma=\{1\}$, i.e.,
is a trivial group, and when $\Gamma$ is a cofinite Fuchsian
group). Using the model $\HH^2\simeq\Del$, tensor of type $(l,m)$
for $\Gamma$ is a $C^\infty$-function $\omega$ on $\Del$
satisfying
\begin{equation*}
\omega(\gamma z)\gamma'(z)^l
\ov{\gamma^\prime(z)}^m=\gamma(z)~\text{for all}~\gamma\in\Gamma.
\end{equation*}
Let $\omega^{\vep}$ be a smooth family of tensors of type $(l,m)$
for $\Gamma_{\vep\mu}=w_{\vep\mu}\circ \Gamma\circ
w_{\vep\mu}^{-1}$, where $\mu \in \Omega^{-1,1} (\Del,\Gamma)$ and
$\vep \in \C$ is sufficiently small. Set
\begin{equation*}
(w_{\vep\mu})^* (\omega^{\vep}) = \omega^{\vep} \circ w_{\vep\mu}
((w_{\vep\mu})_z)^l ((\ov{w}_{\vep\mu})_{\z})^m,
\end{equation*}
which is a tensor of type $(l,m)$ for $\Gamma$ --- a pull-back of
the tensor $\omega^\vep$ by $w_{\vep\mu}$. Lie derivatives of
the family $\omega^{\vep}$ along vector fields $\pa/\pa
\vep_{\mu}$ and $\pa/\pa \bar{\vep}_{\mu}$ are defined in the
standard way,
\begin{align*}
L_{\mu}\omega = \frac{\pa}{\pa \vep}
\Bigr|_{\vep=0} (w_{\vep\mu})^*
(\omega^{\vep})\quad\text{and}\quad L_{\bar{\mu}}\omega  = \frac{\pa}{\pa \bar{\vep}}
\Bigr|_{\vep=0}(w_{\vep\mu})^* (\omega^{\vep}).
\end{align*}
When $\omega$ is a function on $T(\Gamma)$ --- a tensor of type
$(0,0)$, the Lie derivatives reduce to directional derivatives
\begin{equation*}
L_{\mu}\omega =
\pa\omega(\mu)\quad\text{and}\quad L_{\bar{\mu}}\omega=
\bar\pa\omega(\bar\mu)
\end{equation*}
--- the evaluation of the 1-forms $\pa\omega$ and $\bar\pa\omega$ on the
holomorphic and antiholomorphic tangent vectors $\mu$ and
$\bar\mu$ to $T(\Gamma)$ at the origin. Corresponding real vector
fields $\frac{\pa}{\pa t_{\mu}}$ are defined by
\begin{align*}
\frac{\pa}{\pa t_{\mu}}=\frac{\pa}{\pa \vep_{\mu}}+\frac{\pa}{\pa
\bar{\vep}_{\mu}},
\end{align*}
so that
\begin{align*}
\frac{\pa}{\pa\vep_{\mu}}=\frac{1}{2}\left(\frac{\pa}{\pa
t_{\mu}}-i\frac{\pa}{\pa t_{i\mu}}\right)\quad\text{and}\quad
\frac{\pa}{\pa\bar{\vep}_{\mu}}=\frac{1}{2}\left(\frac{\pa}{\pa
t_{\mu}}+i\frac{\pa}{\pa t_{i\mu}}\right).
\end{align*}

For the model $\HH^2\simeq\UU$ we have
\begin{align}\label{variation-1}
&\frac{\pa}{\pa\vep} w_{\vep\mu}(z)= -\frac{1}{\pi}
\iint\limits_{\UU}R\left(w_{\vep\mu}(z), w_{\vep\mu}(u)\right) \mu(u) (w_{\vep\mu})^2_u(u)
 d^2u, \\
&\frac{\pa}{\pa \bar{\vep}} w_{\vep\mu}(z) = -\frac{1}{\pi}
\iint\limits_{\UU}
R\left(w_{\vep\mu}(z),\ov{w_{\vep\mu}(u)}\right)
\ov{\mu(u)(w_{\vep\mu})^2_{u}(u)} d^2u,\nonumber
\end{align}
where the q.c.~mapping $w_{\vep\mu}$ is normalized by fixing
$0,1,\infty$ and the kernel $R$ is
\begin{displaymath} R(z,u) =
\frac{z(z-1)}{(u-z)u(u-1)}=\frac{1}{u-z}+\frac{z-1}{u}
-\frac{z}{u-1}.
\end{displaymath}
Setting
\begin{equation*}
F[\mu] = \left.\frac{\pa}{\pa \vep} \right|_{\vep=0} w_{\vep\mu}
\quad\text{and}\quad\Phi[\mu]=\left.\frac{\pa}{\pa \bar{\vep}}
\right|_{\vep=0} w_{\vep\mu},
\end{equation*}
we get from \eqref{variation-1}
\begin{align} \label{F-and-Phi}
F[\mu](z)=& -\frac{1}{\pi}\iint\limits_{\UU}R(z,u)\mu(u)\, d^2u,\\
\Phi[\mu](z)=& -\frac{1}{\pi}\iint\limits_{\UU}
R(z,\bar{u})\ov{\mu(u)}\,d^2u. \nonumber
\end{align}
The function $\Phi[\mu](z)$ is holomorphic on $\UU$ and satisfies
\begin{displaymath}
\Phi[\mu]_{zzz}(z) =
-\frac{6}{\pi}\iint\limits_{\UU}\frac{\ov{\mu(u)}}{(\bar{u}-z)^4}\,d^2u.
\end{displaymath}
As it follows from \eqref{projection}, the projection
$P:L^\infty(\U)\rightarrow \Omega^{-1,1}(\U)$ is given by
\begin{equation} \label{projection-1}
(P \mu)(z) = -\frac{3
(z-\z)^2}{\pi}\iint\limits_{\U}\frac{\mu(u)}{(u-\z)^4}d^2 u.
\end{equation}
Equivalently, for $\mu(z) = \tfrac{(z-\bar{z})^2}{2}\ov{\phi(z)}$
with $\phi\in A_\infty(\UU)$, $\Phi[\mu]_{zzz}=\phi$ on $\UU$. The
function $F[\mu]$ satisfies $F[\mu]_{\z}=\mu$ on $\UU$, and is
holomorphic on the lower half-plane $\ov{\UU}$ .
\begin{lemma} \label{integral-4-5}
For $\mu\in\Omega^{-1,1}(\UU)$ and $z\in\U$,
\begin{align*}
\lim_{\vep\rightarrow 0}\iint\limits_{\U(z,\vep)}
\frac{\mu(u)}{(u-z)^4}d^2u = \lim_{\vep\rightarrow
0}\iint\limits_{\U(z,\vep)} \frac{\mu(u)}{(u-z)^5}d^2u = 0,
\end{align*}
where $\U(z,\vep)=\U\setminus\{u\in\U : |u-z|<\vep\}$.
\end{lemma}
\begin{proof} The proof of the first formula essentially follows
the classical Ahlfors' proof in \cite[Lemma 2 in Sect. VI D]{Ahl2}
  by using $\mu(u) = \tfrac{(u-\bar{u})^2}{2}\ov{\phi(u)}$ with
   $\phi\in A_\infty(\UU)$, the identity
\begin{displaymath}
\frac{(u-\bar{u})^2}{(u-z)^4}  = \frac{\pa}{\pa u}\left(
-\frac{1}{u-z} +
\frac{\bar{u}-z}{(u-z)^2} -\frac{1}{3}
\frac{(\bar{u}-z)^2}{(u-z)^3}\right),
\end{displaymath}
and Stokes' theorem. The second formula is proved similarly.
\end{proof}
Another classical result of Ahlfors \cite{Ahl1} is the following.
\begin{lemma}\label{Ahlfors} For $\mu\in\Omega^{-1,1}(\UU)$ and $z\in\U$,
\begin{align*}
F[\mu](z)= \frac{(z-\z)^2}{2} \ov{\Phi^{\prime\prime}(z)} + (z-\z)
\ov{\Phi^{\prime}(z)} + \ov{\Phi(z)},
\end{align*}
where $\Phi(z)=\Phi[\mu](z)$.
\end{lemma}
\begin{remark} It follows from Lemma \ref{Ahlfors} that
$F[\mu]_{zzz}=0$ for $\mu\in\Omega^{-1,1}(\UU)$, in agreement with
Lemma \ref{integral-4-5}.
\end{remark}
\begin{corollary}\label{integral-1}
For $\mu\in\Omega^{-1,1}(\UU)$ and $z\in\U$,
\begin{align*}
\iint\limits_{\U} \frac{\mu(u)}{(u-z) (u-\z)^3}d^2u =0.
\end{align*}
\end{corollary}
\begin{proof}
Using \eqref{F-and-Phi}, we have
\begin{align*}
&F[\mu](z)- \frac{(z-\z)^2}{2} \ov{\Phi^{\prime\prime}(z)} -
(z-\z) \ov{\Phi^{\prime}(z)} - \ov{\Phi(z)} \\
=&-\frac{1}{\pi}\iint\limits_{\UU} \mu(u) \frac{(z-\z)^3}{(u-z)
(u-\z)^3}d^2w.
\end{align*}
\end{proof}
For $\mu\in L^\infty(\UU)_1$ set
\begin{displaymath}
K_{\mu}(u,v) =
\frac{(w_{\mu})_u(u)(w_{\mu})_v(v)}{(w_{\mu}(u) -
w_{\mu}(v))^2}\quad\text{and}\quad K_{\mu}(u,\bar{v}) =
\frac{(w_{\mu})_u(u)\ov{(w_{\mu})_v(v)}}{(w_{\mu}(u) -
\ov{w_{\mu}(v)})^2}.
\end{displaymath}
We have from \eqref{variation-1} the following formulas
\cite{Ahlfors3}
\begin{align} \label{ahlfors-K}
\frac{\pa}{\pa\vep}\, K_{\vep\mu}(z,u) = &
-\frac{1}{\pi}\iint\limits_{\UU}\mu(v)K_{\vep\mu}(z,v)K_{\vep\mu}(v,u)\,d^2v,
\\
\frac{\pa}{\pa\bar{\vep}}\, K_{\vep\mu}(z,u) = &
-\frac{1}{\pi}\iint\limits_{\UU}\ov{\mu(v)}K_{\vep\mu}(z,\bar{v})
K_{\vep\mu}(\bar{v},u)\,d^2v, \nonumber
\end{align}
and
\begin{align} \label{ahlfors-K-bar}
\frac{\pa}{\pa\vep}\, K_{\vep\mu}(z,\bar{u}) = &
-\frac{1}{\pi}\iint\limits_{\UU}\mu(v)K_{\vep\mu}(z,v)K_{\vep\mu}(v,\bar{u})\,d^2v,
\\
\frac{\pa}{\pa\bar{\vep}}\, K_{\vep\mu}(z,\bar{u}) = &
-\frac{1}{\pi}\iint\limits_{\UU}\ov{\mu(v)}K_{\vep\mu}(z,\bar{v})
K_{\vep\mu}(\bar{v},\bar{u})\,d^2v, \nonumber
\end{align}
where the integrals are understood in the principal value sense.

For the model $\HH^2\simeq\Del$ the q.c.~mapping $w_{\mu}$ is
normalized by fixing $-1,-i,1$. The kernel $R$ is given by
\begin{displaymath} R(z,u) =
\frac{(z+1)(z+i)(z-1)}{(u-z)(u+1)(u+i)(u-1)},
\end{displaymath}
and formulas similar to \eqref{F-and-Phi} hold for $F$ and $\Phi$.
In particular, let $f$ be a q.c.~mapping such that
$f|_{\Del}\in\mathcal{D}$, and let $\mu$ be a Beltrami differential
supported on the quasi-disk $\Omega^*=f(\Del^*)$. Let
$v_{t\mu}$ be the solution on $\CC$ of the Beltrami equation
\begin{displaymath}
(v_{t\mu})_{\bar{z}} = t\mu(v_{t\mu})_z,
\end{displaymath}
satisfying $v_{t\mu}(0)=0, v_{t\mu}^\prime(0)=1$ and
$v_{t\mu}^{\prime\prime}(0)=0$. Then
\begin{displaymath}
\dot{v}=\left.\frac{d}{dt}\right|_{t=0}v_{t\mu}
\end{displaymath}
is a holomorphic function on $\Omega=f(\Del)$ and
\begin{equation} \label{infinitesimal}
\dot{v}_{zzz}(z) =
-\frac{6}{\pi}\iint\limits_{\Omega^*}\frac{\mu(u)}{(u-z)^4} d^2u.
\end{equation}
\section{$T(1)$ as a Hilbert manifold}
In this section we are going to endow $T(1)$ with a structure of a
complex manifold modeled on the separable Hilbert space
\begin{align*}
A_2(\Del) = \left\{\phi\;\text{holomorphic on}\; \Del :
\|\phi\|_2^2 = \iint\limits_{\Del}|\phi|^2 \rho^{-1}(z)d^2 z
<\infty\right\}
\end{align*}
of holomorphic functions on $\Del$. In the corresponding topology,
the universal \Te space $T(1)$ is a disjoint union of uncountably
many components on which the right translations act transitively.
\subsection{Hilbert space structure on tangent spaces} Let
\begin{align*}
A_2(\Del^*) = \left\{\phi\;\text{holomorphic on}\; \Del^* :
\|\phi\|_2^2 = \iint\limits_{\Del^*}|\phi|^2 \rho^{-1}(z)d^2 z
<\infty \right\}
\end{align*}
be the Hilbert space of holomorphic functions on $\Del^*$.
\begin{lemma}\label{subspace}
The vector spaces $A_2(\Del)$ and $A_2(\Del^*)$ are subspaces of
$A_{\infty}(\Del)$ and $A_{\infty}(\Del^*)$ respectively. The
natural inclusion maps $A_2(\Del)\hookrightarrow A_{\infty}(\Del)$
and $A_2(\Del^*)\hookrightarrow A_{\infty}(\Del^*)$ are bounded
linear mappings of Banach spaces.
\end{lemma}
\begin{proof}
It is sufficient to consider only the spaces of holomorphic
functions on $\Del$. For every $\phi\in A_2(\Del)$, let $\phi=
\sum_{n=2}^{\infty} (n^3-n) a_n z^{n-2}$ be the power series
expansion. Then
\begin{align*}
\|\phi\|_2^2 = \iint\limits_{\Del}|\phi|^2 \rho^{-1}d^2z =
\frac{\pi}{2}\sum_{n=2}^{\infty} (n^3-n)|a_n|^2,
\end{align*}
and by Cauchy-Schwarz inequality,
\begin{align*}
 |\phi(z)| &= \left|\sum_{n=2}^{\infty}
(n^3-n) a_n z^{n-2}\right| \\ &\leq
\left(\sum_{n=2}^{\infty}(n^3-n)|a_n|^2 \right)^{1/2}\left(
\sum_{n=2}^{\infty}(n^3-n)|z|^{2n-4}\right)^{1/2}
\end{align*}
for every $z\in\Del$. Since
\begin{align*}
\sum_{n=2}^{\infty}(n^3-n)|z|^{2n-4}=\frac{6}{(1-|z|^2)^4},
\end{align*}
we have
\begin{align*}
\|\phi\|_{\infty} = \sup_{z\in \Del}
\left|(1-|z|^2)^2\phi(z)\right| \leq
\sqrt{\frac{12}{\pi}}\,\|\phi\|_2.
\end{align*}
\end{proof}
Let
\begin{align*}
H^{-1,1}(\Del)=&\left\{\mu=\rho^{-1}\bar{\phi},\,\phi\;
\text{holomorphic on}\; \Del :
\|\mu\|^2_2=\iint\limits_{\Del}|\mu|^2 \rho(z)d^2z
<\infty\right\}\\ \intertext{and} \\
H^{-1,1}(\Del^*)=&\left\{\mu=\rho^{-1}\bar{\phi},\,\phi\;
\text{holomorphic on}\; \Del^* :
\|\mu\|^2_2=\iint\limits_{\Del^*}|\mu|^2 \rho(z)d^2z
<\infty\right\}
\end{align*}
be the Hilbert spaces of harmonic Beltrami differentials on $\Del$
and $\Del^*$ respectively. It follows from Lemma \ref{subspace}
that the natural inclusion maps $H^{-1,1}(\Del) \hookrightarrow
\Omega^{-1,1}(\Del)$ and $H^{-1,1}(\Del^*) \hookrightarrow
\Omega^{-1,1}(\Del^*)$ are bounded and under the linear mapping
$D_0\beta$, $ H^{-1,1}(\Del^*)\xrightarrow{\sim}A_2(\Del)$.
\begin{remark} It follows from the proof of Lemma \ref{subspace}
that every $\mu\in H^{-1,1}(\Del)$ (respectively in
$H^{-1,1}(\Del^*)$) satisfies
\begin{align*}
\lim_{|z|\rightarrow 1}\mu(z) =0.
\end{align*}
Indeed, for given $\vep>0$ let $N$ be such that
\begin{align*}
\sum_{n=N}^{\infty}(n^3-n) |a_n|^2 <\vep.
\end{align*}
Then
\begin{align*}
 |\mu(z)| \leq & \frac{(1-|z|^2)^2}{4}\left|\sum_{n=2}^{N-1}
(n^3-n) a_n z^{n-2}\right|\\ &+\frac{(1-|z|^2)^2}{4}
\left(\sum_{n=N}^{\infty}(n^3-n)|a_n|^2 \right)^{1/2}\left(
\sum_{n=2}^{\infty}(n^3-n)|z|^{2n-4}\right)^{1/2}\\ \leq
&\frac{(1-|z|^2)^2}{4}\left|\sum_{n=2}^{N-1} (n^3-n) a_n
z^{n-2}\right|+\frac{\sqrt{6 \vep}}{4},
\end{align*}
so that
\begin{align*}
\limsup_{|z|\rightarrow 1} |\mu(z)| \leq \frac{\sqrt{6\vep}}{4}.
\end{align*}
Since $\vep$ is arbitrary this proves the assertion.
\end{remark}

For every $[\mu]\in T(1)$ let $D_0
R_{[\mu]}\left(H^{-1,1}(\Del^*)\right)$ be the subspace of the
tangent space $T_{[\mu]}T(1)=D_0
R_{[\mu]}\left(\Omega^{-1,1}(\Del^*)\right)$ with a Hilbert space
structure isomorphic to $H^{-1,1}(\Del^*)$. Let $\mathfrak{D}_T$
be the distribution on $T(1)$, defined by the assignment
\begin{displaymath}
T(1)\ni [\mu]\mapsto D_0R_{[\mu]}\left(H^{-1,1}(\Del^*)\right)
\subset T_{[\mu]} T(1).
\end{displaymath}
 Similarly, let $\mathfrak{D}_A$ be the distribution on $A_{\infty}(\Del)$,
defined by
\begin{displaymath}
A_{\infty}(\Del)\ni\phi\mapsto A_2(\Del)\subset T_{\phi}
A_{\infty}(\Del) \simeq A_{\infty}(\Del).
\end{displaymath}

The next statement asserts that under the Bers embedding $\beta:
T(1)\rightarrow A_{\infty}(\Del)$ the distribution
$\mathfrak{D}_T$ is isomorphic to the restriction of the
distribution $\mathfrak{D}_A$ to $\beta(T(1))$.
\begin{theorem} \label{gen-subspace} For every $[\mu]\in T(1)$ the
linear mapping
\begin{equation*}
D_0\left(\beta\circ R_{[\mu]}\right): H^{-1,1}(\Del^*)\rightarrow
A_2(\Del)
\end{equation*}
is a topological isomorphism.
\end{theorem}
\begin{proof}
Let $\nu\in H^{-1,1}(\Del^*)$. Set $w_t =
w_{t\nu\ast\mu}=w_{t\nu}\circ w_{\mu}$ and let
$w_t=\mathrm{g}^{-1}_t\circ \mathrm{f}^t$ be the conformal welding
associated with the q.c. mapping $w_t$ by \eqref{newcon}. Let $v_t
= \rf^t\circ \rf^{-1}$, where $ w_{\mu} = \rg^{-1}\circ \rf$ is
the factorization for $w_{\mu}$, and set
$\Omega=\rf(\Del)=\rg(\Del)$, $\Omega^*=\rf(\Del^*)=\rg(\Del^*)$.
Since
\begin{equation*}
\beta([t\nu\ast\mu]) = \mathcal{S}(\rf^t)= \mathcal{S}(v_t)\circ
\rf \,\rf_z^2 +\mathcal{S}(\rf),
\end{equation*}
we have
\begin{equation*}
D_0\left(\beta\circ R_{[\mu]}\right)(\nu)=
\frac{d}{dt}\Bigr\vert_{t=0}\mathcal{S}(\rf^t)= \dot{v}_{zzz}\circ
\rf\,
\rf_z^2,\quad\text{where}\quad\dot{v}=\frac{d}{dt}\Bigr\vert_{t=0}v_t.
\end{equation*}
The q.c.~mapping $v_t$ is holomorphic on $\Omega$ and satisfies
$v_t\circ \rg = \rg_t \circ w_{t\nu}$. Since $\rg_t$ and $\rg$ are
holomorphic on $\Del^*$, the Beltrami differential of $v_t$ is
given by
\begin{align*}
t\tilde{\nu}(z) = \begin{cases} 0, & z\in\Omega, \\ t(\nu\circ
\rg^{-1})(z)\frac{\ov{\rg^{-1}_z(z)}}{\rg^{-1}_z(z)},&
z\in\Omega^*.
\end{cases}
\end{align*}
It follows from \eqref{infinitesimal} that
\begin{align}\label{varied}
D_0\left(\beta\circ
R_{[\mu]}\right)(\nu)\left(\rf^{-1}(z)\right)\left(\rf^{-1}_z(z)\right)^2
=\dot{v}_{zzz}(z)= -\frac{6}{\pi}
\iint\limits_{\Omega^*}\frac{\tilde{\nu}(u)}{(u-z)^4} d^2u.
\end{align}
Let $\rho_1(z) = (\rho\circ \rf^{-1})(z)|\rf^{-1}_z(z)|^2$ and
$\rho_2(z) =(\rho\circ \rg^{-1})(z)|\rg^{-1}_z(z)|^2$ be the
hyperbolic metric densities on the domains $\Omega$ and $\Omega^*$
respectively. Classical inequalities (see e.g., \cite{Lehto,Nag2})
\begin{align*}
\frac{1}{4}\hspace{0.1cm}\leq \hspace{0.1cm}\eta_i^2(z) \rho_i(z)
\hspace{0.1cm}\leq \hspace{0.1cm} 4, \hspace{1cm} i=1,2,
\end{align*}
where  $\eta_1(z)$ and $\eta_2(z)$ stand, respectively, for the
distances of $z \in\Omega$ and $z\in\Omega^*$ to the quasi-circle
$\rf(S^1)$, yield the following estimates (cf.~\cite[Sect.~
3.4.5]{Nag2})
\begin{align*}
\iint\limits_{\Omega} \frac{d^2z}{|u-z|^4}\leq \iint\limits_{|z-u|
\geq \eta_2(u)}\frac{d^2z}{|u-z|^4}=\frac{\pi}{ \eta_2(u)^2}\leq
4\pi\rho_2(u),\;\; u\in\Omega^*,
\end{align*}
and
\begin{align*}
\iint\limits_{\Omega^*} \frac{d^2u}{|u-z|^4}\leq
4\pi\rho_1(z),\;\; z\in\Omega.
\end{align*}
From here it follows
\begin{align*}
&\Vert D_0\left(\beta\circ R_{[\mu]}\right)(\nu)\Vert^2_2=
\iint\limits_{\Del} \bigl\vert D_0 \left(\beta\circ
R_{[\mu]}\right)(\nu)\bigr\vert^2 \rho^{-1} d^2z =
\iint\limits_{\Omega}\bigl\vert \dot{v}_{zzz}\bigr\vert^2
\rho_1^{-1} d^2z\\ \leq &\frac{6^2}{\pi^2}
\iint\limits_{\Omega}\rho_1(z)^{-1}
\iint\limits_{\Omega^*}\frac{d^2v}{|v-z|^4}\iint\limits_{\Omega^*}
\frac{|\tilde{\nu}(u)|^2 d^2u}{|u-z|^4}\,\,d^2z\\ \leq & \frac{6^2
\cdot 4}{\pi} \iint\limits_{\Omega} \iint\limits_{\Omega^*}
\frac{|\tilde{\nu}(u)|^2 d^2u}{|u-z|^4}\,\,d^2z \leq 6^2\cdot 4^2
\iint\limits_{\Omega^*}|\tilde{\nu}(u)|^2 \rho_2(u)d^2u \\ =&
6^2\cdot 4^2 \iint\limits_{\Del^*} |\nu|^2 \rho(u)d^2u
=576\Vert\nu\Vert^2_2.
\end{align*}

To prove that the mapping $D_0\left(\beta\circ R_{[\mu]}\right)$
is onto, we adapt to our case Bers' arguments, as presented in
\cite[Sect.~3.5]{Nag2}. For $\phi\in A_2(\Del)$ set $q =
(\phi\circ \rf^{-1})(\rf^{-1}_z)^2$, and choose $\mu$ in the
equivalence class of $[\mu]\in T(1)$ to be the conformally natural
extension of $(\rg^{-1}\circ \rf)\vert_{S^1}$, constructed by
Douady and Earle \cite{DE}. Let $h$ be the corresponding
quasiconformal reflection \cite{EN} on $\C$ which fixes the
quasi-circle $\rf(S^1)$. According to the Bers reproducing formula
\cite{Bers2},
\begin{align}\label{reproducing}
q(z) = -\frac{3}{\pi}\iint\limits_{\Omega^*}
\frac{q(h(u))(u-h(u))^2 h_{\bar{u}}(u)}{(u-z)^4}d^2u.
\end{align}
Analogous to $L^{\infty}(\Del^*)$ and $\Omega^{-1,1}(\Del^*)$,
consider the Banach spaces $L^{\infty}(\Omega^*)$ and
\begin{displaymath}
\Omega^{-1,1}(\Omega^*) = \left\{ \mu\in L^{\infty}(\Omega^*): \mu
=  \rho_2^{-1}\bar{q},\; q\;\;\text{is holomorphic on}\;\;
\Omega^* \right\}.
\end{displaymath}
Denote by $\tilde{P}$ the corresponding projection $\tilde{P}:
L^{\infty}(\Omega^*)\rightarrow \Omega^{-1,1}(\Omega^*)$. The
mapping
\begin{displaymath}
\mu\mapsto (\rg^{\ast})^{-1}(\mu) = \mu\circ \rg^{-1}
\frac{\ov{(\rg^{-1})^\prime}}{(\rg^{-1})^\prime}
\end{displaymath}
establishes the isomorphisms $L^{\infty}(\Del^*)\simeq
L^{\infty}(\Omega^*)$ and $\Omega^{-1,1}(\Del^*)\simeq
\Omega^{-1,1}(\Omega^*)$, and $\tilde{P} = (\rg^{\ast})^{-1} \circ
P \circ \rg^{\ast}$. Define $\nu\in \Omega^{-1,1}(\Del^*)$ by
\begin{align*}
(\rg^{\ast})^{-1}(\nu)(z)
=\tilde{P}\left(\frac{1}{2}(q(h(z))(z-h(z))^2h_{\bar{z}}(z)\right)\in
\Omega^{-1,1}(\Omega^*).
\end{align*}
The comparison between \eqref{reproducing} and \eqref{varied}
shows that
\begin{displaymath}
D_0\left(\beta\circ R_{[\mu]}\right)(\nu) = \phi.
\end{displaymath}
To prove that $\nu\in H^{-1,1}(\Del^*)$ we use the Earle-Nag \cite{EN}
estimate,
\begin{align}\label{metric}
\frac{1}{C} \leq |z- h(z)|^4 \rho_1(h(z))\rho_2(z)\leq C,\;\; z\in
\Omega^*,
\end{align}
where the constant $C$ depends only on $\Vert \mu \Vert_{\infty}$.
Since the operator $\tilde{P}$  gives the orthogonal projection of
$L^2(\Omega^*,\rho_2(z)\,d^2z)$ onto $H^{-1,1}(\Omega^*)$, we get
by the Earle-Nag inequality
\begin{align*}
\iint\limits_{\Del^*} |\nu|^2 \rho(z)d^2z=&\iint\limits_{\Omega^*}
|(\rg^{\ast})^{-1}(\nu)(z)|^2\rho_2(z)d^2z\\
 \leq & \frac{1}{4}\iint\limits_{\Omega^*} |q(h(z))(z-h(z))^2 h_{\bar{z}}(z)|^2
\rho_2(z)d^2z\\ \leq & \frac{C}{4} \iint\limits_{\Omega^*}
|q(h(z)) h_{\bar{z}}(z)|^2\rho_1^{-1}(h(z))d^2z.
\end{align*}
Since $h$ is sense reversing, for
\begin{align*}
\kappa = \frac{h_z^{-1}}{h^{-1}_{\bar{z}}}
\end{align*}
we have $\Vert \kappa \Vert_{\infty} < 1$. Now
\begin{align*}
&\iint\limits_{\Omega^*} |q(h(z)) h_{\bar{z}}(z)|^2
\rho_1(h(z))^{-1}d^2z \\ =&\iint\limits_{\Omega} |q(z)|^2
\rho_1(z)^{-1} |(h_{\bar{z}}\circ h^{-1})(z)
h_{\bar{z}}^{-1}(z)|^2(1-|\kappa(z)|^2)d^2z\\ =&
\iint\limits_{\Omega} \frac{|q(z)|^2}{1-|\kappa(z)|^2}
\rho_1(z)^{-1}d^2z\\
 \leq &C_1\iint\limits_{\Omega} |q(z)|^2
\rho_1 (z)^{-1}d^2z =C_1\Vert\phi\Vert^2_2 < \infty,
\end{align*}
so that $\Vert\nu\Vert_2\leq C_2\Vert\phi\Vert_2$. This also
proves that the inverse map to $D_0\left(\beta\circ
R_{[\mu]}\right)$ is bounded, so that $D_0\left(\beta\circ
R_{[\mu]}\right)$ is a topological isomorphism.
\end{proof}
\begin{remark} \label{estimate-bers}
It follows from the proof of the first part of Theorem
\ref{gen-subspace} that $D_0\left(\beta\circ
R_{[\mu]}\right)=D_0\left(\beta\circ\Phi\circ R_\mu\right)$
extends to be a bounded linear operator on
$L^2(\Del^*,\rho(z)d^2z)$ and the estimate
\begin{displaymath}
\Vert D_0\left(\beta\circ\Phi\circ R_\mu\right)(\nu)\Vert_2=\Vert
D_\mu \left(\beta\circ\Phi\right)\left(D_0
R_\mu(\nu)\right)\Vert_2 \leq 24\Vert\nu\Vert_2
\end{displaymath}
holds for all $\nu\in L^2(\Del^*,\rho(z)d^2z)$ and $\mu\in
L^\infty(\Del^*)_1$.
\end{remark}
\subsection{The $L^2$-estimates}
The lemmas below are needed for the rigorous definition of a
complex Hilbert manifold structure on $T(1)$.
\begin{lemma}\label{appr}
For every $\vep >0$ there exists $0<\delta<1$  such that for all
$\mu \in \Omega^{-1,1}(\Del^*)$ with $\Vert \mu\Vert_{\infty} <
\delta$,
\begin{displaymath}
\left\vert\frac{|(w_{\mu})_z(z) |^2}{(1-|w_{\mu}(z)|^2)^2} -
\frac{1}{(1-|z|^2)^2}\right\vert < \frac{\vep}{(1-|z|^2)^2}
\end{displaymath}
for all $z\in\Del\cup\Del^*$. The same inequality holds for
$w_{\mu^{-1}}=w_{\mu}^{-1}$.
\end{lemma}
\begin{proof}
Using the isomorphism
$\Omega^{-1,1}(\Del^*)\xrightarrow{\sim}\Omega^{-1,1}(\Del)$ given
by reflection \eqref{sym} and  property \eqref{reflection}, it is
sufficient to prove the estimate for $z\in\Del$. Since
$\g_\mu=\alpha\circ w_\mu$, where $\alpha\in\PSU(1,1)$, the
estimate holds for $w_\mu$ if and only if it holds for $\g_\mu$.
By Lemma \ref{g(0)=0} $\g_\mu$ fixes $0$ and $\infty$, and by the
result of Ahlfors and Bers in \cite{AB} (see also the remark of
Bers in \cite{Bers1}) the functional
\begin{align*}
\mathcal{L}^\infty(\Del^*)_1\ni\mu  \mapsto (\g_{\mu})_z (0)\in\CC
\end{align*}
is real-analytic at $\mu=0$. In particular, for every $\vep
>0$ there exists $0<\delta<1$ such that for all
$\mu\in\Omega^{-1,1}(\Del^*)_1$ with $\Vert \mu \Vert_{\infty} <
\delta$,
\begin{align*}
\Bigl\vert\vert(\g_{\mu})_z (0)\vert^2-1\Bigr\vert < \vep.
\end{align*}
For $z\in \Del$, let $\tilde{\mu} = \mu \circ \sigma_z $, where
$\sigma_z (w) = \frac{w+ z}{1+ \ov{z} w}$. Then $\g_{\tilde{\mu}}
= \tilde{\sigma}_{z} \circ \g_{\mu} \circ \sigma_z$ for some
$\tilde{\sigma}_z \in \PSU(1,1)$. Since
$\tilde{\mu}\in\Omega^{-1,1}(\Del^*)_1$, it follows from Lemma
\ref{g(0)=0} that $\g_{\tilde{\mu}}(0)=0$. Therefore from
$\tilde{\sigma}_z(\g_{\mu}(z))=0$, one obtains
\begin{align*}
\frac{|(\g_{\mu})_z(z) |^2}{(1-|\g_{\mu}(z)|^2)^2} (1-|z|^2)^2 =
\left\vert (\g_{\tilde{\mu}})_z(0)\right\vert^2.
\end{align*}
Since $\Vert \tilde{\mu} \Vert_{\infty} = \Vert\mu
\Vert_{\infty}$, the assertion follows. On the other hand, if
$\mu\in\Omega^{-1,1}(\Del^*)_1$, then
$\mu^{-1}\in\mathcal{L}^\infty(\Del^*)_1$ and
$\Vert\mu^{-1}\Vert_\infty = \Vert\mu\Vert_\infty$. Hence the
assertion also holds for $w_{\mu^{-1}}$.
\end{proof}
\begin{corollary} \label{righttran} Let
$\mu\in\Omega^{-1,1}(\Del^*),\,\Vert\mu\Vert_{\infty}< \delta$,
where $\delta$ corresponds to $\vep=1$ in the previous lemma. Then
for every $\lambda\in L^{\infty}(\Del^*)_1$ the linear mapping
$D_{\lambda} R_{\mu}$ extends to an invertible bounded linear
operator on the Hilbert space $L^2(\Del^*, \rho(z)d^2z)$.
Moreover,
\begin{align*} \Bigl\Vert D_{\lambda}
R_{\mu}(\nu)\Bigr\Vert_2\leq\frac{\sqrt{2}}{(1-\Vert\mu\Vert_{\infty})^2}
\Vert \nu \Vert_2,
\end{align*}
for all $\nu\in L^2(\Del^*,\rho(z)d^2z)$ and $\lambda\in
L^{\infty}(\Del^*)_1$. The same inequality holds for $D_{\lambda}
R_{\mu^{-1}}$.
\end{corollary}
\begin{proof} Since
\begin{equation*}
D_{\lambda}R_{\mu}(\nu) =  \frac{(1-|\mu|^2)\,\nu \circ
w_{\mu}}{\left(1+ \ov{\mu} \lambda\circ
w_{\mu}\frac{\ov{(w_{\mu})_z}}{(w_{\mu})_z}\right)^2}\,
\frac{\ov{(w_{\mu})_z}}{(w_{\mu})_z},
\end{equation*}
and
\begin{equation*}
\left\Vert\frac{1-|\mu|^2}{\left(1+ \ov{\mu} \lambda\circ
w_{\mu}\frac{\ov{(w_{\mu})_z}}{(w_{\mu})_z}\right)^2}\right\Vert_{\infty}\leq
\frac{1}{(1-\Vert\mu\Vert_{\infty})^2}
\end{equation*}
for all $\lambda\in L^{\infty}(\Del^*)_1$, we have by using Lemma
\ref{appr} and
$\Vert\mu\Vert_{\infty}=\Vert\mu^{-1}\Vert_{\infty}$,
\begin{align*}
&\iint\limits_{\Del^*} \Bigl|D_{\lambda}R_{\mu}(\nu)\Bigr|^2
\rho(z)d^2z \leq(1-\Vert\mu\Vert_\infty)^{-4}
\iint\limits_{\Del^*} \left\vert\nu\circ w_{\mu}
\frac{\ov{(w_{\mu})_z}}{(w_{\mu})_z}\right\vert^2 \rho(z) d^2 z\\
=& (1-\Vert\mu\Vert_\infty)^{-4} \iint\limits_{\Del^*} \left\vert
\nu\right\vert^2 \frac{4|(w_{\mu^{-1}})_z|^2}{(1-
|w_{\mu^{-1}}|^2)^2} (1-|\mu^{-1}|^2) d^2z\\ \leq &
2(1-\Vert\mu\Vert_\infty)^{-4} \iint\limits_{\Del^*} |\nu|^2
\rho(z) d^2z =2(1-\Vert\mu\vert_\infty)^{-4}\Vert\nu\Vert^2_2.
\end{align*}
Replacing everywhere $\mu$ by $\mu^{-1}$ we get the same estimate
for $D_{\lambda} R_{\mu^{-1}}$.
\end{proof}

Denote by $\mathcal{O}(\Del^*)_1$ the subgroup of
$L^{\infty}(\Del^*)_1$ generated by $\mu\in\Omega^{-1,1}(\Del^*),
\Vert\mu\Vert_\infty<\delta$, where $\delta$ is as in Corollary
\ref{righttran}.
\begin{lemma} \label{right1}
For every $\mu\in \mathcal{O}(\Del^*)_1$ there exists $C>0$ such
that
\begin{equation*}
\Vert R_{\mu}(\lambda_1)-R_{\mu}(\lambda_2)\Vert_2<C\Vert
\lambda_1-\lambda_2\Vert_2
\end{equation*}
for all $\lambda_1, \lambda_2 \in L^{\infty}(\Del^*)_1$ satisfying
$\lambda_1 - \lambda_2\in L^2(\Del^*, \rho(z)d^2z)$.
\end{lemma}
\begin{proof}
Suppose first that $\Vert\mu\Vert_{\infty}<\delta$. Set
$\lambda(t)=\lambda_1 + t\nu$, where $\nu =\lambda_2-\lambda_1$,
so that $\lambda(t)\in L^{\infty}(\Del^*)_1$, $0\leq t\leq 1$. By
fundamental theorem of calculus,
\begin{align*}
R_{\mu}(\lambda_1)-R_{\mu}(\lambda_2) = &\int^{1}_0 \frac{d}{dt}
R_{\mu}(\lambda(t))dt \\ = &\int^{1}_0 D_{\lambda(t)}R_{\mu}(\nu)
dt.
\end{align*}
Using Corollary \ref{righttran},
\begin{align*}
\Vert R_{\mu}(\lambda_1)-R_{\mu}(\lambda_2)\Vert^2_2
=&\iint\limits_{\Del^*} \left\vert\int^{1}_0
D_{\lambda(t)}R_{\mu}(\nu)(z)\right\vert^2 \rho(z) d^2 z \\ \leq
&\int^1_0\left(\iint\limits_{\Del^*} \left\vert
D_{\lambda(t)}R_{\mu}(\nu)(z)\right\vert^2 \rho(z) d^2 z\right) dt
\\
\leq & C^2\Vert\nu\Vert^2_2 = C^2\Vert \lambda_1 -
\lambda_2\Vert^2_2.
\end{align*}
The same estimate also holds for $R_{\mu}^{-1}$.

Since every $\mu\in\mathcal{O}(\Del^*)_1$ can be written as
$\mu^{\vep_n}_n\ast\dots\ast\mu^{\vep_1}_1$, where
$\mu_i\in\Omega^{-1,1}(\Del^*),\,\Vert\mu_i\Vert_{\infty}<\delta$,
and $\vep_i=\pm 1,\;i=1,\dots,n$, we have
\begin{displaymath}
R_{\mu} = R^{\vep_1}_{\mu_1}\circ\dots\circ R^{\vep_n}_{\mu_n},
\end{displaymath}
and the assertion of the lemma follows.
\end{proof}
\begin{remark} \label{R-estimate}
Applying the same argument, we get from Corollary \ref{righttran}
that for every $\mu\in\mathcal{O}(\Del^*)_1$ there exists $C>0$,
depending only on $\Vert\mu\Vert_{\infty}$ such that
\begin{displaymath}
\Bigl\Vert D_{\lambda} R_{\mu}(\nu)\Bigr\Vert_2\leq C \Vert \nu
\Vert_2,
\end{displaymath}
for all $\nu\in L^2(\Del^*, \rho(z)d^2z)$ and $\lambda\in
L^{\infty}(\Del^*)_1$.
\end{remark}
\begin{lemma}\label{rightinv}
For every $\mu \in \mathcal{O}(\Del^*)_1$ there exists $C>0$ such
that
\begin{align*}
\Vert (\beta\circ\Phi)(\lambda\ast\mu) -
(\beta\circ\Phi)(\mu)\Vert_2 \leq C\Vert\lambda\Vert_2
\end{align*}
for all $\lambda \in L^2(\Del^*,\,\rho(z)d^2z)\cap
L^{\infty}(\Del^*)_1$.
\end{lemma}
\begin{proof}
Set $\phi(t)=(\beta\circ\Phi)(t\lambda\ast\mu)$. By fundamental
theorem of calculus,
\begin{align*}
(\beta\circ\Phi)(\lambda\ast\mu)-(\beta\circ\Phi)(\mu) =
\int^{1}_0 \frac{d\phi}{dt}(t) dt,
\end{align*}
where
\begin{align*}
\frac{d\phi}{dt}(t) = D_{t\lambda}\left(\beta\circ\Phi\circ
R_{\mu}\right)(\lambda) =
\left(D_{t\lambda\ast\mu}(\beta\circ\Phi)\circ
D_{t\lambda}R_{\mu}\right)(\lambda)
\end{align*}
by chain rule. Since $(D_0 R_\mu)^{-1} = D_\mu R_{\mu^{-1}}$, it
follows from Remarks \ref{estimate-bers} and \ref{R-estimate} that
\begin{equation*}
\Vert D_{t\lambda\ast\mu}(\beta\circ\Phi)(\nu)\Vert_2\leq 24\Vert
D_{t\lambda\ast\mu} R_{(t\lambda\ast\mu)^{-1}}(\nu)\Vert_2\leq
C_1\Vert\nu\Vert_2.
\end{equation*}
Using Remark \ref{R-estimate} again, we get
\begin{align*}
\left\Vert\frac{d\phi}{dt}(t)\right\Vert_2=\Vert\left(D_{t\lambda\ast\mu}
(\beta\circ\Phi)\circ D_{t\lambda}R_{\mu}\right)(\lambda)\Vert_2
\leq C_2 \Vert\lambda\Vert_2,\;\;0\leq t\leq 1.
\end{align*}
Therefore,
\begin{align*}
\bigl\Vert (\beta\circ\Phi)(\lambda \ast\mu)-
(\beta\circ\Phi)(\mu) \bigr\Vert_2^2 =&
\iint\limits_{\Del}\left|\int_0^1 \frac{d\phi}{dt}(t,z)dt\right|^2
\rho(z)^{-1} d^2z\\ \leq &\int_0^1 \left(\iint\limits_{\Del}
\bigl|\frac{d\phi}{dt}(t,z)\bigr|^2\rho(z)^{-1} d^2z\right) dt\\
\leq & C_2^2 \bigl\Vert \lambda\bigr\Vert_2^2,
\end{align*}
which concludes the proof.
\end{proof}
\subsection{The Hilbert manifold structure of $T(1)$}
\label{hilbert-manifold}
 For every $\mu\in\mathcal{O}(\Del^*)_1$ let
$V_\mu\subset U_\mu\subset T(1)$ be the image under the map
$h_\mu^{-1} =\Phi\circ R_{\mu}\circ\Lambda$ of the open ball of
radius $\sqrt{\pi/3}$ about the origin in $A_2(\Del)$, which by
Lemma \ref{subspace} is contained in the ball of radius $2$ in
$A_{\infty}(\Del)$. Here $(U_\mu,h_\mu)$ is the coordinate chart
$U_{\mu}$ of the complex-analytic atlas for $T(1)$ as a complex
Banach manifold (see Section \ref{complex1}). Let
\begin{displaymath}
\tilde{h}_{\mu}= h_{\mu}\bigr\vert_{V_{\mu}}: V_\mu\rightarrow
A_2(\Del).
\end{displaymath}
The main result of this subsection is the following.
\begin{theorem} \label{Hilbert} For every $\mu, \nu\in
\mathcal{O}(\Del^*)_1$ the sets $\tilde{h}_\mu(V_\mu\cap V_\nu)$
and $\tilde{h}_\nu(V_\mu\cap V_\nu)$ are open in $A_2(\Del)$ and
the map
\begin{displaymath}
\tilde{h}_{\mu\nu}=\tilde{h}_\mu\circ \tilde{h}_\nu^{-1}:
\tilde{h}_\nu(V_\mu \cap V_\nu)\longrightarrow
\tilde{h}_\mu(V_\mu\cap V_\nu)\subset A_2(\Del)
\end{displaymath}
is a biholomorphic function on the Hilbert space $A_2(\Del)$.
\end{theorem}
\begin{proof} First we prove that the sets $\tilde{h}_\mu(V_\mu\cap V_\nu)$
and $\tilde{h}_\nu(V_\mu\cap V_\nu)$ are open in $A_2(\Del)$.
Since $V_\mu\cap V_\nu\neq\emptyset$ (otherwise there is nothing
to prove), there exist $\phi_1\in\tilde{h}_\mu(V_\mu\cap V_\nu)$
and $\phi_2\in\tilde{h}_\nu(V_\mu\cap V_\nu)$,
$\Vert\phi_1\Vert_2, \Vert\phi_2\Vert_2< \sqrt{\pi/3}$, such that
$\tilde{h}_{\mu}^{-1}(\phi_1)=\tilde{h}_{\nu}^{-1}(\phi_2)$, i.e.,
\begin{displaymath}
(\Phi\circ R_\mu\circ\Lambda)(\phi_1) = (\Phi\circ
R_\nu\circ\Lambda)(\phi_2).
\end{displaymath}
Setting $\lambda_1=\Lambda(\phi_1), \lambda_2=\Lambda(\phi_2)$ and
$\kappa=\nu\ast\mu^{-1}$, we get
\begin{displaymath}
\Phi(\lambda_1) = \Phi(\lambda_2\ast\kappa).
\end{displaymath}

The sets $h_\mu(U_\mu\cap U_\nu)$ and $h_\nu(U_\mu\cap U_\nu)$ are
open in $A_\infty(\Del)$, so that there exists $\delta_1>0$ such
that $h_\mu(U_\mu\cap U_\nu)$ contains a ball of radius $\delta_1$
about $\phi_1$ in $A_\infty(\Del)$. The mapping $h_{\mu\nu}:
h_\nu(U_\mu\cap U_\nu)\rightarrow h_\mu(U_\mu\cap U_\nu)\subset
A_\infty(\Del)$ is a continuous function in the Banach space
$A_\infty(\Del)$, so that there exists $\delta_2>0$ such that the
inverse image by $h_{\mu\nu}$ of the ball of radius $\delta_1$
about $\phi_1$ in $A_\infty(\Del)$ contains the ball of radius
$\delta_2$ about $\phi_2$ in $A_\infty(\Del)$. According to Lemma
\ref{subspace}, the latter ball contains any ball of radius
$\delta_3<\sqrt{\pi/12}\,\delta_2$ about $\phi_2$ in $A_2(\Del)$.
Now for every $\varphi_2\in A_2(\Del)$ satisfying $\Vert\varphi_2
-\phi_2\Vert_2<\delta_3$ set
\begin{displaymath}
\varphi_1 =h_{\mu\nu}(\varphi_2) = (\beta\circ\Phi\circ
R_\kappa\circ\Lambda)(\varphi_2).
\end{displaymath}
We claim that $\delta_3>0$ can be chosen such that $\varphi_1\in
A_2(\Del)$ and $\Vert\varphi_1\Vert_2 <\sqrt{\pi/3}$, which
implies that $\tilde{h}_\nu(V_\mu\cap V_\nu)$ contains the ball of
radius $\delta_3$ about $\phi_2$ in $A_2(\Del)$. Indeed, set
$\lambda = \Lambda(\varphi_2)$, so that $\varphi_1 =
(\beta\circ\Phi)(\lambda\ast\kappa)$. Since $\lambda -
\lambda_2\in L^2(\Del^*,\rho(z)d^2z)$, we have by Lemmas
\ref{rightinv} and \ref{right1},
\begin{align*}
\Vert\varphi_1 -\phi_1\Vert_2 & = \Vert
(\beta\circ\Phi)(\lambda\ast\kappa) -
(\beta\circ\Phi)(\lambda_2\ast\kappa)\Vert_2 \\ & \leq
C\Vert\lambda\ast\lambda_2^{-1}\Vert_2 \leq C^2\Vert\lambda -
\lambda_2\Vert_2 \\ & = 2C^2\Vert\varphi_2 - \phi_2\Vert_2
<2C^2\delta_3,
\end{align*}
where the constant $C>0$ (chosen to be the same for both Lemmas
\ref{right1} and \ref{rightinv}) depends only on $\lambda_2$ and
$\kappa$. Choosing $\delta_3$ small enough we have
$\Vert\varphi_1\Vert_2 <\sqrt{\pi/3}$.

The same argument applied to the map
$\tilde{h}_{\nu\mu}=\tilde{h}_{\mu\nu}^{-1}$ proves that
$\tilde{h}_\mu(V_\mu\cap V_\nu)$ is open in $A_2(\Del)$.

It remains to prove that the map $\tilde{h}_{\mu\nu}$ is a
holomorphic function in the Hilbert space $A_2(\Del)$. It is
bounded, so according to \cite{bourbaki} it is sufficient to prove
that for every $\varphi\in\tilde{h}_\nu(V_\mu\cap V_\nu)$ and
every $\eta\in A_2(\Del)$ the mapping $\CC\ni
t\mapsto\phi(t)=\tilde{h}_{\mu\nu}(\varphi+ t\eta)\in A_2(\Del)$
is a holomorphic function in some neighborhood of $0$ in $\CC$.
For this purpose we use the standard argument based on the fact
that the map $h_{\mu\nu}$ is already a holomorphic function in the
Banach space $A_{\infty}(\Del)$ and the mapping $\CC\ni t\mapsto
\phi(t) \in A_\infty(\Del)$ is a holomorphic function in some
neighborhood of $0$ in $\CC$. Thus there exists $\delta>0$ such
that for every $|t_0|<\delta$,
\begin{equation*}
\left\Vert \phi(t)-\phi(t_0) - (t-t_0)\frac{d
\phi}{dt}(t_0)\right\Vert_\infty =o(|t-t_0|)\quad\text{as}\quad
t\rightarrow t_0.
\end{equation*}
Moreover, $\delta$ can be chosen such that $\varphi +
t\eta\in\tilde{h}_\nu(V_\mu\cap V_\nu)$ for $|t|<\delta$. Then for
every $z\in\Del$ the complex-valued function $\phi(t)(z)$ is
holomorphic on $|t|<\delta$ and
\begin{gather*}
\phi(t,z) - \phi(t_0,z) - (t-t_0)\frac{d\phi}{dt}(t_0,z) \\=
\frac{1}{2\pi i}\oint_{|w-t_0|=\delta_1}\phi(w,z)
\left(\frac{1}{w-t}-\frac{1}{w-t_0} -
\frac{t-t_0}{(w-t_0)^2}\right) dw\\ =\frac{(t-t_0)^2}{2\pi i}
\oint_{|w-t_0| =\delta_1} \frac{\phi(w,z)}{(w-t_0)^2(w-t)}dw,
\end{gather*}
where $\delta_1>0$ is such that the disk of radius $\delta_1$
about $t_0$ is inside the disk of radius $\delta$ about the
origin, and $t$ satisfies $|t-t_0|<\delta_1$. Since $\phi(t)\in
\tilde{h}_\mu(V_\mu\cap V_\nu)$, $\Vert\phi(t)\Vert^2_2<\pi/3$ for
all $|t|<\delta$, and we have
\begin{align*}
&\left\Vert \frac{\phi(t) - \phi(t_0)}{t-t_0}
-\frac{d\phi}{dt}(t_0)\right\Vert_2^2\\ \leq
&\frac{|t-t_0|^2}{4\pi^2} \oint_{|w|=\delta_1}
\frac{|dw|}{|w|^4|w-(t-t_0)|^2}\oint_{|w|=\delta_1} \bigl\Vert
\phi(w + t_0)\bigr\Vert_2^2 |dw|\\ =&
O(|t-t_0|^2)\quad\text{as}\quad t\rightarrow t_0.
\end{align*}
\end{proof}
According to \cite{bourbaki}, Theorem \ref{Hilbert} justifies the
following definition.
\begin{definition} The covering
\begin{displaymath}
T(1)=\bigcup_{\mu\in\mathcal{O}(\Del^*)_1}  V_\mu
\end{displaymath}
with the coordinate maps $\tilde{h}_\mu: V_{\mu}\longrightarrow
A_2(\Del)$ and the transition maps
\begin{displaymath}
\tilde{h}_{\mu\nu}=\tilde{h}_\mu\circ \tilde{h}_\nu^{-1}:
\tilde{h}_\nu(V_\mu\cap V_\nu)\longrightarrow
\tilde{h}_\mu(V_\mu\cap V_\nu)
\end{displaymath}
is a complex-analytic atlas which endows $T(1)$ with the structure
of a complex Hilbert manifold modeled on the Hilbert space
$A_2(\Del)$.
\end{definition}
\begin{corollary} \label{righttrans-hilbert} The right translations
are biholomorphic mappings on the Hilbert manifold $T(1)$.
\end{corollary}
\begin{proof}
Representing a point in $T(1)$ by $\mu\in\mathcal{O}(\Del^*)_1$ we
have $R_{[\mu]}(V_\lambda)=V_{\lambda\ast\mu}$, so that
$\tilde{h}_{\lambda\ast\mu}\circ
R_{[\mu]}\circ\tilde{h}_{\lambda}^{-1}$ is the identity mapping on
$\tilde{h}_{\lambda}(V_\lambda)\subset A_2(\Del)$.
\end{proof}
We will continue to use the name Bers coordinates for the complex coordinates $(V_\mu,\tilde{h}_\mu)$ on the Hilbert manifold
$T(1)$. As in Section \ref{complex1}, the vector
field $\tfrac{\del}{\del \vep_\nu}$ corresponding to $\nu\in
H^{-1,1}(\Del^*)$ at a point $[\mu]\in V_0$ in terms of the Bers
coordinates on $V_\mu$ has the same form \eqref{bers-change},
i.e.,
\begin{equation*}
\left.\dfrac{\del}{\del \vep_\nu}\right\vert_{[\mu]} = P\left(
\left(\dfrac{\nu}{1 -
|\mu|^2}\:\dfrac{(w_\mu)_z}{(\ov{w}_\mu)_{\z}} \right)\circ
w_\mu^{-1}\right),
\end{equation*}
where $P: L^2(\Del^*,\rho(z)d^2 z)\rightarrow H^{-1,1}(\Del^*)$ is
the orthogonal projector given by \eqref{projection}.

\subsection{Integral manifolds of the distribution
$\mathfrak{D}_T$}\label{integralmanifold} Finally, we introduce a
Hilbert manifold structure on the Banach space $A_\infty(\Del)$ by
defining the coordinate chart at every $\phi\in A_\infty(\Del)$ to
be $\phi \,+ A_2(\Del)$. By Lemma \ref{subspace} the Hilbert
manifold topology on $A_\infty(\Del)$ is stronger than the Banach
space topology. The Hilbert manifold $A_\infty(\Del)$ is not
connected. Rather $A_\infty(\Del)$ is the union of uncountably
many components $\phi \,+ A_2(\Del)$ with $\phi\in
A_\infty(\Del)/A_2(\Del)$, which are integral manifolds of the
distribution $\mathfrak{D}_A$.
\begin{theorem} \label{bers-hilbert} The Bers embedding
$\beta: T(1)\rightarrow\beta(T(1))\subset A_\infty(\Del)$ is a
biholomorphic mapping of Hilbert manifolds.
\end{theorem}
\begin{proof}
To prove that the Bers embedding is holomorphic it is sufficient
to show that for every $\mu\in\mathcal{O}(\Del^*)_1$ the image of
the ball of radius $\sqrt{\pi/3}$ about $0$ in $A_2(\Del)$ by the
mapping $\beta\circ\tilde{h}^{-1}_\mu=\beta\circ\Phi\circ
R_\mu\circ\Lambda$  is inside a translate by
$(\beta\circ\Phi)(\mu)$ of some ball about $0$ in $A_2(\Del)$.
This immediately follows from Lemma \ref{rightinv},
\begin{align*}
\left\Vert\left(\beta\circ\tilde{h}^{-1}_\mu\right)(\varphi) -
\left(\beta\circ\Phi\right)(\mu)\right\Vert_2& =
\Vert\left(\beta\circ\Phi\right)(\lambda\ast\mu) -
\left(\beta\circ\Phi\right)(\mu)\Vert_2
\\ & <C\Vert\lambda\Vert_2,
\end{align*}
where $\lambda=\Lambda(\varphi)\in L^2(\Del^*,\rho(z)d^2z) \cap
L^{\infty}(\Del^*)_1$ and the constant $C>0$ depends only on
$\Vert\mu\Vert_\infty$. Since the Bers embedding is a holomorphic
mapping of Banach manifolds, the standard argument used in the
proof of Theorem \ref{Hilbert} works for this case, so that the
mapping $\beta\circ\tilde{h}^{-1}_\mu - (\beta\circ\Phi)(\mu)$ is
a holomorphic function on the Hilbert space $A_2(\Del)$.

Finally, the image $\beta(T(1))$ is open in the Hilbert manifold
$A_\infty(\Del)$ since it is open in a weaker Banach manifold
topology. Using Theorem \ref{gen-subspace} and the inverse
function theorem for Hilbert manifolds \cite{Lang2} we see that
the Bers embedding is biholomorphic.
\end{proof}

 Theorem \ref{bers-hilbert} allows us to conclude that the distribution
$\mathfrak{D}_T$ on $T(1)$ is equivalent to the restriction of the
distribution $\mathfrak{D}_A$ on $\beta(T(1))\subset
A_\infty(\Del)$, and therefore is integrable. Its integral
manifolds are inverse images by the Bers embedding $\beta$ of the
integral manifolds of the distribution $\mathfrak{D}_A$ on
$\beta(T(1))$, i.e., of the components $\left(\phi\, +
A_2(\Del)\right)\cap\beta(T(1))$. For every $[\mu]\in T(1)$ denote
by $T_{[\mu]}(1)$ the component of the Hilbert manifold $T(1)$
containing $[\mu]$. It follows from Theorems \ref{gen-subspace}
and \ref{bers-hilbert} that the Hilbert manifold $T_{[\mu]}(1)$ is
the integral manifold of the distribution $\mathfrak{D}_T$ passing
through $[\mu]\in T(1)$. The right translations act transitively
on the set of components, i.e.
$R_{[\nu]}(T_{[\mu]}(1))=T_{[\mu\ast\nu]}(1)$ for all
$[\mu],[\nu]\in T(1)$.

We will denote the component of $0\in T(1)$ by $T_0(1)$. We will
prove in Appendix A that $T_0(1)$ is the inverse image of
$\beta(T(1))\cap A_2(\Del)$ under the Bers embedding. Moreover, we
will show that $T_0(1)$ is a topological group which is the
closure of $\Mob(S^1)\bk \Diff_+(S^1)$ in $T(1)$ in the Hilbert
manifold topology.

\section{Velling-Kirillov and Weil-Petersson metrics}
\subsection{Velling-Kirillov metric on the universal \Te curve}
The Velling-Kirillov metric is a right-invariant Hermitian metric
on $\mathcal{T}(1)$, defined at the origin of $\mathcal{T}(1)$ by
\begin{equation} \label{VK-metric}
\Vert \text{v} \Vert_{VK}^2 = \sum^\infty_{n=1}n |c_n|^2,
\end{equation}
where
\begin{displaymath}
\text{v}=\dfrac{\text{u} - i J\,\text{u}}{2}\quad\text{and}\quad
\text{u}=\sum_{n\in\ZZ\setminus\{0\}} c_n
e^{in\theta}\frac{d}{d\theta}\in T_0^{\RR}
S^1\backslash\text{Homeo}_{qs}(S^1).
\end{displaymath}
The convergence of the series is guaranteed by the property
\textbf{TS4} (with $s=1/2$) in Section \ref{TS}. The
Velling-Kirillov metric is a smooth right-invariant \Ka metric on
the complex Banach manifold $\mathcal{T}(1)$. Its symplectic form
$\omega_{VK}$ at the origin of $\mathcal{T}(1)$ is given by
\begin{displaymath}
\omega_{VK}(\text{v},\bar{\text{v}}) = \frac{i}{2}\Vert \text{v}
\Vert_{VK}^2,
\end{displaymath}
In the following section we prove that the Velling-Kirillov metric
is real-analytic on $\mathcal{T}(1)$ by presenting its
real-analytic \Ka potential.
\begin{remark}
For the homogeneous space $S^1\backslash\Diff_+(S^1)$ the metric
was introduced in this form by A.A.~Kirillov \cite{Ki} and it has
been studied by A.A.~Kirillov and D.~Yuriev. In particular, in
\cite{KY} it was shown to be \Ka. In \cite{Ve}, J. Velling has
introduced a Hermitian metric for the space $\mathcal{T}(1)$ using
arguments from  geometric function theory. In \cite{Teo},
Kirillov's definition was extended to $\mathcal{T}(1)$ and it was
shown that the resulting metric coincides with the metric
introduced by Velling. The Velling-Kirillov metric is the unique
right-invariant \Ka metric on the universal \Te curve
$\mathcal{T}(1)$ \cite{Teo}.
\end{remark}
\subsection{Weil-Petersson metric on the universal
Teichm\"uller space}\label{WP-metric} In this section we consider
$T(1)$ as a Hilbert manifold. The Weil-Petersson metric on $T(1)$
is a Hermitian metric defined by the Hilbert space inner product
on tangent spaces, which are identified with the Hilbert space
$H^{-1,1}(\Del^*)$ by right translations (see Section
\ref{hilbert-manifold}). Thus the Weil-Petersson metric is
a right-invariant metric on $T(1)$ defined at the origin of $T(1)$
by
\begin{equation} \label{scalarproduct}
\langle \mu,\nu\rangle_{WP} =
\iint\limits_{\Del^*}\mu\bar{\nu}\rho(z)d^2z,\;\;\mu,\nu\in
H^{-1,1}(\Del^*) = T_0 T(1).
\end{equation}

To every $\mu\in H^{-1.1}(\Del^*)$ there corresponds a vector
field $\tfrac{\pa}{\pa\vep_\mu}$ over $V_0$, given by
\eqref{bers-change}-\eqref{projection}. We set for every
$\kappa\in V_0$,
\begin{equation} \label{WP-explicit}
g_{\mu\bar{\nu}}(\kappa) = \left\la \left.\tfrac{\pa}{\pa
\vep_\mu}\right|_{\kappa}, \left.\tfrac{\pa}{\pa
\vep_\nu}\right|_{\kappa}\right\ra_{WP} = \iint\limits_{\Del^*}
P(R(\mu, \kappa)) \ov{P(R(\nu, \kappa))} \rho(z)d^2z.
\end{equation}
This formula explicitly defines the Weil-Petersson metric on the
coordinate chart $V_0$. The Weil-Petersson metric extends to other
charts $V_{\mu}$ by right translations.

The following statement is an easy consequence of Lemma
\ref{appr}.
\begin{lemma}\label{neighbour} The Weil-Petersson metric is
continuous on $T(1)$.
\end{lemma}
\begin{proof}
As it follows from Corollary \ref{righttrans-hilbert}, it is
sufficient to prove that for every $\mu\in H^{-1,1}(\Del^*)$ the
function $g_{\mu\bar{\mu}}$ is continuous on $V_0$ at $0$. Since
the embedding $V_0\hookrightarrow U_0$ is continuous by Lemma
\ref{subspace}, it is sufficient to prove that the function
$g_{\mu\bar{\mu}}$ is defined on a neighborhood of $0$ in $U_0$
and is continuous at $0$.

Since the projector $P$ is norm-decreasing,
\begin{align*}
g_{\mu\bar{\mu}} (\kappa) &= \iint\limits_{\Del^*} P(R(\mu,
\kappa)) \ov{P(R(\mu, \kappa))} \rho(z)d^2z \\ &\leq
\iint\limits_{\Del^*} R(\mu, \kappa) \ov{R(\mu, \kappa)}
\rho(z)d^2z
\\ &=4\iint\limits_{\Del^*}
\frac{|\mu|^2}{1-|\kappa|^2}\frac{|(w_{\kappa})_z
|^2}{(1-|w_{\kappa}|^2)^2} d^2 z.
\end{align*}
According to Lemma \ref{appr}, for every $\vep>0$ there exists
$0<\delta<1$ such that for all $\kappa\in U_0$ satisfying
$\Vert\kappa \Vert_{\infty} < \delta$ we have
\begin{align*}
\Bigl|g_{\mu\bar{\mu}}(\kappa)-g_{\mu\bar{\mu}}(0) \Bigr| &\leq
4\iint\limits_{\Del^*}
\frac{|\mu|^2}{1-|\kappa|^2}\left\vert\frac{|(w_{\kappa})_z
|^2}{(1-|w_{\kappa}|^2)^2} -\frac{1}{(1-|z|^2)^2}\right\vert \\
&\hspace{2cm}+
\frac{|\mu|^2}{(1-|z|^2)^2}\left(\frac{1}{1-|\kappa|^2}-1\right)
d^2z\\ &\leq \frac{\vep+ \delta^2}{1-\delta^2} \iint\limits_{\Del}
|\mu|^2 \rho d^2z.
\end{align*}
Thus, for $\delta$ small enough
$\Bigl|g_{\mu\bar{\mu}}(\kappa)-g_{\mu\bar{\mu}}(0) \Bigr|\leq
2\vep\Vert \mu\Vert_2^2.$
\end{proof}
\begin{remark}
Using the basic properties of the q.c.~mappings, it can be shown
that the Weil-Petersson metric is real-analytic on $T(1)$. In
fact, it is sufficient to prove that for every $\mu, \nu\in
H^{-1.1}(\Del^*)$ the mapping $V_0\ni\kappa \mapsto g_{\mu
\bar{\nu}}(\kappa)\in\CC$ is real-analytic on $V_0$. Since this
result will not be used later, we omit the proof. Explicit
curvature computations in Section 7 will show that the
Weil-Petersson metric on $T(1)$ is twice differentiable.
\end{remark}
We will prove in Section 7 that the Weil-Petersson metric is \Ka.
Its symplectic form $\omega_{WP}$ is a right-invariant
$(1,1)$ form on the Hilbert manifold $T(1)$. At the origin of
$T(1)$,
\begin{align*}
\omega_{WP}(\mu,\bar{\nu})=\tfrac{i}{2} \langle \mu, \nu
\rangle_{WP},\;\;\mu,\nu\in H^{-1,1}(\Del^*).
\end{align*}
\begin{remark}
The Weil-Petersson metric on the distribution $\mathfrak{D}_T$
(without defining the Hilbert manifold structure) was introduced
by S. Nag and A. Verjovsky \cite{NV} as a direct generalization of
the Weil-Petersson metric on the finite-dimensional Teichm\"uller
spaces. It was proved in \cite{NV} that the embedding
$\Mob(S^1)\bk\Diff_+(S^1) \hookrightarrow T(1)$ is holomorphic and
the pull-back of the Weil-Petersson metric on the distribution
$\mathfrak{D}_T$ coincides, up to a constant, with the right invariant
\Ka metric introduced by Kirillov \cite{Ki} by the orbit method. At the
tangent space of the origin the latter metric is defined by (cf.
\eqref{VK-metric})
\begin{align*}
\left\Vert \text{v}\right\Vert^2 =\sum_{n=2}^{\infty} (n^3-n)
|c_n|^2,
\end{align*}
where
\begin{displaymath}
\text{v}=\dfrac{\text{u} - i J\,\text{u}}{2}\quad\text{and}\quad
\text{u}=\sum_{n\in\ZZ\setminus\{0,\pm 1\}} c_n
e^{in\theta}\frac{d}{d\theta}\in T_0^{\RR}
\Mob(S^1)\backslash\Diff_+(S^1).
\end{displaymath}
\end{remark}

\section{Characteristic forms of the universal \Te curve}
Let $V=T_v\mathcal{T}(1)$ be the vertical tangent bundle of the
fibration
\begin{displaymath}
\pi: \mathcal{T}(1) \rightarrow T(1).
\end{displaymath}
It is a holomorphic line bundle over the complex Banach manifold
$\mathcal{T}(1)$, the fiber over a point $([\mu],z)\in
\mathcal{T}(1)$ is the holomorphic tangent bundle to the
quasi-disk $w^\mu(\Del^*)$. The hyperbolic metric on
$w^\mu(\Del^*)$ defines a Hermitian metric on $V$, and we denote
by $c_1(V)$ the first Chern form of $V$ corresponding to this
metric.

The Hilbert manifold structure on $T(1)$ naturally induces a
Hilbert manifold structure on $\mathcal{T}(1)$, such that the
projection $\pi: \mathcal{T}(1) \rightarrow T(1)$ is a holomorphic
mapping of Hilbert manifolds. Similar to the Hilbert manifold
$T(1)$, the Hilbert manifold $\mathcal{T}(1)$ is also a disjoint
union of uncountably many components. We will prove in Appendix A
that the component containing the identity $\mathcal{T}_0(1)$ is a
topological group.

 The vertical tangent bundle is also a
holomorphic line bundle over the Hilbert manifold
$\mathcal{T}(1)$. We will continue to denote corresponding line
bundle by $V$, and by $c_1(V)$ --- the first Chern form
corresponding to the hyperbolic metric on the fibers, specifying
explicitly which topology we are using. Since the Hilbert manifold
topology is stronger than the Banach manifold topology, the form
$c_1(V)$ for the Banach manifold structure on $\mathcal{T}(1)$
naturally restricts onto $\mathcal{T}(1)$ considered as a Hilbert
manifold.

Similar to Wolpert's work \cite{Wol} on finite dimensional
Teichm\"uller spaces, we define the analogs of
Mumford-Morita-Miller characteristic forms as the following
$(n,n)$-forms on the Hilbert manifold $T(1)$,
\begin{equation}\label{mumford}
\kappa_n = (-1)^{n+1}\pi_*(c_1(V)^{n+1}),
\end{equation}
where
$\pi_*:\Omega^*(\mathcal{T}(1))\rightarrow\Omega^{*-2}(T(1))$ is
the operation of ``integration over the fibers'' of the fibration
$\pi: \mathcal{T}(1) \rightarrow T(1)$. As we will see in Section
\ref{fiber-integration}, it is the passage from the Banach manifold
structure to the Hilbert manifold structure that makes the
operation $\pi_*$ well-defined (i.e., the integrals over
non-compact fibers become convergent).
\subsection{The form $c_1(V)$ as Velling-Kirillov symplectic form}
In this section we work with the Banach manifold structure on
$\mathcal{T}(1)$. Let $z$ be the complex coordinate on
$\hat{\C}\setminus\{0\}$. The assignment
$\mathcal{T}(1)\ni([\mu],z)\mapsto -z^2 \pa_z$ defines a
holomorphic section of the line bundle $V$ over $\mathcal{T}(1)$
\footnote{Under the conformal map $z \mapsto \frac{1}{z}$ the
vector field $-z^2 \pa_z\mapsto\pa_z$.}. The hyperbolic metric on
$w^\mu(\Del^*)$ is the pull-back of the hyperbolic metric on
$\Del^*$ by the conformal map $g_\mu^{-1}$, so that
\begin{equation*}
\left\Vert z^2\pa_z\right\Vert^2_{([\mu],z)} =
 \frac{4 \vert z^2 (g_{\mu}^{-1})'(z)\vert^2}
{ (\vert g_{\mu}^{-1}(z)\vert^2 -1)^2}.
\end{equation*}
The first Chern form of the line bundle $V$ is
\begin{equation*}
c_1(V) = \frac{i}{2\pi}\Theta =\frac{i}{2\pi}\bar{\pa}\pa \log
\left\Vert z^2\pa_z\right\Vert^2,
\end{equation*}
where $\pa$ and $\bar{\pa}$ are, respectively, the $(1,0)$ and
$(0,1)$ components of the de Rham differential on
$\mathcal{T}(1)$.

Let
\begin{displaymath}
K=\log \left\Vert z^2 \pa_z \right\Vert_{([\mu],z)} -
\,\log 2.
\end{displaymath}
\begin{lemma}\label{functionK} The function $K:\mathcal{T}(1)\rightarrow\RR$
is real-analytic. Under the correspondence
$\mathcal{T}(1)\ni([\mu],z)\mapsto\g\in S^1\backslash
\text{Homeo}_{qs}(S^1)$, where $\gamma = \left.(g^{-1}\circ
f)\right|_{S^1}$,
\begin{equation*}
K(\gamma) = \log |g'(\infty)|.
\end{equation*}
\end{lemma}
\begin{proof}
Using the formulas $g = \lambda_w \circ g_{\mu} \circ
\sigma_w^{-1}$ and $w=g_\mu^{-1}(z)$ from Section
\ref{conformalwelding}, it is straightforward to compute
\begin{displaymath}
g'(\infty) = \frac{z^2 (g_{\mu}^{-1})^\prime(z)(1-\bar{w})}{(|w|^2
-1)(1-w)}.
\end{displaymath}
Now it easily follows from the general properties of q.c.~mappings
that for $z\in\Del^*$ the functional $T(1)\ni[\mu]\mapsto
g_\mu^{-1}(z)\in\CC$ is real-analytic so that $|g^\prime(\infty)|$
is a real-analytic function on $\mathcal{T}(1)$.
\end{proof}
\begin{remark}
The quantity $|g^\prime(\infty)|$ is the capacity of the
quasi-circle $g(S^1)$ corresponding to
$\gamma\in\text{Homeo}_{qs}(S^1)$.
\end{remark}
\begin{theorem} \label{Chern-V-K}
The first Chern form of the vertical tangent bundle to the
universal \Te curve $\mathcal{T}(1)$ is proportional to the
symplectic form of the Velling-Kirillov metric,
\begin{displaymath}
c_1(V) =-\tfrac{2}{\pi}\,\omega_{VK}.
\end{displaymath}
Equivalently, the function $K$ is a \Ka potential for the
Velling-Kirillov metric on $\mathcal{T}(1)$.
\end{theorem}
\begin{proof} It is based on the following lemmas.
\begin{lemma} \label{right-invariance} The $(1,1)$-form $c_1(V)$ on $\mathcal{T}(1)$ is
right-invariant.
\end{lemma}
\begin{proof}
We need to prove that for every $\g_0\in
S^1\backslash\text{Homeo}(S^1)\cong \mathcal{T}(1)$ the difference
$R_{\gamma_0}^*K -K$, where $R_{\gamma_0}:
\mathcal{T}(1)\rightarrow\mathcal{T}(1)$ is the right translation
by $\gamma_0$ and
$(R_{\gamma_0}^*K)(\gamma)=K(\gamma\circ\gamma_0)$, is a harmonic
function on $\mathcal{T}(1)$.

For every $\gamma = g^{-1}\circ f\in\mathcal{T}(1)$ let
$\tilde{\gamma}= \gamma \circ \gamma_0=\tilde{g}^{-1}\circ
\tilde{f} $. Since $\tilde{g}=\tilde{f}\circ\g^{-1}_0\circ
f^{-1}\circ g$, we have
\begin{displaymath}
\left( R_{\gamma_0}^*K -K \right)(\gamma) = \log \vert
\tilde{g}'(\infty)\vert -\log \vert g'(\infty)\vert= \log\vert
(\tilde{f}\circ\gamma_0^{-1}\circ f^{-1})'(\infty)\vert.
\end{displaymath}
In \cite{Bers1}, Bers has proved that the function
\begin{align*}
(\gamma, z) \mapsto \left(\tilde{f} \circ \gamma_0^{-1} \circ
f^{-1}\right)(z)= h(\gamma ,z)
\end{align*}
depends holomorphically on $\gamma$ and $z$, which
implies that $\left(\tilde{f} \circ \gamma_0^{-1} \circ
f^{-1}\right)'(\infty)$ depends holomorphically on $\gamma$ and
our assertion follows.
\end{proof}
\begin{lemma}
Let $\gamma = g^{-1} \circ f \in \mathcal{T}(1)$, where
$f\vert_{\Del}(z) = \sum_{n=0}^{\infty} a_n z^{n+1},~a_0 =1$, and
$g\vert_{\Del^*}(z) = \sum_{n=0}^{\infty} b_n z^{1-n}$. Then
\begin{equation} \label{boterm}
|b_0|^2 = \sum_{n=0}^{\infty} (n+1) |a_n|^2 +
\sum_{n=1}^{\infty}(n-1) |b_n|^2.
\end{equation}
\end{lemma}
\begin{proof}
Evaluate the Euclidean area of the domain $\Omega = f(\Del) =
g(\Del)$ in two different ways. First,
\begin{align*}
A_E(\Omega) &= \lim_{r \rightarrow 1^-} \iint\limits_{f(\Del_r)}
d^2z
   = \lim_{r \rightarrow 1^-} \iint\limits_{\Del_r} |f'|^2 d^2z
   = \pi \sum_{n=0}^{\infty} (n+1) |a_n|^2,
\end{align*}
where $\Del_r =\{ z\in\CC : |z| < r \}$. On the other hand, the classical
area theorem gives
\begin{displaymath}
A_E(\Omega) = \pi \sum_{n=0}^{\infty} (1-n) |b_n|^2,
\end{displaymath}
and we obtain \eqref{boterm}.
\end{proof}
Now we complete the proof of the theorem. Let
\begin{displaymath}
\text{u} = \sum_{n\in\ZZ\setminus\{0\}}
c_n e^{in\theta}\frac{d}{d\theta} \in T^{\RR}_0 \mathcal{T}(1),
\end{displaymath}
and let $\gamma_t = g_t^{-1} \circ f^t$, $\gamma_0 = \id$, be the
corresponding smooth curve in $\mathcal{T}(1)$. Using notations
from the previous lemma, differentiate the relation \eqref{boterm}
with respect to $t$ and set $t=0$, and using
$b_0(0)=1,\,a_n(0)=b_n(0)=0$ for $n\geq 1$, we get
\begin{equation} \label{b-0}
\dot{b}_0 +\bar{\dot{b}}_0=0.
\end{equation}
Here we denote
\begin{displaymath}
\dot{a}_n(\text{u})= \dot{a}_n = \frac{d a_n}{dt}(0) \quad
\text{and} \quad \dot{b}_n(\text{u}) =  \dot{b}_n = \frac{d
b_n}{dt}(0).
\end{displaymath} Differentiating \eqref{boterm} twice with respect
to $t$ and setting $t=0$ we get
\begin{align*}
\frac{d^2b_0}{dt^2}(0) + \frac{d^2\bar{b}_0}{dt^2}(0)
+ 2 \left\vert \frac{d b_0}{dt}(0) \right\vert^2
&= 2\sum_{n=1}^{\infty} (n+1)
|\dot{a}_n|^2 + 2\sum_{n=1}^{\infty} (n-1) |\dot{b}_n|^2 \\
 &= 4\sum_{n=1}^{\infty} n |\dot{a}_n|^2,\\
\end{align*}
where we have also used the property \textbf{TS1} in Section \ref{TS}.
Since $g_t^\prime(\infty) = b_0(t)$, using \eqref{b-0} we get
\begin{align*}
\frac{d^2}{dt^2} \log |g_t^\prime(\infty)| \bigl\vert_{t=0}
&=2\sum_{n=1}^{\infty} n |\dot{a}_n|^2.
\end{align*}

Let $\text{v} =\tfrac{1}{2}(\text{u} - iJ\,\text{u})$ be the holomorphic
tangent vector to $\mathcal{T}(1)$. Since
$\dot{a}_n(J\,\text{u})=i\dot{a}_n(\text{u})$, using \eqref{de-Rham}
we finally get
\begin{equation} \label{K}
(\pa\bar\pa K)(\text{v},\bar{\text{v}})= \frac{1}{2}\sum_{n=1}^\infty
n\left(|\dot{a}_n(\text{u})|^2 + |\dot{a}_n(J\,\text{u})|^2\right)=
\sum_{n=1}^\infty n|\dot{a}_n(\text{u})|^2.
\end{equation}
This proves that $\Theta=4i\omega_{VK}$ at the origin of
$\mathcal{T}(1)$. Since both these $(1,1)$-forms on $\mathcal{T}(1)$ are
right-invariant, the assertion follows.
\end{proof}
\begin{remark}
In \cite{KY}, A.A.~Kirillov and D.~Yuriev have stated that the
function $K$, restricted to the space $S^1\bk\Diff_+(S^1)$, is a
K\"ahler potential of the Velling-Kirillov metric. Theorem
\ref{Chern-V-K} extends this result to $\mathcal{T}(1) \simeq
S^1\bk\text{Homeo}_{qs}(S^1)$ and gives its geometric
interpretation.
\end{remark}
\subsection{The Chern form $c_1(V)$ and the resolvent kernel}
Let $\mu\in\Omega^{-1,1}(\Del^*)$ be a horizontal holomorphic
tangent vector to $\mathcal{T}(1)$ at the origin. According to the
property \textbf{TV1} in Section \ref{h-v-subspaces}, the vector field
$\tau_\mu$ --- the horizontal lift of the vector field
$\tfrac{\pa}{\pa\varepsilon_\mu}$ on $U_0 \subset T(1)$ to the
point $(0,z)\in\pi^{-1}(0)$, is identified with $(\sigma_z^{-1})^*
\mu\in\Omega^{-1,1}(\Del^*)$.
\begin{proposition} \label{V-K-onfiber}
\begin{itemize}
\item[(i)] On the fiber $\pi^{-1}(0)\subset \mathcal{T}(1)$ the
Velling-Kirillov metric is given by
\begin{align*}
\left\la \pa_z, \pa_z \right\ra_{VK} (0,z)
&= \frac{1}{(1-|z|^2)^2}, \\ \left\la \pa_z,\tau_\mu \right\ra_{VK}
 (0,z)& = 0, \\
\left\la \tau_\mu, \tau_\mu
\right \ra_{VK} (0,z) & = \frac{1}{2}\iint\limits_{\Del^*} G(z,u)
|\mu(u)|^2 \rho(u)d^2 u.
\end{align*}
\item[(ii)] On the fiber $\pi^{-1}(0)\subset\mathcal{T}(1)$ the
$(1,1)$-form $\Theta$ is given by
\begin{align*}
\Theta \left( \pa_z, \pa_{\bar{z}}\right)(0,z) &= -\frac{2}{(1-|z|^2)^2}, \\
\Theta\left(\pa_z, \tau_{\bar{\mu}} \right)
 (0,z)& = 0, \\
\Theta\left( \tau_\mu, \tau_{\bar{\mu}}\right) (0,z) &= -\iint\limits_{\Del^*} G(z,u)|\mu(u)|^2 \rho(u)d^2 u.
\end{align*}
\item[(iii)] The vertical holomorphic tangent bundle $V\rightarrow
\mathcal{T}(1)$ of the fibration $\pi:\mathcal{T}(1)\rightarrow
T(1)$ is a negative line bundle.
\end{itemize}
\end{proposition}
\begin{proof}
It follows from the property \textbf{TV2} in Section
\ref{h-v-subspaces} that the vector field $\pa_z$ at
$(0,z)\in\pi^{-1}(0)$ corresponds to the tangent vector
\begin{displaymath}
\text{u}=\sum_{n\in\ZZ\setminus\{0\}}c_n
e^{in\theta}\dfrac{d}{d\theta}\in T_0^{\RR}
S^1\backslash\text{Homeo}_{qs}(S^1)
\end{displaymath}
with $c_1=\tfrac{1-\z}{(1-z)(1-|z|^2)}$ and $c_n=0$ for $n\geq 2$.
This proves the first formula in part (i). The second formula
follows from the fact that, according to the property \textbf{TS1}
in Section \ref{TS}, the tangent vector $\text{u}\in T_0^{\RR}
S^1\backslash\text{Homeo}_{qs}(S^1)$ which corresponds to the
horizontal lift $\tau_\mu$ of the vector field
$\tfrac{\pa}{\pa\varepsilon_\mu}$ to $(0,z)\in\pi^{-1}(0)$, has
$c_1=0$. The last formula follows from the following lemma.
\begin{lemma} Let $\mu \in \Omega^{-1,1}(\Del^*)$
and
\begin{displaymath}
\phi(z)=D_0\beta(\mu)(z) =  \sum_{n=2}^{\infty} (n^3 -n) a_n z^{n-2}
\in A_{\infty}(\Del).
\end{displaymath}
Then
\begin{align*}
\iint\limits_{\Del^*} G(z,u) |\mu(u)|^2\rho(u) d^2
u= 2 \sum_{n=2}^{\infty} n |a_n^z|^2,
\end{align*}
where
\begin{displaymath}
(\sigma_z^{-1})^*(\phi)(u) =\sum_{n=2}^{\infty} (n^3 -n) a_n^z u^{n-2}
\end{displaymath}
is the power series expansion of $(\sigma_z^{-1})^*(\phi) =
\phi\circ\sigma_z^{-1}\left((\sigma_z^{-1})^\prime\right)^2\in
A_\infty(\Del)$.
\end{lemma}
\begin{proof}
Since $G(z,u)$ is a point-pair invariant, it is sufficient to
prove the formula for $z=\infty$. In this case, using the relation
$G(\infty,u)=G(0,1/\bar{u})$ between the resolvent kernels on
$\Del^*$ and $\Del$ and the formula $\mu =\Lambda(\phi)$ we get
\begin{align*}
\iint\limits_{\Del^*} G(\infty,u) |\mu(u)|^2\rho(u) d^2u
  = \iint\limits_{\Del} G(0,u) (1-|u|^2)^2 |\phi(u)|^2 d^2u.
\end{align*}
It follows from the explicit formula \eqref{resolvent-kernel} that
\begin{displaymath}
(1-|u|^2)^2 G(0,u) = h(r) = \left(\frac{1}{2 \pi} \frac{1+r^2}{1-r^2} \log
\frac{1}{r^2} - \frac{1}{\pi}\right) (1-r^2)^2,\;\;r=|u|.
\end{displaymath}
Since $h$ is an integrable function on $[0,1)$, we have
\begin{align*}
\iint\limits_{\Del} G(0,u) (1-|u|^2)^2 |\phi(u)|^2 d^2u = 2 \pi
\sum_{n=2}^{\infty} (n^3-n)^2 |a_n|^2 \int_0^1 h(r) r^{2n-4} r dr.
\end{align*}
A straightforward computation gives
\begin{displaymath}
 \int_0^1 h(r) r^{2n-4} r dr = \frac{1}{\pi}
 \frac{1}{(n^3-n)(n^2-1)},
\end{displaymath}
which proves the lemma.
\end{proof}
Using this Lemma, the property \textbf{TV1} in Section
\ref{h-v-subspaces} and the property \textbf{TS2} in Section
\ref{TS}, we get the third formula in part (i). Now part (ii)
follows from Theorem \ref{Chern-V-K}, and part (iii) follows from
part (ii) and the property \textbf{RK2} in Section
\ref{resolvent}.
\end{proof}
\begin{remark} Parts (ii) and (iii) of Proposition
\ref{V-K-onfiber} generalize Wolpert's computation of the
$(1,1)$-form $\Theta$ for finite-dimensional \Te spaces (see
Theorem 5.5 and formula (5.3) in \cite{Wol}).
\end{remark}
\subsection{Mumford-Morita-Miller characteristic forms} \label{fiber-integration}
In this section we consider $\pi:\mathcal{T}(1)\rightarrow T(1)$
as a holomorphic fibration of Hilbert manifolds and evaluate the
Mumford-Morita-Miller forms $\kappa_n$ on $T(1)$.
\begin{theorem}\label{Mumford} The characteristic forms $\kappa_n$
are right-invariant on the Hilbert manifold $T(1)$ and for
$\mu_1,\dots,\mu_n,\nu_1,\dots,\nu_n\in H^{-1,1}(\Del^*)\simeq T_0
T(1)$
\begin{align*}
&\kappa_n(\mu_1,\dots,\mu_n,\bar\nu_1,\dots,\bar\nu_n)\\
=&\frac{i^n(n+1)!}{(2\pi)^{n+1}}\sum_{\sigma\in S_n}
sgn(\sigma)\iint\limits_{\Del^*}
 G\left(\mu_1 \bar{\nu}_{\sigma(1)}\right)\ldots G\left(\mu_n
\bar{\nu}_{\sigma(n)}\right)\rho(z)d^2z,
\end{align*}
where the sum goes over the permutation group $S_n$ on $n$
elements and $sgn(\sigma)$ is the sign of the permutation
$\sigma$.
\end{theorem}
\begin{proof}
It is straightforward computation of the integral
\begin{gather*}
\kappa_n(\mu_1,\dots,\mu_n,\bar\nu_1,\dots,\bar\nu_n) \\ =
\left(\frac{-i}{2\pi}\right)^{n+1}\iint\limits_{\Del^*}
\Theta^{n+1}\left(\pa_z,\pa_{\z},
\tau_{\mu_1},\tau_{\bar{\nu}_1},\dots,
\tau_{\mu_n},\tau_{\bar{\nu}_n}\right)dz\wedge
d\z
\end{gather*}
using Part (ii) of Proposition \ref{V-K-onfiber}. We need only to
verify that the integral is convergent. This follows from the
properties of the resolvent kernel in Section \ref{resolvent}.
Indeed, the property \textbf{RK3} assures that $G(\mu\bar{\nu})$
is bounded on $\Del^*$ for $\mu,\nu\in\Omega^{-1,1}(\Del^*)$, and
properties \textbf{RK2} and \textbf{RK4} imply that for
$\mu,\nu\in H^{-1,1}(\Del^*)$,
\begin{align*}
\left|\iint\limits_{\Del^*}G(\mu\bar{\nu}) \rho(z)d^2z\right| &
\leq \iint\limits_{\Del^*}\iint\limits_{ \Del^*}G(z,u)|\mu(u)\nu(u)|\rho(u)\rho(z)
d^2ud^2z\\ &= \iint\limits_{\Del^*}\iint\limits_{ \Del^*}G(z,u)|\mu(u)\nu(u)|\rho(z)\rho(u)
d^2zd^2u\\
&=\iint\limits_{\Del^*}|\mu(u)\nu(u)|\rho(u)d^2u<\infty.
\end{align*}
\end{proof}
\begin{corollary} \label{K-1-WP}
On the Hilbert manifold $T(1)$
\begin{displaymath}
\kappa_1= \tfrac{1}{\pi^2}\,\omega_{WP}.
\end{displaymath}
\end{corollary}
\begin{proof}
We have, using again the property \textbf{RK4} in Section
\ref{resolvent},
\begin{align*}
\kappa_1(\mu, \bar{\nu})=&
\frac{i}{2\pi^2}\iint\limits_{\Del^*}G(\mu\bar{\nu})
\rho(z)d^2z =\frac{i}{2\pi^2}\iint\limits_{\Del^*}\mu(z)
\ov{\nu(z)}\rho(z)d^2z\\ = & \frac{1}{\pi^2}\,\omega_{WP}(\mu,
\bar{\nu}).
\end{align*}
\end{proof}
\begin{remark} Combining Corollary \ref{K-1-WP}, Part (ii) of
Proposition \ref{V-K-onfiber} and Theorem \ref{Chern-V-K}, we get
another proof of Theorem 4.3 in \cite{Teo}.
\end{remark}
\begin{remark}
Theorem \ref{Mumford} generalizes Wolpert's result for
finite-dimensional \Te spaces (see Lemma 5.9 and Lemma 5.10 in
\cite{Wol}) to the universal \Te space.
\end{remark}
\section{First and second variations of the hyperbolic metric}\label{VM}
Here we present a concise formula for the second variation of the
hyperbolic metric in terms of the resolvent kernel. We are using
the model $\HH^2\simeq \UU$, so that the density of the hyperbolic
metric is the  $(1,1)$ -- tensor $\rho(z) =y^{-2}$ on $\UU$.
\subsection{The first variation}
It is a classical result of Ahlfors \cite{Ahl1} that the first
variation of the hyperbolic metric at the origin of $T(1)$ is
identically zero.
\begin{lemma}\label{Ahlfors2} For every $\mu\in
\Omega^{-1,1}(\U)$,
\begin{displaymath}
L_\mu\rho=0.
\end{displaymath}
\end{lemma}
\begin{proof}
Since
\begin{align*}
\rho^{\vep \mu}= w_{\vep\mu}^*(\rho) = -\,\frac{4\left\vert
(w_{\vep\mu})_z(z)\right\vert^2}{(w_{\vep\mu}(z) -
\ov{w_{\vep\mu}(z)})^2},
\end{align*}
we have
\begin{align*}
L_\mu\rho(z) = \frac{\pa}{\pa
\vep}\Bigr\vert_{\vep=0} \rho^{\vep \mu}(z) = -4\frac{F_z(z) +
\ov{\Phi^{\prime}(z)}}{(z-\z)^2} + 8\frac{F(z) -\ov{\Phi(z)}
}{(z-\z)^3},
\end{align*}
where $F=F[\mu], \Phi=\Phi[\mu]$, and the result follows from
Lemma \ref{Ahlfors}.
\end{proof}
\begin{remark}
For the case $\mu\in \Omega^{-1,1}(\U,\Gamma)$, where $\Gamma$ is
a cofinite Fuchsian group, another proof of the Ahlfors result was
given by Wolpert \cite{Wol}.
\end{remark}
\subsection{The second variation}
Set
\begin{displaymath}
L_\mu L_{\bar{\mu}}\rho = \frac{\pa^2}{\pa
\vep \pa \bar{\vep} }\Bigr\vert_{\vep=0}\,\rho^{\vep\mu}.
\end{displaymath}
We have
\begin{proposition}\label{secondvarmetric}
For every  $\mu\in \Omega^{-1,1}(\U)$,
\begin{align*}
L_\mu L_{\bar{\mu}} \rho = \rho G(|\mu|^2).
\end{align*}
\end{proposition}
\begin{proof} Using the representation
\begin{equation} \label{metric-pullback}
\rho^{\vep\mu}(z) = - 4K_{\vep\mu}(z,\z)
\end{equation}
and the first formula in \eqref{ahlfors-K-bar} we get
\begin{align}\label{firstvar}
\frac{\pa}{\pa \vep}\rho^{\vep\mu}(z) = &\frac{4}{\pi}
\iint\limits_{\U}\mu(u)K_{\vep\mu}(z,u)K_{\vep\mu}(u,\z)d^2u,
\end{align}
where the integral is understood in the principal value sense.
Setting $\vep=0$ in \eqref{firstvar} and using Lemma \ref{Ahlfors2}, we obtain
\begin{equation} \label{integral-zero}
\iint\limits_{\U} \frac{\mu(u)}{(u-z)^2(u-\z)^2} d^2u
=0\;\;\text{for all}\;\; z\in\U.
\end{equation}
Using formulas \eqref{ahlfors-K} and
\eqref{ahlfors-K-bar}, we get from \eqref{firstvar} the following
integral representation for the second variation of the hyperbolic
metric
\begin{align} \label{sing-integral}
L_\mu L_{\bar{\mu}}\rho(z)=&-\frac{4}{\pi^2} \iint\limits_{\U}\iint\limits_{\U}\mu(u)
\ov{\mu(v)}\left(\frac{1}{(u-\bar{v})^2(z-\bar{v})^2(u-\z)^2}\right. \nonumber
\\ &\hspace{2cm}+
\left.\frac{1}{(u-z)^2(u-\bar{v})^2(\z-\bar{v})^2}\right)d^2ud^2v.
\end{align}
The differentiation under the integral sign is justified by the
same argument as in \cite{Ahlfors3}. We transform the principal
value integrals in \eqref{sing-integral} into the ordinary
integrals by using the identity
\begin{equation} \label{identity}
\frac{1}{u-\bar{v}} = \frac{z-\z}{(u-\z)(z-\bar{v})} +
\frac{(\bar{v}-\z)(z-u)}{(u-\z)(z-\bar{v})(u-\bar{v})},
\end{equation}
which gives
\begin{align*}
&\frac{1}{(u-z)^2(u-\bar{v})^2(\z- \bar{v})^2} \\ =&\,
\frac{1}{(u-\z)^2(z-\bar{v})^2(u- \bar{v})^2}\, +
\,\frac{2(z-\z)}{(u-\z)^3(z-\bar{v})^3(u-\bar{v})} \\ +&
\,\frac{(z-\z)^2}{(u-z)^2(u -\z)^2 (z-\bar{v})^2(\z -\bar{v})^2}\,
+\,\frac{2(z-\z)^2}{(u-z)(u-\z)^3 (z-\bar{v})^3(\z -\bar{v})}.\\
\end{align*}
Using \eqref{integral-zero} and Corollary \ref{integral-1} we see
that last two terms in this formula do not contribute into the
representation \eqref{sing-integral} and we obtain
\begin{gather*}
L_\mu L_{\bar{\mu}}\rho(z) =
-\frac{8}{\pi^2} \iint\limits_{\U}\iint\limits_{\U}\mu(u)
\ov{\mu(v)}\left(\frac{1}{(u-\bar{v})^2(z-\bar{v})^2(u-\z)^2}\right.
\\ + \left.\frac{(z-\z)}{(u-\z)^3(u-\bar{v})(z-\bar{v})^3}\right)d^2ud^2v.
\end{gather*}
Now we apply the operator $2(\Delta_0 + \tfrac{1}{2})$ to the
bounded function $\rho^{-1}L_\mu L_{\bar{\mu}}\rho$ on $\U$. Using \eqref{identity}
it is straightforward to compute that
\begin{align*}
\left(2\Delta_0+ 1\right)\left(\frac{(z-\z)^2}
{(u-\bar{v})^2(z-\bar{v})^2(u-\z)^2}
+\frac{(z-\z)^3}{(u-\z)^3(u-\bar{v})(z-\bar{v})^3}\right)\\ =
\frac{9}{2} \frac{(z-\z)^4}{(u-\z)^4(z-\bar{v})^4},
\end{align*}
which, together with \eqref{projection-1}, gives
\begin{gather*}
\left(2\Delta_0+1\right)\left(\rho^{-1}L_\mu L_{\bar{\mu}}\rho\right)(z)\\ =
\frac{9}{\pi^2} \iint\limits_{\U}\iint\limits_{\U}
\mu(u)\ov{\mu(v)} \frac{(z-\z)^4}{(u-\z)^4(z-\bar{v})^4}d^2ud^2v =
|\mu(z)|^2.
\end{gather*}
Using the property \textbf{RK3} in Section \ref{resolvent}
completes the proof.
\end{proof}
\begin{corollary} \label{G-formula} For every
$\mu,\nu\in\Omega^{-1,1}(\Del^*)$,
\begin{align*}
G(\mu\bar{\nu})(z)& =
\frac{2}{\pi^2} \iint\limits_{\U}\iint\limits_{\U}\mu(u)
\ov{\nu(v)}\left(\frac{(z-\z)^2}{(u-\bar{v})^2(z-\bar{v})^2(u-\z)^2}\right.
\\& + \left.\frac{(z-\z)^3}{(u-\z)^3(u-\bar{v})(z-\bar{v})^3}\right)d^2ud^2v.
\end{align*}
\end{corollary}
\begin{remark} It follows from Proposition \ref{secondvarmetric} by polarization that
\begin{equation*}
L_{\mu}L_{\bar{\nu}}\rho =\left.\frac{\pa^2}{\pa\vep_1\pa\bar{\vep}_2}
\right\vert_{\vep_1=\vep_2=0}\rho^{\vep_1\mu +\vep_2\nu} = \rho G(\mu\bar{\nu}).
\end{equation*}
\end{remark}
\begin{remark} For the case $\mu\in\Omega^{-1,1}(\U,\Gamma)$,
where $\Gamma$ is a cofinite Fuchsian group, the formula for the
second variation of the hyperbolic metric in Proposition
\ref{secondvarmetric} was first proved by Wolpert \cite{Wol}.
However, the method in \cite{Wol} does not work for the universal
Teichm\"uller space $T(1)$. The proof of Proposition
\ref{secondvarmetric} shows that the  original singular integral
representation \eqref{sing-integral} of Ahlfors can be easily
transformed to a closed form using the resolvent kernel.
\end{remark}
\section{Riemann curvature tensor} \label{Riemann-tensor}
In this section we consider $T(1)$ as a Hilbert manifold equipped
with the Weil-Petersson metric. We prove that the Weil-Petersson
metric is \Ka, compute its Riemann and Ricci tensors, and show
that the Ricci, holomorphic, and  sectional curvatures are all
negative. Since the Weil-Petersson metric is right-invariant, it
is sufficient to compute these tensors at the origin of $T(1)$.
\subsection{The first variation of the Weil-Petersson metric}
For $\mu,\nu\in\Omega^{-1,1}(\Del^*)$ set
\begin{align*}
Q(\mu, \nu)= P(R(\mu,\nu))\circ w_{\nu}\,
\frac{(\ov{w_{\nu}})_{\z}}{(w_{\nu})_z}.
\end{align*}
\begin{proposition} \label{varQ1}
For $\mu , \nu \in \Omega^{-1,1}(\Del^*)$,
\begin{align*}
&\frac{\pa}{\pa \vep} Q(\mu, \vep\nu)
\Bigr\vert_{\vep=0}=0,\\ &\frac{\pa}{\pa \bar{\vep}}
Q(\mu, \vep\nu) \Bigr\vert_{\vep=0} =-2\frac{\pa}{\pa
\z}\rho^{-1}\frac{\pa}{\pa \z} G(\mu\bar{\nu}).
\end{align*}
\end{proposition}
\begin{proof} We will be using canonical complex anti-linear
isomorphism $\Omega^{-1,1}(\Del^*)\simeq\Omega^{-1,1}(\Del)$,
given by the reflection \eqref{sym}, and the model
$\HH^2\simeq\UU$ of the hyperbolic plane. From \eqref{R} we get
\begin{equation} \label{Q}
\rho^{\vep\nu}(z)\, Q(\mu,\vep\nu)(z) = \frac{12}{\pi} \iint\limits_{\U}
\mu(u) K_{\vep\nu}(u,\bar{z})^2 d^2u.
\end{equation}
It follows from equations \eqref{ahlfors-K-bar} that
\begin{equation}\label{varQ-e}
\frac{\pa}{\pa\vep}\Bigr\vert_{\vep=0} Q(\mu, \vep\nu)(z)
=\frac{6(z-\z)^2}{\pi^2}\iint\limits_\U\iint\limits_\U
\frac{\mu(u)\nu(v)}{(u-\z)^2(u-v)^2(v-\z)^2}d^2ud^2v
\end{equation}
and
\begin{equation} \label{varQ-be}
\left.\frac{\pa}{\pa\bar{\vep}}\right|_{\vep=0} Q(\mu, \vep\nu)(z)=
\frac{6(z-\z)^2}{\pi^2}\iint\limits_\U\iint\limits_\U\frac{\mu(u)\ov{\nu(v)}}{(u-
\z)^2(u-\bar{v})^2(\bar{v}-\z)^2}d^2ud^2v.
\end{equation}
The integrals are understood in the principal value sense and
differentiation under the integral sign in \eqref{Q} is justified
as in \cite{Ahlfors3}.

To prove that the integral \eqref{varQ-e} is zero, we use the identity
\begin{displaymath}
\frac{1}{(u-v)(u-\z)(v-\z)} =\frac{1}{(u-v)^2}\left(\frac{1}{v-\z}-\frac{1}{u-\z}\right),
\end{displaymath}
and rewrite the integral \eqref{varQ-e} as follows
\begin{gather} \label{4I}
\iint\limits_\U\iint\limits_\U
\frac{\mu(u)\nu(v)}{(u-\z)^2(u-v)^2(v-\z)^2}d^2ud^2v \\ \nonumber
=\lim_{\vep\rightarrow 0}\iiiint\limits_{\U\times\U\setminus\{|u-v|<\vep\}}
\left(\frac{1}{(u-v)^4(v-\z)^2} - \frac{2}{(u-v)^5(v-\z)}\right. \\ \nonumber
+\left.\frac{2}{(u-v)^5(u-\z)} +\frac{1}{(u-v)^4(u-\z)^2}\right)
\mu(u)\nu(v)d^2ud^2v\\ \nonumber
=I_1 + I_2 + I_3 + I_4. \nonumber
\end{gather}
Applying Lemma \eqref{integral-4-5} to the principal value
integrals over $u$ in the terms $I_1$, $I_2$ and $I_4$, we
conclude that these terms vanish. Changing the order of
integrations in $I_3$ (which is legitimate since domain of
integration is invariant under the involution $(u,v)\mapsto
(v,u)$) and applying Lemma  \eqref{integral-4-5} to the integral
over $v$ we conclude that the term $I_3$ also vanishes. This
proves that the holomorphic variation of $Q(\mu,\nu)$ vanishes.

To prove the formula for the antiholomorphic variation, we again
use the identity \eqref{identity}, which gives
\begin{align*}
&\frac{1}{(u-\z)^2(u-\bar{v})^2(\bar{v}-\z)^2} =\frac{(u-z)^2}
{(u-\z)^4(\bar{v}-z)^2(u-\bar{v})^2} \\
& + \frac{(z-\z)^2}{(u-\z)^4(\bar{v}-\z)^2(\bar{v}-z)^2} +
\frac{2(z-\z)^2(z-u)}{(u-\z)^5(z - \bar{v})^3(\bar{v}-\z)} \\
& + \frac{2(z-\z)(z-u)^2}{(u-\z)^5(z - \bar{v})^3(u-\bar{v})}.
\end{align*}
Using formula \eqref{integral-zero} and Corollary \ref{integral-1}
we see that the second and third terms do not contribute to
\eqref{varQ-be}, and we get
\begin{align} \label{varQ-be-final}
&\left.\frac{\pa}{\pa\bar{\vep}}\right|_{\vep=0} Q(\mu, \vep\nu)(z)=
\frac{6(z-\z)^2}{\pi^2}\iint\limits_\U\iint\limits_\U\mu(u)\ov{\nu(v)}\\ &\left(
\frac{(u-z)^2}{(u-\z)^4(\bar{v}-z)^2(u-\bar{v})^2} +
\frac{2(z-\z)(z-u)^2}{(u-\z)^5(z - \bar{v})^3(u-\bar{v})}\right)d^2u^2v. \nonumber
\end{align}
Now applying $\pa_{\z}\rho^{-1}\pa_{\z}$ to the integral
representation in  Corollary \ref{G-formula} (the differentiation
under the integral sign being trivially justified) we get the
formula for the antiholomorphic variation.
\end{proof}
 Set
\begin{displaymath}
\dot{Q}(\mu)[\nu] =\left.\frac{\pa}{\pa \bar{\vep}}\right|_{\vep=0}Q(\mu,
\vep\nu) .
\end{displaymath}
\begin{proposition} \label{var-in-L^2} Let $\mu\in H^{-1,1}(\Del^*)$
 and $\nu\in\Omega^{-1,1}(\Del^*)$. Then $\dot{Q}(\mu)[\nu]\in L^2(\Del^*,\rho(z)d^2z)$ and
\begin{align*}
\left\Vert\dot{Q}(\mu)[\nu]\right\Vert^2_2= \Vert\mu\bar{\nu}\Vert^2_2
- (\mu\bar{\nu}, G(\mu\bar{\nu})),
\end{align*}
where $(~,~)$ stands for the inner product in the Hilbert space $L^2(\Del^*,\rho(z)d^2z)$.
\end{proposition}
\begin{proof}  As in the proof of Proposition \ref{varQ1}, it is
 convenient to use the isomorphism $\Omega^{-1,1}(\Del^*)\simeq\Omega^{-1,1}(\Del)$.
  For  $\mu\in BC^\infty(\Del) \cap L^2(\Del,\rho(z)d^2z)$ with
  compact support and $\nu\in BC^\infty(\Del)$ we set
\begin{displaymath}
\mathcal{Q}_\nu(\mu) = 2\frac{\pa}{\pa
\z}\rho^{-1}\frac{\pa}{\pa \z} G(\mu\bar{\nu}).
\end{displaymath}
We will prove that
$\Vert\mathcal{Q}_\nu(\mu)\Vert^2_2=\Vert\mu\bar{\nu}\Vert^2_2 - (\mu\bar{\nu}, G(\mu\bar{\nu}))$,
 so that
$\mathcal{Q}_\nu$ extends to a bounded linear operator on
$L^2(\Del,\rho(z)d^2z)$. Since, according to Proposition
\ref{varQ1}, $\mathcal{Q}_\nu(\mu)=-\dot{Q}(\mu)[\nu]$ for $\mu\in
H^{-1,1}(\Del^*)$ and $\nu\in\Omega^{-1,1}(\Del^*)$, the assertion
follows from this fact.

From the explicit formula \eqref{resolvent-kernel} we get the following estimates
\begin{align} \label{estimate}
G(z,w) &= O((1-|z|^2)),\quad (\pa_z\rho^{-1}\pa_z)G(z,w)=O((1-|z|^2)) , \\
\pa_zG(z,w) &=O(1),\quad
\left(\pa_{\z}\left(\pa_z\rho^{-1}\pa_z\right)\right)G(z,w)=O(1)
\quad\text{as}\quad |z|\rightarrow 1, \nonumber
\end{align}
uniformly as $w$ varies on compact subsets of $\Del$. Using
Stokes' theorem and the identity
\begin{displaymath}
\rho^{-1}\pa_z\,\rho^{-1}\pa_z\,\rho\,\pa_{\z}\,\rho^{-1}\pa_{\z} = \Delta_0(\Delta_0 + \tfrac{1}{2})
\end{displaymath}
we get
\begin{align*}
\iint\limits_\Del \left|\mathcal{Q}_\nu(\mu)\right|^2\rho(z)d^2z
&=4\iint\limits_\Del\pa_z\left(\rho^{-1}\pa_zG(\bar{\mu}\nu)\right)
\pa_{\z}\left(\rho^{-1}\pa_{\z}G(\mu\bar{\nu})\right) \rho(z)d^2z \\
&=4\iint\limits_\Del \Delta_0(\Delta_0 + \tfrac{1}{2})G(\mu\bar{\nu})
G(\bar{\mu}\nu)\rho(z)d^2z \\
&= 2\iint\limits_\Del \Delta_0 (\mu\bar{\nu})
G(\bar{\mu}\nu)\rho(z)d^2z,
\end{align*}
where in the last line we have used property \textbf{RK3} from
Section \ref{resolvent}. Due to the estimates \eqref{estimate} the
boundary terms arising in the Stokes' formula vanish. Using
Stokes' theorem once again we finally  get
\begin{align*}
\iint\limits_\Del \left|\mathcal{Q}_\nu(\mu)\right|^2\rho(z)d^2z & =
 2\iint\limits_\Del \mu\bar{\nu}\Delta_0 G(\bar{\mu}\nu)\rho(z)d^2z \\
 &=\iint\limits_\Del |\mu\bar{\nu}|^2\rho(z)d^2z  - \iint\limits_\Del
  \mu\bar{\nu}G(\bar{\mu}\nu)\rho(z)d^2z \\
 & = \Vert\mu\bar{\nu}\Vert^2_2 - (\mu\bar{\nu}, G(\mu\bar{\nu})).
\end{align*}
The boundary terms again vanish due to \eqref{estimate} and Remark
3.2.
\end{proof}
\begin{corollary} \label{scalar-product}
For $\mu,\nu\in H^{-1,1}(\Del^*)$ and $\kappa\in\Omega^{-1,1}(\Del^*)$,
\begin{displaymath}
\left(\dot{Q}(\mu)[\kappa],\dot{Q}(\nu)[\kappa]\right) = (\mu\bar{\kappa},
\nu\bar{\kappa}) - \left(\mu\bar{\kappa}, G(\nu\bar{\kappa})\right).
\end{displaymath}
\end{corollary}
\begin{theorem} \label{metricfirst}
For $\mu, \nu \in H^{-1,1}(\Del^*)$ and $\kappa \in
\Omega^{-1,1}(\Del^*)$,
\begin{align*}
\left.\frac{\pa}{\pa \vep}g_{\mu \bar{\nu}}(\vep\kappa)\right|_ {\vep=0} =0.
\end{align*}
\end{theorem}
\begin{proof} Since $P^2 = P$, we get from  \eqref{WP-explicit},
\begin{equation} \label{WP-explicit-two}
g_{\mu\bar{\nu}}(\kappa)=\iint\limits_{\Del^*}Q(\mu,\kappa)
\ov{Q(\nu,\kappa)}(1-|\kappa|^2)w_{\kappa}^*(\rho)(z)d^2z =
\iint\limits_{\Del^*}\mu\ov{Q(\nu,\kappa)}
w_{\kappa}^*(\rho)(z)d^2z,
\end{equation}
so that
\begin{align*}
\left.\frac{\pa}{\pa \vep} g_{\mu \ov{\nu}}(\vep\kappa)\right|_{\vep=0}
= \left(\dot{Q}(\mu)[\kappa],\nu\right) + \left(\mu, \dot{Q}(\nu)[\kappa]\right)
=\left(\mu,\dot{Q}(\nu)[\kappa]\right).
\end{align*}
Differentiation under the integral sign is justified as in
\cite{Ahlfors3}. Thus for all $\mu, \nu \in H^{-1,1}(\Del^*)$ and
$\kappa \in \Omega^{-1,1}(\Del^*)$
\begin{displaymath}
\left(\dot{Q}(\mu)[\kappa],\nu\right) = 0,
\end{displaymath}
and the theorem follows.
\end{proof}
Let $\{\mu_n\}_{n=2}^{\infty}$ be an orthonormal basis for the
Hilbert space $H^{-1,1}(\Del^*)$,
\begin{equation*}
\mu_n(z) = -\sqrt{\frac{n^3-n}{8\pi}}(1-|z|^2)^2 \z^{-n-2},\;\;n=2,3,\dots,
\end{equation*}
and let $\{\vep_n\}_{n=2}^\infty$ be the corresponding Bers
coordinates on the chart $V_0$. Since $\Vert\mu\Vert_2 =2\Vert
D_0\beta(\mu)\Vert_2$, it follows from Section
 \ref{hilbert-manifold} that $\sum_{n=2}^\infty|\vep_n|^2<\tfrac{4\pi}{3}$.
 Denote by $\tfrac{\pa}{\pa\vep_n}$ the corresponding directional
derivatives --- the vector field $\tfrac{\pa}{\pa\vep_{\mu_n}}$ on
$V_0$, and set $g_{m\bar{n}} =g_{\mu_m\bar{\mu}_n}$. Since the
basis $\{\mu_n\}_{n=2}^{\infty}$ is orthonormal,
$g_{m\bar{n}}=\delta_{mn}$ at the origin of $T(1)$.
\begin{corollary} \label{kahler}
The Weil-Petersson metric  is  a K\"{a}hler metric on the Hilbert
manifold $T(1)$, and the Bers coordinates are geodesic coordinates
at the origin of $T(1)$.
\end{corollary}
\begin{proof}
It follows from Theorem \ref{metricfirst} that
\begin{displaymath}
\frac{\pa g_{m\bar{n}}}{\pa\vep_l}(0)=0.
\end{displaymath}
\end{proof}
\begin{remark}
Propositions \ref{varQ1} and \ref{var-in-L^2} and Theorem
\ref{metricfirst} generalize Wolpert's results for
finite-dimensional \Te spaces (see Lemma 2.7 and Theorem 2.9 in
\cite{Wol}) to the universal \Te space. In particular, our proof
of Theorem \ref{metricfirst} (after Proposition \ref{var-in-L^2}
has been established) is the same as in \cite{Wol}.
\end{remark}
\subsection{The second variation of the Weil-Petersson metric}

Due to Corollary \ref{kahler}, the Riemann tensor of the
Weil-Petersson metric at the origin of $T(1)$ is given by
\begin{equation*}
R_{k\bar{l} m\bar{n}} = - \frac{\pa^2
g_{k \bar{l} }}{\pa \vep_m \pa
\bar{\vep}_n}(0),
\end{equation*}
where we are using conventions of Yano and Bochner
\cite{Yano-Bochner} in Hermitian geometry.
\begin{theorem}\label{metricsecond}
For $\mu, \nu \in H^{-1,1}(\Del^*)$ and $\kappa \in
\Omega^{-1,1}(\Del^*)$,
\begin{align*}
\left.\frac{\pa^2}{\pa \vep\pa\bar{\vep}} g_{\mu \bar{\nu}}(\vep\kappa)
\right\vert_{\vep=0} =  \left(\mu\bar{\kappa}, G(\nu
\bar{\kappa})\right) +\left(\mu\bar{\nu}, G(|\kappa|^2)\right).
\end{align*}
\end{theorem}
\begin{proof}
Differentiating the representation \eqref{WP-explicit-two} for
$g_{\mu\bar{\nu}}(\vep\kappa)$ with respect to $\vep$ and
$\bar{\vep}$ we get
\begin{align*}
&\left.\frac{\pa^2}{\pa\vep\pa\bar{\vep}}\,
g_{\mu \bar{\nu}}(\vep\kappa)\right\vert_{\vep=0} \\ =&\left(
\left.\frac{\pa^2}{\pa\vep\pa\bar{\vep}}\,Q(\mu,\vep\kappa)\right\vert_{\vep=0},
\nu\right) +\left( \mu\bar{\nu}, \rho^{-1}
\left.\frac{\pa^2}{\pa\vep\pa\bar{\vep}}\,\rho^{\vep\kappa}\right\vert_{\vep=0}\right) \\
=&\left(
\left.\frac{\pa^2}{\pa\vep\pa\bar{\vep}}\,Q(\mu,\vep\kappa)\right\vert_{\vep=0},
\nu\right) +\left(\left.\frac{\pa}{\pa
\bar{\vep}}\right\vert_{\vep=0}Q(\mu,
\vep\kappa),\left.\frac{\pa}{\pa
\bar{\vep}}\right\vert_{\vep=0}Q(\nu, \vep\kappa)\right)\\ +
&\left(\mu,
\left.\frac{\pa^2}{\pa\vep\pa\bar{\vep}}\,Q(\nu,\vep\kappa)\right\vert_{\vep=0}
\right)+ \left( \mu\bar{\nu},
\rho^{-1}\left.\frac{\pa^2}{\pa\vep\pa\bar{\vep}}\,
\rho^{\vep\kappa}\right\vert_{\vep=0}\right)
-\left( \mu\bar{\nu},|\kappa|^2\right).
\end{align*}
The differentiation under the integral sign is justified as in
\cite{Ahlfors3}, provided that all integrals above are absolutely
convergent. This follows from Proposition \ref{secondvarmetric},
property \textbf{RK3} in Section \ref{resolvent}, Proposition
\ref{var-in-L^2} and the following
\begin{lemma} \label{second-in-L^2}
For $\mu\in H^{-1,1}(\Del^*)$ and $\nu\in\Omega^{-1,1}(\Del^*)$,
\begin{displaymath}
\left.\frac{\pa^2}{\pa\vep\pa\bar{\vep}}\,Q(\mu,\vep\nu)\right\vert_{\vep=0}
\in L^2(\Del^*,\rho(z)d^2z).
\end{displaymath}
\end{lemma}
We relegate the proof of this lemma to  Appendix B. Now comparing
the two expressions for the second variation of $g_{\mu\bar{\nu}}$
and using Corollary \ref{scalar-product} we get
\begin{align*}
\left(\mu,
\left.\frac{\pa^2}{\pa\vep\pa\bar{\vep}}Q(\nu,\vep\kappa)\right|_{\vep=0}
\right)=& -\left(\left.\frac{\pa}{\pa
\bar{\vep}}\right\vert_{\vep=0}Q(\mu,
\vep\kappa),\left.\frac{\pa}{\pa
\bar{\vep}}\right\vert_{\vep=0}Q(\nu,
\vep\kappa)\right)+\left( \mu\bar{\nu}, |\kappa|^2\right)\\ =&
\left(\mu\bar{\kappa}, G(\nu \bar{\kappa})\right) .
\end{align*}
Using Proposition \ref{secondvarmetric}, we finally obtain
\begin{align*}
\left.\frac{\pa^2}{\pa\vep\pa\bar{\vep}}\right\vert_{\vep=0}
g_{\mu \bar{\nu}}(\vep\kappa) =&\left( \mu\bar{\kappa}, G(\nu
\bar{\kappa})\right) +\left( \mu\bar{\nu}, G(|\kappa|^2)\right).
\end{align*}
\end{proof}

\begin{corollary} \label{Riemann-tensor-formula} At the origin of $T(1)$,
\begin{align*}
R_{k\bar{l} m\bar{n}} = -\left(\mu_k\bar{\mu}_l,
G(\bar{\mu}_m\mu_n)\right) - \left(\bar{\mu}_l\mu_m,
G(\bar{\mu}_k\mu_n)\right).
\end{align*}
\end{corollary}
\begin{proof} It follows from Theorem \ref{metricsecond} by polarization that
\begin{align*}
R_{\kappa\bar{\lambda}\mu\bar{\nu} }& = -\left.\frac{\pa^2}{\pa\vep_1\pa\bar{\vep}_2}
\right\vert_{\vep_1=
\vep_2=0}g_{\mu\bar{\nu}}(\vep_1\kappa+\vep_2\lambda)\\
&=-(\kappa\bar{\lambda},G(\bar{\mu}\nu)) - (\bar{\lambda}\mu,G(\bar{\kappa}\nu)).
\end{align*}
\end{proof}
\begin{remark}
For finite-dimensional \Te spaces this result was proved by
Wolpert \cite{Wol}. Except Lemma \ref{second-in-L^2}, our
derivation is the same as in \cite{Wol}.
\end{remark}
\subsection{Ricci and sectional curvatures}

The Ricci tensor at the origin of $T(1)$ for the orthonormal basis
$\{\mu_n\}_{n=2}^\infty$ of $H^{-1,1}(\Del^*)$ is defined by the
following series
\begin{align*}
\mathcal{R}_{k\bar{l}} = \sum_{n=2}^{\infty} R_{k\bar{n} n\bar{l}}.
\end{align*}

\begin{theorem} \label{Ricci}
The Ricci tensor at the origin of $T(1)$ is well-defined and is given by
\begin{align*}
\mathcal{R}_{k\bar{l}} = -\frac{13}{12\pi}\delta_{kl}.
\end{align*}
\end{theorem}
\begin{proof}
Set  $\mu=\mu_k, \nu=\mu_l$, and $\mathcal{R}_{\mu\bar{\nu}}=\mathcal{R}_{k\bar{l}}$. We have
\begin{gather*}
\mathcal{R}_{\mu\bar{\nu}} =- \frac{2}{\pi}\sum_{n=2}^{\infty}(n^3-n)
\left(\iint\limits_{\Del^*}\iint\limits_{\Del^*} G(z,w) \mu(z) z^{n-2}
\ov{\nu(w)} \bar{w}^{n-2} d^2w d^2z \right.\\
+\left. \iint\limits_{\Del^*}\iint\limits_{\Del^*}G(z,w)
 \bar{z}^{n-2}z^{n-2} \mu(w)\ov{\nu(w)}  \frac{(1-|z|^2)^2}{(1-|w|^2)^2} d^2wd^2z \right)\\ =-\frac{12}{\pi}
\iint\limits_{\Del^*}\iint\limits_{\Del^*}G(z,w)
\frac{\mu(z) \ov{\nu(w)}}{(1-z \bar{w})^4} d^2w d^2z \\
-\frac{3}{4\pi} \iint\limits_{\Del^*}\iint\limits_{\Del^*}G(z,w)
  \mu(w)\ov{\nu(w)}\rho(z)\rho(w)d^2w d^2z \\
  =I_1 + I_2.
\end{gather*}
For the second integral, we use property \textbf{RK4} in Section \ref{resolvent} and get
\begin{align*}
I_2 = -\frac{3}{4\pi} \iint\limits_{\Del} \mu(w)\ov{\nu(w)}\rho(w)d^2w
 =-\frac{3}{4\pi}g_{\mu\bar{\nu}}.
\end{align*}
For the first integral, we use projection formula
\eqref{projection} and get
\begin{gather*}
I_1=-\frac{36}{\pi^2} \iint\limits_{\Del^*}\iint\limits_{\Del^*}\iint\limits_{\Del^*}
G(z,w) \mu(z) \ov{\nu(v)} \frac{(1-|w|^2)^2}{(1-z\bar{w})^4
(1-w \bar{v})^4}  d^2vd^2w d^2z.
\end{gather*}
Let
\begin{align*}
B(z,v) =  \iint\limits_{\Del^*} G(z,w)
\frac{(1-|w|^2)^2}{(1-w\bar{v})^4(1-z\bar{w} )^4} d^2 w.
\end{align*}
The kernel $B(z,v)$ satisfies
\begin{equation} \label{B-formula}
B(z, v) = B\Bigl( \sigma z, \sigma v\Bigr) \sigma'(z)^2
\ov{\sigma'(v)}^2 \quad\text{for all}\quad \sigma \in \PSU(1,1),
\end{equation}
and
\begin{align*}
B(0,v) &=  \iint\limits_{\Del^*}\left( \frac{1}{2\pi}
\frac{1+|w|^2}{1-|w|^2} \log \frac{1}{|w|^2} - \frac{1}{\pi}
\right)\frac{(1-|w|^2)^2}{(1-\bar{v}w)^4}d^2w =\frac{1}{9}.
\end{align*}
Hence
\begin{align*}
B(z,v) = \frac{1}{9(1-z\bar{v})^4}
\end{align*}
and
\begin{align*}
I_1 &= -\frac{4}{\pi^2}\iint\limits_{\Del^*}\iint\limits_{\Del^*} \mu(z)
\ov{\nu(v)} \frac{1}{(1-z\bar{v})^4 }  d^2z d^2w \\
&=
-\frac{1}{3\pi} \iint\limits_{\Del^*} \mu(z) \ov{\nu(z)}
\rho(z)d^2z =-\frac{1}{3\pi}g_{\mu\bar{\nu}}.
\end{align*}
Therefore.
\begin{align*}
\mathcal{R}_{\mu\bar{\nu}} = -\left( \frac{3}{4\pi} +
\frac{1}{3\pi} \right)g_{\mu\bar{\nu}}= -\frac{13}{12\pi}g_{\mu\bar{\nu}}.
\end{align*}
\end{proof}

Since the Weil-Petersson metric on $T(1)$ is right-invariant, it
follows from Theorem \ref{Ricci} that the Ricci tensor is
well-defined everywhere on $T(1)$. Denote by $Ric_{WP}$
corresponding Ricci $(1,1)$-form on $T(1)$. In terms of Bers
coordinates $\{\vep_n\}_{n=2}^\infty$ on the coordinate chart
$V_\mu$ the Ricci form is given by
\begin{displaymath}
Ric_{WP} =\frac{i}{2}\sum_{k,l=2}^\infty \mathcal{R}_{k\bar{l}}\,d\vep_k\wedge d\bar{\vep}_l.
\end{displaymath}
\begin{corollary} The universal \Te space $T(1)$ is a K\"ahler-Einstein
 manifold with negative constant Ricci curvature,
\begin{displaymath}
Ric_{WP} = -\frac{13}{12\pi}\,\omega_{WP}.
\end{displaymath}
\end{corollary}
\begin{proof}
Since $(1,1)$-- forms $\omega_{WP}$ and $Ric$ are right-invariant,
the result immediately follows from Theorem \ref{Ricci}
\end{proof}
\begin{remark} For the dense submanifold $\Mob(S^1)\bk\Diff_+(S^1)$ of $T_0(1)$
the statement of Theorem \ref{Ricci} was established by different
methods in \cite{KY} and \cite{BR2,BR1}. The ``magic ratio''
$\tfrac{13}{12\pi}$ is omnipresent in mathematics related to
string theory.
\end{remark}

Let $\frac{\pa}{\pa t_{\mu}}, \frac{\pa}{\pa t_{\nu}}\in
T_0^{\RR}T(1)$ be real tangent vectors. According to
\cite{Yano-Bochner}, the sectional curvature of the section
spanned by these vectors is given by $R/g$, where
\begin{align}\label{curve1}
R&=R_{\mu\bar{\nu}\nu\bar{\mu}} +
R_{\nu\bar{\mu}\mu\bar{\nu}} -R_{\mu\bar{\nu}\mu\bar{\nu}} -
R_{\nu\bar{\mu}\nu\bar{\mu}},\\
g&=4g_{\mu\bar{\mu}}g_{\nu\bar{\nu}}-2|g_{\mu\bar{\nu}}|^2-2\re(g_{\mu\bar{\nu}})^2.\nonumber
\end{align}
Similarly, the holomorphic sectional curvature of the section spanned by the
holomorphic tangent vector $\frac{\pa}{\pa \vep_{\mu}}$, where
$g_{\mu\bar{\mu}}=1$, is given by $R_{\mu\bar{\mu}\mu\bar{\mu}}$.

As in the finite-dimensional case \cite{Wol}, we have
\begin{theorem}
The sectional and holomorphic sectional curvatures of $T(1)$ are
negative.
\end{theorem}
\begin{proof}
For  a section spanned by $\frac{\pa}{\pa\vep_\mu}$, the
holomorphic sectional curvature is obviously
negative: $\mu\neq 0$ so that $G(|\mu|^2)> 0$, and $(|\mu|^2,G(|\mu|^2))> 0$.

For a section spanned by the real tangent vectors $\frac{\pa }{\pa
t_{\mu}}$ and $\frac{\pa}{\pa t_{\nu}}$, using Cauchy-Schwarz
inequality, it is easy to see that $g$ is positive. Using
Corollary \ref{Riemann-tensor-formula}, the properties
\textbf{RK1} and \textbf{RK2}, we get
\begin{align*}
R= 4\text{Re}\; (\mu\bar{\nu}, G(\bar{\mu}\nu) ) -2(\mu\bar{\nu},
G(\mu\bar{\nu})) -2(|\mu|^2, G(|\nu|^2)).
\end{align*}
From the property \textbf{RK2} and Cauchy-Schwarz inequality we
have
\begin{gather*}
|G(\mu\bar{\nu})(z)|\leq \iint\limits_{\Del^*}G(z,w)^{1/2}|\mu(w)|G(z,w)^{1/2}|\nu(w)|\rho(w)d^2w \\
\leq \left(\iint\limits_{\Del^*}G(z,w)|\mu(w)|^2\rho(w)d^2w\right)^{1/2}
\left(\iint\limits_{\Del^*}G(z,w)|\nu(w)|^2\rho(w)d^2w\right)^{1/2}
\end{gather*}
so that
\begin{align*}
|(\bar{\mu}\nu, G(\mu\bar{\nu}))|& \leq
\iint\limits_{\Del^*}|\mu\nu| G(|\mu|^2)^{1/2}
G(|\nu|^2)^{1/2} \rho(z)d^2z\\
&\leq
\left(\iint\limits_{\Del^*}|\mu|^2G(|\nu|^2)\rho(z)d^2z\right)^{1/2}
\left(\iint\limits_{\Del^*}|\nu|^2G(|\mu|^2)\rho(z)d^2
z\right)^{1/2}\\
&=\bigl(|\mu|^2, G(|\nu|^2)\bigr).
\end{align*}
On the other hand, if we let $\mu\bar{\nu}=\alpha+i\beta$, where
$\alpha$ and $\beta$ are real--valued functions, then
$\text{Re}\;(\mu\bar{\nu}, G(\bar{\mu}\nu) )=(\alpha,
G(\alpha))-(\beta, G(\beta))$ and $(\mu\bar{\nu},
G(\mu\bar{\nu}))=(\alpha, G(\alpha))+(\beta, G(\beta))$. Since for
a bounded real-valued smooth function $h$, we can deduce from the
property \textbf{RK3} and the positivity of the operator $\Delta_{0}+\tfrac{1}{2}$ that $(h, G(h))\geq 0$,
it follows that $\text{Re}\,(\mu\bar{\nu}, G(\bar{\mu}\nu) )\leq
(\mu\bar{\nu}, G(\mu\bar{\nu}))$. Hence $R$ is negative by
Cauchy-Schwarz inequality.
\end{proof}

\section{Finite-dimensional Teichm\"uller spaces} \label{finite-dim}
Curvature properties of finite-dimensional \Te spaces were
extensively studied by Ahlfors \cite{Ahlfors3}, Royden
\cite{Royden}, and especially by Wolpert \cite{Wol}. Here, for the
\Te space $T(\Gamma)$ of a cofinite Fuchsian group $\Gamma$  we
show how to get Wolpert's explicit formulas from the curvature
formulas of the Hilbert manifold $T(1)$, derived in Section
\ref{Riemann-tensor}.

First note that the canonical embedding  of the finite-dimensional
complex manifold $T(\Gamma)$ into $T(1)$ is holomorphic with
respect to the Banach manifold structure on $T(1)$ but not with
respect to the Hilbert manifold structure on $T(1)$. Indeed, for a
cofinite Fuchsian group $\Gamma$ the finite-dimensional vector
space $\Omega^{-1,1}(\Del^*,\Gamma)$ is not a subspace of the
Hilbert space $H^{-1,1}(\Del^*)$ , but rather
\begin{displaymath}
\Omega^{-1,1}(\Del^*,\Gamma)\cap H^{-1,1}(\Del^*) =\{0\}.
\end{displaymath}
Thus the Weil-Petersson metric on $T(\Gamma)$, defined in Section
\ref{fuchsian-groups}, is not a pull-back of the Weil-Petersson
metric on $T(1)$. However, due to Lemma \ref{regularization} we
can represent  the Petersson inner product on the tangent space at
the origin of $T(\Gamma)$ as an  ``average'' of the inner products
in $T(1)$. Namely, using the canonical complex anti-linear
isomorphisms $\Omega^{-1,1}(\Del^*,\Gamma)\simeq
\Omega^{-1,1}(\Del,\Gamma)$ and $H^{-1,1}(\Del^*)\simeq
H^{-1,1}(\Del)$, we have
\begin{align*}
\la\mu, \nu\ra_{WP} =\iint\limits_{\Gamma\bk\Del} \mu\bar{\nu}\rho(z)d^2z
 &=\lim_{r \rightarrow 1^-}
\frac{A(\Gamma\bk \Del)}{A(\Del_r)} \iint\limits_{\Del_r}  \mu\bar{\nu}\rho(z)d^2z\\
&=\lim_{r \rightarrow 1^-}
\frac{A(\Gamma\bk \Del)}{A(\Del_r)} \iint\limits_{\Del}  \mu_r\bar{\nu}_r\rho(z)d^2z.
\end{align*}
Here $\mu,\nu\in\Omega^{-1,1}(\Del,\Gamma)$ and $\mu_r =\chi_r\mu,
\nu_r=\chi_r\nu$, where $\chi_r$ is the characteristic function of
$\Del_r=\{z\in\Del : |z|\leq r\}$. In what follows we will denote
by $(~,~)_\Gamma$ the Petersson inner  product $\la~,~\ra_{WP}$ in
$\Omega^{-1,1}(\Del, \Gamma)$, as well as the inner product for
the Hilbert space $L^2(\Gamma\bk\Del,\rho(z)d^2z)$. Since they are
given by the same formula, there would be no confusion. Moreover,
for $\mu\in\Omega^{-1,1}(\Del,\Gamma)$, $|\mu|\in
L^2(\Gamma\bk\Del,\rho(z)d^2z)$.

\begin{lemma} \label{projappr} Let $\mu\in\Omega^{-1,1}(\Del)$ and
$\nu\in L^\infty(\Del)\cap L^1(\Del,\rho(z)d^2z)$. Then
\begin{itemize}
\item[(i)] For $0<r<1$,
\begin{displaymath}
P(\mu_r)\in H^{-1,1}(\Del)\;\;\text{and}\;\;
P(\mu_r)(z) = O\left((1-|z|^2)^2\right)\;\;\text{as}\;\; |z|\rightarrow 1.
\end{displaymath}
\item[(ii)]
\begin{displaymath}
\lim_{r\rightarrow 1^-}\iint\limits_{\Del} P(\mu_r)\bar{\nu}\,\rho(z)d^2z =
\iint\limits_{\Del} \mu\bar{\nu}\,\rho(z)d^2z.
\end{displaymath}
\end{itemize}
\end{lemma}
\begin{proof}
Since
\begin{align*}
\iint\limits_{\Del} \left|P(\mu_r)\right|^2 \rho(z)d^2z
=\iint\limits_{\Del}
P(\mu_r)\bar{\mu}_r\,\rho(z)d^2z=\iint\limits_{\Del_r}
P(\mu_r)\bar{\mu}_r\,\rho(z)d^2z<\infty,
\end{align*}
$P(\mu_r)\in H^{-1,1}(\Del)$. Using
$\mu(z)=-\tfrac{(1-|z|^2)^2}{2}\sum_{n=2}^{\infty} (n^3-n) a_n
\z^{n-2}$ and \eqref{projection}, we get
\begin{align*}
P(\mu_r)(z) = -\frac{(1-|z|^2)^2}{4}\sum_{n=2}^{\infty} (n^3-n) a_n
\left(\frac{r^{2n+2}}{n+1}-\frac{2r^{2n}}{n}+\frac{r^{2n-2}}{n-1}\right)\z^{n-2},
\end{align*}
so that $(1-|z|^2)^{-2}P(\mu_r)(z)$ is continuous on $|z|=1$.

To prove part (ii), consider the estimate
\begin{align*}
\left|\left(\mu-P(\mu_r)\right)(z)\right| \leq &\frac{3(1-|z|^2)^2}{\pi}\Vert\mu\Vert_{\infty}\;
\iint\limits_{\Del\setminus \Del_r}
\frac{d^2u}{|1-u\z|^4}\\
=&3\Vert\mu\Vert_{\infty}(1-|z|^2)^2\sum_{n=1}^{\infty} n
|z|^{2n-2}(1-r^{2n})\\
=&3\Vert\mu\Vert_{\infty}\left(1-\frac{r^2(1-|z|^2)^2}{(1-r^2|z|^2)^2}\right).
\end{align*}
For fixed $r$ the right hand side of this estimate is an increasing
function of $|z|$, so that
\begin{align*}
\sup_{|z|\leq s} \left|\left(\mu-P(\mu_r)\right)(z)\right| \leq
3\Vert\mu\Vert_{\infty}\left(1-\frac{r^2(1-s^2)^2}{(1-r^2s^2)^2}\right),
\end{align*}
and for fixed $s$,
\begin{align*}
\lim_{r\rightarrow 1^-} \sup_{|z|\leq s}
\left|\left(\mu-P(\mu_r)\right)(z)\right|=0.
\end{align*}
Also for fixed $r$ we have the estimate
\begin{displaymath}
\Vert\mu - P(\mu_r)\Vert_{\infty}\leq 3\Vert\mu\Vert_{\infty}.
\end{displaymath}
Now since $\nu\in L^1(\Del,\rho(z)d^2z)$, for every $\vep>0$ there
exists $0<s<1$ such that
\begin{align*}
\iint\limits_{\Del\setminus \Del_s} \left| \nu \right|\rho(z)d^2z \leq
\vep,
\end{align*}
and we obtain,
\begin{gather*}
\left| \iint\limits_{\Del} \left(\mu -P(\mu_r)\right) \bar{\nu}
\rho(z)d^2z\right| \\
\leq \iint\limits_{\Del_s}
\left|\left(\mu-P(\mu_r)\right)\right||\nu|\,\rho(z)d^2z
+\iint\limits_{\Del\setminus \Del_s}
\left|\left(\mu-P(\mu_r)\right)\right||\nu|\,\rho(z)d^2z\\
\leq \sup_{|z|\leq s}
\left|\left(\mu-P(\mu_r)\right)(z)\right|\iint\limits_{\Del}|\nu|\,\rho(z)d^2z
+ 3\vep\Vert\mu\Vert_{\infty}.
\end{gather*}
Passing to the limit $r\rightarrow 1^-$, we get
\begin{align*}
\lim_{r\rightarrow 1^-}\left| \iint\limits_{\Del} \left(\mu
-P(\mu_r)\right)\bar{\nu}\,\rho(z)d^2z \right|\leq 3\vep\Vert\mu
\Vert_{\infty}.
\end{align*}
Since $\vep$ is arbitrary, the result follows.
\end{proof}

\begin{lemma}\label{regularmetric}
For $\mu ,\nu\in\Omega^{-1,1}(\Del,\Gamma)$,
\begin{align*}
(\mu, \nu)_{\Gamma} &=\lim_{r \rightarrow
1^-}\lim_{s\rightarrow 1^-} \frac{ A(\Gamma\bk \Del)}
{A(\Del_r)}\iint\limits_{\Del} P(\mu_s)\ov{P(\nu_r)}\rho(z)d^2z.
\end{align*}
\end{lemma}
\begin{proof}
Since $\mu\in\Omega^{-1,1}(\Del,\Gamma)\subset\Omega^{-1,1}(\Del)$,
\begin{align*}
\iint\limits_{\Del} \mu_r\bar{\nu}_r\rho(z)d^2z = \iint\limits_{\Del}\mu
\bar{\nu}_r\rho(z)d^2z=\iint\limits_{\Del} \mu \ov{P(\nu_r)}\rho(z)d^2z.
\end{align*}
According to part (i) of Lemma \ref{projappr}, $P(\nu_r)\in
L^1(\Del,\rho(z)d^2z)$ for $0<r<1$, so that the result follows
from part (ii) of Lemma \ref{projappr}.
\end{proof}

\begin{remark}
The limits in Lemma \ref{regularmetric} can not be interchanged.
Indeed, it follows from part (ii) of Lemma \ref{projappr} that for
fixed $s<1$  the limit $r\rightarrow 1$ is always zero.
\end{remark}
In a neighborhood of the origin in $T(\Gamma)$ the Weil-Petersson
metric is given by
\begin{align*}
g_{\mu\bar{\nu}} (\kappa) =\iint\limits_{\Gamma_{\kappa}\bk \Del}
P(R(\mu, \kappa))\ov{P(R(\nu, \kappa))}\rho(z)d^2z,
\end{align*}
where
$\kappa\in\Omega^{-1,1}(\Del,\Gamma),\,\Vert\kappa\Vert_\infty$ is
sufficiently small, and $\Gamma_{\kappa}
=w_{\kappa}\circ\Gamma\circ w_{\kappa}^{-1}$.
\begin{lemma}\label{reg2}
Let $\mu,\nu\in\Omega^{-1,1}(\Del,\Gamma)$. For $\kappa\in\Omega^{-1,1}
(\Del,\Gamma),\,\Vert\kappa\Vert_{\infty}$ sufficiently small,
\begin{align*}
g_{\mu\bar{\nu}} (\kappa) =&\lim_{r \rightarrow
1^-}\lim_{s\rightarrow 1^-} \frac{A(\Gamma\bk \Del)}{A(\Del_r)}
\iint\limits_{\Del} P(R(P(\mu_s), \kappa))\ov{P(R(P(\nu_r), \kappa))}
\rho(z)d^2z.
\end{align*}
\end{lemma}
\begin{proof}
First, we have
\begin{gather*}
\iint\limits_{\Del} P(R(P(\mu_s), \kappa))\ov{P(R(P(\nu_r), \kappa))}
\rho(z)d^2z \\
 =\frac{12}{\pi} \iint\limits_{\Del}\iint\limits_{\Del}\frac{P(\mu_s)(u)(w_{\kappa})_u(u)^2
\ov{P(\nu_r)(z) (w_{\kappa})_z(z)^2}
}{(1-w_{\kappa}(u)\ov{w_{\kappa}(z)})^4}d^2zd^2u.
\end{gather*}
Since $\rho P(\nu_r)$ is bounded on $\Del $, and for $\Vert \kappa
\Vert_{\infty}$ sufficiently small $(1/2)\rho \leq
w_{\kappa}^*\rho \leq (3/2) \rho$, we conclude that $\rho
R(P(\nu_r), \kappa)$ is also bounded on $\Del $. As a result,
\begin{gather*}
\iint\limits_{\Del}\left|\iint\limits_{\Del}\frac{(w_{\kappa})_u(u)^2
\ov{P(\nu_r)(z) (w_{\kappa})_z(z)^2}
}{(1-w_{\kappa}(u)\ov{w_{\kappa}(z)})^4}d^2z\right|d^2u \\
=\iint\limits_{\Del} \left|\iint\limits_{\Del}\frac{
\ov{R(P(\nu_r),\kappa)(z)} }{(1-u\z)^4}d^2z\right| \frac{d^2u}{1-
|\kappa(u)|^2} \\
\leq  C\iint\limits_{\Del}\iint\limits_{\Del}\frac{ (1-|z|^2)^2
}{|1-u\z|^4}d^2zd^2u =\pi^2 C<\infty.
\end{gather*}
It follows from part (ii) of Lemma \ref{projappr} that
\begin{gather*}
\lim_{s\rightarrow 1^-}
\iint\limits_{\Del} P(R(P(\mu_s), \kappa))\ov{P(R(P(\nu_r), \kappa))}
\rho(z)d^2z\\
=\frac{12}{\pi}
\iint\limits_{\Del}\iint\limits_{\Del}\frac{\mu(u)(w_{\kappa})_u(u)^2
\ov{P(\nu_r)(z)(w_{\kappa})_z(z)^2}
}{(1-w_{\kappa}(u)\ov{w_{\kappa}(z)})^4}d^2zd^2u \\
=\iint\limits_{\Del} P(R(\mu, \kappa))\ov{P(R(P(\nu_r), \kappa))}\rho(z)d^2z.
\end{gather*}
Now
\begin{gather*}
\iint\limits_{\Del} P(R(\mu, \kappa))\ov{P(R(P(\nu_r), \kappa))}
\rho(z)d^2z\\
=\frac{144}{\pi^2}\iint\limits_{\Del}\iint\limits_{\Del_r}\iint\limits_{\Del}
 \frac{\mu(u)(w_{\kappa})_u(u)^2\ov{\nu(v)(w_{\kappa})_z(z)^2}}
 {(1-w_{\kappa}(u)\ov{w_{\kappa}(z)})^4(1-z\bar{v})^4}\rho(z)^{-1}d^2ud^2vd^2z\\
=\frac{144}{\pi^2}\iint\limits_{\Del_r}\lambda(v)\ov{\nu(v)}\rho(v)d^2v,
\end{gather*}
where
\begin{align*}
\lambda(v)=\rho(v)^{-1}\iint\limits_{\Del}\iint\limits_{\Del}\frac{\mu(u)(w_{\kappa})_u(u)^2
\ov{(w_{\kappa})_z(z)^2}\rho(z)^{-1}}{(1-w_{\kappa}(u)
\ov{w_{\kappa}(z)})^4(1-z\bar{v})^4}d^2ud^2z.
\end{align*}
Since $\mu\in\Omega^{-1,1}(\Del,\Gamma)$ and
$w_{\kappa}\circ\Gamma\circ w_{\kappa}^{-1}
=\Gamma_{\kappa}\subseteq \PSU(1,1)$, it is easy to see that
$\lambda\in\Omega^{-1,1}(\Del,\Gamma)$. Using Lemma
\ref{regularization}, we get
\begin{gather*}
\lim_{r\rightarrow 1^{-}}\frac{A(\Gamma\bk\Del)}{A(\Del_r)}
\iint\limits_{\Del} P(R(\mu, \kappa))\ov{P(R(P(\nu_r), \kappa))}
\rho(z)d^2z \\
=\frac{144}{\pi^2}\iint\limits_{\Gamma\bk
\Del}\lambda(v)\ov{\nu(v)}\rho(v)d^2v\\
=\frac{144}{\pi^2}\iint\limits_{\Del}\iint\limits_{\Gamma\bk\Del}
\iint\limits_{\Del}\frac{\mu(u)(w_{\kappa})_u(u)^2\ov{\nu(v)
(w_{\kappa})_z(z)^2}}{(1-w_{\kappa}(u)
\ov{w_{\kappa}(z)})^4(1-z\bar{v})^4}
\rho(z)^{-1}d^2ud^2vd^2z\\
=\frac{144}{\pi^2}\iint\limits_{\Gamma\bk\Del}\iint\limits_{\Del}
\iint\limits_{\Del}\frac{\mu(u)(w_{\kappa})_u(u)^2\ov{\nu(v)
(w_{\kappa})_z(z)^2}}{(1-w_{\kappa}(u)
\ov{w_{\kappa}(z)})^4(1-z\bar{v})^4}
\rho(z)^{-1}d^2ud^2vd^2z\\
=\iint\limits_{\Gamma\bk\Del} P(R(\mu, \kappa))\ov{P(R(\nu, \kappa))}
\rho(z)d^2z,
\end{gather*}
where we have used the fact that the integrals above do not change
if we let any one of the integration variables to range over
$\Gamma\bk\Del$ while others range over $\Del$ (cf.
\cite{Ahlfors3}). The latter property follows from the fact that
$\mu,\nu$ and $\kappa$ are $(-1,1)$ tensors for $\Gamma$, and the
representation
$\Del=\bigcup_{\gamma\in\Gamma}\,\gamma(\Gamma\bk\Del)$.
\end{proof}
\begin{theorem} \label{WPvariation-Gamma}
For $\mu,\nu,\kappa\in\Omega^{-1,1}(\Del,\Gamma)$,
\begin{align*}
\left.\frac{\pa}{\pa\vep}\right|_{\vep=0}g_{\mu\bar{\nu}}(\vep\kappa)=& 0,\\
\left.\frac{\pa^2}{\pa
\vep\pa\bar{\vep}}\right|_{\vep=0}g_{\mu\bar{\nu}}(\vep\kappa)
=&\left(\mu\bar{\kappa}, G_{\Gamma}(\nu\bar{\kappa})\right)_{\Gamma} +\left(\mu\bar{\nu},
G_{\Gamma}(|\kappa|^2)\right)_{\Gamma}.
\end{align*}
\end{theorem}
\begin{proof}
We will use Lemma \ref{reg2} and Theorems \ref{metricfirst} and
\ref{metricsecond}, provided one can interchange $\tfrac{\pa}{\pa
\vep},\,\tfrac{\pa^2}{\pa\vep\pa\bar{\vep}}$ with the limits. This
can be done as in \cite{Ahlfors3} by showing that limits of
corresponding derivatives converge uniformly on $\vep$ in  a
neighborhood of $0$. We omit these standard arguments and
concentrate on actual computations.

For the first variation of the Weil-Petersson metric we get
\begin{align*}
\left.\frac{\pa}{\pa\vep}\right|_{\vep=0}g_{\mu\bar{\nu}}(\vep\kappa)=
\lim_{r\rightarrow 1^-}\lim_{s\rightarrow 1^-} \frac{ A(\Gamma\bk
\Del)}{A(\Del_r)} \left.\frac{\pa}{\pa\vep}\right|_{\vep=0}g_{P(\mu_s)\ov{P(\nu_r)}}(\vep
\kappa).
\end{align*}
Since $P(\mu_s), P(\nu_r)\in H^{-1,1}(\Del)$, we conclude from
Theorem \ref{metricfirst} that this is identically zero.

Similarly, for the second variation we have
\begin{align*}
\left.\frac{\pa^2}{\pa
\vep\pa\bar{\vep}}\right|_{\vep=0}g_{\mu\bar{\nu}}(\vep\kappa)
=&\lim_{r \rightarrow 1^-}\lim_{s\rightarrow 1^-} \frac{A(\Gamma\bk \Del)}
{A(\Del_r)} \left.\frac{\pa^2}{\pa
\vep\pa\bar{\vep}}\right|_{\vep=0}g_{P(\mu_s)\ov{P(\nu_r)}}
(\vep\kappa).
\end{align*}
Since $P(\nu_r), P(\mu_s)\in H^{-1,1}(\Del)$ and
$\kappa\in\Omega^{-1,1}(\Del,\Gamma)\subset\Omega^{-1,1}(\Del)$,
 we get from Theorem \ref{metricsecond}\footnote{It is for this case that we need the condition
$\kappa\in\Omega^{-1,1}(\Del)$ in Theorem \ref{metricsecond}.},
\begin{gather*}
\left.\frac{\pa^2}{\pa
\vep\pa\bar{\vep}}\right|_{\vep=0}g_{\mu\bar{\nu}}(\vep\kappa)  \\
= \lim_{r \rightarrow 1^-}\lim_{s\rightarrow 1^-} \frac{
A(\Gamma\bk \Del)}{A(\Del_r)}\left(\Bigl(P(\mu_s)\bar{\kappa},
G(P(\nu_r)\bar{\kappa})\Bigr)
 +\left( P(\mu_s)\ov{P(\nu_r)},
G(|\kappa|^2)\right)\right).
\end{gather*}
By properties \textbf{RK2} and \textbf{RK4} in Section \ref{resolvent},
\begin{align*}
\iint\limits_{\Del} \left| G(P(\nu_r)\kappa)\right|
|\kappa|\rho(z)d^2z\leq &\Vert\kappa\Vert_{\infty}^2
\iint\limits_{\Del}\iint\limits_{\Del}
G(z,w)\left| P(\nu_r)(w)\right|\rho(z)\rho(w)d^2wd^2z\\
=&\Vert\kappa\Vert_{\infty}^2\iint\limits_{\Del}
\left| P(\nu_r)(w)\right|\rho(w)d^2w  <\infty,
\end{align*}
and by property \textbf{RK3} $G(|\kappa|^2)$ is bounded on $\Del$,
so that it follows from Lemma \ref{projappr} that
\begin{align*}
\left.\frac{\pa^2}{\pa
\vep\pa\bar{\vep}}\right|_{\vep=0}g_{\mu\bar{\nu}}(\vep\kappa) =
\lim_{r \rightarrow 1^-} \frac{ A(\Gamma\bk
\Del)}{A(\Del_r)}\left(\Bigl(\mu\bar{\kappa},
G(P(\nu_r)\bar{\kappa})\Bigr)
 +\left(\mu\ov{P(\nu_r)},
G(|\kappa|^2)\right)\right).
 \end{align*}
We have
\begin{align*}
\Bigl(\mu\bar{\kappa}, G(P(\nu_r) \bar{\kappa})\Bigr)
=&\iint\limits_{\Del_r}\lambda_1(v)\ov{\nu(v)}\rho(v)d^2v, \\
\left(\mu\ov{P(\nu_r)}, G(|\kappa|^2)\right)=&\iint\limits_{\Del_r}\lambda_2(v)\ov{\nu(v)}\rho(v)d^2v,
\end{align*}
where
\begin{align*}
\lambda_1(v)=&\frac{12}{\pi}\rho(v)^{-1}\iint\limits_{\Del}\iint\limits_{\Del} \frac{\mu(z)
\ov{\kappa(z)}G(z,u) \kappa(u)}{(1-u\bar{v})^4}\rho(z)d^2u d^2z,\\
\lambda_2(v)=&\frac{12}{\pi}\rho(v)^{-1}\iint\limits_{\Del}\iint\limits_{\Del}
\frac{\mu(z)G(z,u)|\kappa(u)|^2}{(1-z\bar{v})^4}\rho(u)d^2u d^2z,
\end{align*}
and $\lambda_1, \lambda_2\in\Omega^{-1,1}(\Del,\Gamma)$. It follows from Lemma \ref{regularization} that
\begin{gather*}
\left.\frac{\pa^2}{\pa
\vep\pa\bar{\vep}}\right|_{\vep=0}g_{\mu\bar{\nu}}(\vep\kappa)
= \iint\limits_{\Gamma\bk\Del}\lambda_1(v)\ov{\nu(v)}\rho(v)d^2v + \iint\limits_{\Gamma\bk\Del}
\lambda_2(v)\ov{\nu(v)}\rho(v)d^2v \\
=\frac{12}{\pi}\iint\limits_{\Del}\iint\limits_{\Gamma\bk\Del}
\iint\limits_{\Del} \mu(z) \ov{\kappa(z)}G(z,u)
\kappa(u)\frac{\ov{\nu(v)}}{(1-u\bar{v})^4}\rho(z)d^2u d^2vd^2z\\
+\frac{12}{\pi}
\iint\limits_{\Del}\iint\limits_{\Gamma\bk\Del}\iint\limits_{\Del}
\mu(z)G(z,u)|\kappa(u)|^2\frac{\ov{\nu(v)}}{(1-z\bar{v})^4}\rho(u)d^2u
d^2vd^2z\\
=I_1 + I_2.
\end{gather*}
As before, the integrals above do not change if we let any one of
the integration variables to range over $\Gamma\bk\Del$ while
others range over $\Del$. We have, using
property \textbf{RK1} , \eqref{projection} and \eqref{G-Gamma}, 
\begin{gather*}
I_1  =\frac{12}{\pi}\iint\limits_{\Del}\iint\limits_{\Del}
\iint\limits_{\Gamma\bk\Del} \mu(z) \ov{\kappa(z)}G(z,u)
\kappa(u)\frac{\ov{\nu(v)}}{(1-u\bar{v})^4}\rho(z)d^2u d^2vd^2z \\
=\iint\limits_{\Del}\iint\limits_{\Gamma\bk\Del} \mu(z)
\ov{\kappa(z)}G(z,u) \kappa(u)\ov{\nu(u)}\rho(u)\rho(z)d^2ud^2z \\
=\iint\limits_{\Gamma\bk\Del}\iint\limits_{\Gamma\bk\Del} \mu(z)
\ov{\kappa(z)}G_{\Gamma}(z,u) \kappa(u)\ov{\nu(u)}\rho(u)\rho(z)d^2ud^2z \\
=\left(\mu\bar{\kappa},G_{\Gamma}(\nu\bar{\kappa})\right)_{\Gamma}.
\end{gather*}
Similarly,
\begin{gather*}
I_2 = \iint\limits_{\Del}\iint\limits_{\Gamma\bk\Del}\mu(z)\ov{\nu(z)}G(z,u)
\left|\kappa(u)\right|^2\rho(u)\rho(z)d^2ud^2z \\
=\iint\limits_{\Gamma\bk\Del}\iint\limits_{\Gamma\bk\Del}\mu(z)\ov{\nu(z)}G_{\Gamma}(z,u)
\left|\kappa(u)\right|^2\rho(u)\rho(z)d^2ud^2z \\
=\left(\mu\bar{\nu},G_{\Gamma}(|\kappa|^2)\right)_{\Gamma},
\end{gather*}
and the assertion follows.
\end{proof}
\begin{remark} Theorem \ref{WPvariation-Gamma} was proved by Wolpert \cite{Wol},
 and all results on Ricci,  sectional, and scalar
curvatures for finite-dimensional \Te spaces follow from it.
\end{remark}

We conclude this section by deriving a formula for Ricci tensor different from \cite{Wol},
 and indicating its application.  Let $\mu_1, \ldots, \mu_d$ be an orthonormal basis of
$\Omega^{-1,1}(\Del, \Gamma)$, which is a subspace of the Hilbert
space $L^2(\Del,\Gamma)$ of Beltrami differentials $\mu$ for
$\Gamma$ such that $|\mu|\in L^2(\Gamma\bk\Del,\rho(z)d^2z)$. Let
$P:L^2(\Del,\Gamma)\rightarrow\Omega^{-1,1}(\Del,\Gamma)$ be the
orthogonal projector. It follows from the definition and
representation \eqref{projection} that $P$ is an integral operator
with kernel
\begin{align*}
P(z,w)=\sum_{n=1}^d
\mu_n(z)\ov{\mu_n(w)}=\frac{12}{\pi}\rho(z)^{-1}\rho(w)^{-1}\sum_{\gamma\in
\Gamma} \frac{\gamma'(w)^2}{(1-\z\gamma(w))^4}.
\end{align*}
The Ricci tensor at the origin of $T(\Gamma)$ is given by
\begin{gather} \label{ricci-gamma}
\mathcal{R}_{\mu \bar{\nu}} =\sum_{n=1}^{d} R_{\mu
\bar{\mu}_n\mu_n \bar{\nu}}=-\sum_{n=1}^{d} \left(\left(
\mu\bar{\nu},G_{\Gamma}(|\mu_n|^2)\right)_{\Gamma}+\left(\mu
\bar{\mu}_n,G_{\Gamma}(\nu\bar{\mu}_n)\right)_{\Gamma}\right)\\ \nonumber
=-\sum_{n=1}^d\left(\iint\limits_{\Gamma\bk\Del}\iint\limits_{\Gamma\bk\Del}
\mu(z)\ov{\nu(z)}G_{\Gamma}(z,w) |\mu_n(w)|^2 \rho(w)\rho(z)d^2wd^2z\right. \\ \nonumber
+\left.\iint\limits_{\Gamma\bk\Del}\iint\limits_{\Gamma\bk \Del}
\mu(z)\ov{\mu_n(z)}G_{\Gamma}(z,w)\mu_n(w)\ov{\nu(w)}\rho(w)
\rho(z)d^2wd^2z\right) \\ \nonumber
=-\frac{12}{\pi}\iint\limits_{\Del}\iint\limits_{\Gamma\bk\Del}
\mu(z)\ov{\nu(z)}G(z,w) \sum_{\gamma\in\Gamma}\frac{\gamma'(w)^2}
{(1-\bar{w}\gamma(w))^4} \rho(w)^{-1}\rho(z)d^2wd^2z\\ \nonumber
-\frac{12}{\pi}\iint\limits_{\Gamma\bk \Del}\iint\limits_{\Del }
\mu(z)\ov{\nu(w)}G(z,w)\sum_{\gamma\in
\Gamma}\frac{\gamma'(z)^2}{(1-\bar{w}\gamma(z))^4}d^2wd^2z,
\nonumber
\end{gather}
where nothing is changed if we let any of the integration
variables to range over $\Gamma\bk\Del$ while other range over
$\Del$.

It is instructive to compare the Ricci curvatures of the finite-dimensional \Te space $T(\Gamma)$ and that of
the universal \Te space $T(1)$.

First, $T(\Gamma)$ is no longer a \Ka-Einstein manifold. Second,
the sum over $\Gamma$ in \eqref{ricci-gamma} can be transformed
into a sum over the conjugacy classes  of $\Gamma$. As is in
\cite{TZ91}, using variational formulas for the Selberg
zeta-function, we find that the contribution of the hyperbolic
conjugacy  classes is the second variation of the Selberg
zeta-function at $s=2$. The contribution of parabolic conjugacy
classes (if they are present) yields a new \Ka metric on
$T(\Gamma)$, introduced in \cite{TZ91}. The contribution of the
identity element, as it follows from Theorem \ref{Ricci}, is
\begin{gather*}
-\frac{3}{4\pi}\iint\limits_{\Del}\iint\limits_{\Gamma\bk\Del}
\mu(z)\ov{\nu(z)}G(z,w)\rho(w)\rho(z)d^2wd^2z\\
- \frac{12}{\pi}\iint\limits_{\Gamma\bk \Del}\iint\limits_{\Del }
\mu(z)\ov{\nu(w)}G(z,w)\frac{1}{(1-z\bar{w})^4}d^2wd^2z= - \frac{13}{12\pi}(\mu,\nu)_{\Gamma}.
\end{gather*}
As a result, we obtain a local index theorem for families of
$\bar{\pa}$-operators acting on quadratic differentials on Riemann
surfaces, proved in \cite{TZ91}. The above arguments interpret it
as an ``averaged form" of Theorem \ref{Ricci}. Detailed derivation
of the local index theorem for families from \eqref{ricci-gamma}
will be presented elsewhere.

\appendix
\section{}
Here we prove that the connected component $T_{0}(1)$ of $0\in T(1)$ in the Hilbert manifold topology is the inverse image of $\beta(T(1))\cap A_2(\Del)$ under the Bers embedding, and that $T_0(1)$ and $\mathcal{T}_0(1)=\pi^{-1}(T_{0}(1))$ are topological groups.

It follows from Theorem \ref{bers-hilbert} that to prove $T_{0}(1)=\beta^{-1}(\beta(T(1))\cap A_2(\Del))$, it is sufficient to show that $\beta(T(1))\cap A_2(\Del)$ is connected.

\begin{theorem}
The submanifold $\beta(T(1))\cap A_{2}(\Del)$ of $\beta(T(1))$ is
connected in the Hilbert manifold topology of $\beta(T(1))$.
\end{theorem}
\begin{proof}
First, in the Hilbert manifold topology the set
$\beta(\Mob(S^1)\bk\Diff_+(S^1))$ is dense in $\beta(T(1))\cap A_2(\Del)$. Indeed, since
$\beta(T(1))\cap A_2(\Del)$ is open in $A_2(\Del)$, for every $\phi\in\beta(T(1))\cap A_2(\Del)$ there exists
$\delta$ such that $\tilde{\phi}\in \beta(T(1))\cap A_2(\Del)$ for all $\tilde{\phi}$ satisfying
$\Vert\tilde{\phi}-\phi\Vert_{2}<\delta$.
Since
$\phi(z)=\sum_{k=2}^{\infty}(k^3-k)a_k z^{k-2}\in
A_2(\Del)$, for every $n\in\mathbb{N}$ there exists $N\in \mathbb{N}$ such that
$\frac{\pi}{2}\sum_{k=N+1}^{\infty}(k^3-k)|a_k|^2<\tfrac{1}{n^{2}}$, so that for
$\phi_{n}(z)=\sum_{k=2}^{N}(k^3-k)a_kz^{k-2}$ we have $\Vert
\phi_{n}-\phi\Vert_2<\tfrac{1}{n}$. Thus for $\tfrac{1}{n}<\delta$ there exists
$\gamma_{n}=g^{-1}_{n}\circ f_{n}\in \Mob(S^1)\bk
\text{Homeo}_{qs}(S^1)$ such that
$\mathcal{S}(f_{n})=\phi_{n}$. Since $\phi_{n}$ is
analytic on an open domain containing $\Del\cup S^1$, so is
the function $f_{n}$ and, consequently, $f_{n}(S^1)$ is
an analytic curve. Hence $\left.\gamma_{n}\right\vert_{S^1}$ is
smooth and $\phi_{n}\in\beta(\Mob(S^1)\bk\Diff_+(S^1))$.

Secondly, the inclusion
$$\Mob(S^1)\bk\Diff_+(S^1)\hookrightarrow T(1)$$
is a continuous mapping between the Frechet and Hilbert manifolds. Since the Frechet manifold $\Mob(S^1)\bk\Diff_+(S^1)$ is connected, its image in the Hilbert manifold $T(1)$ is also
connected. Finally, since Bers embedding is a continuous mapping (actually biholomorphic) and the closure of a connected set is connected, we conclude that $\beta(T(1))\cap A_{2}(\Del)$ is connected in $\beta(T(1))$.
\end{proof}

\begin{corollary}\label{Hilbertcon}
The Hilbert submanifold $T_0(1)$ of $T(1)$ is characterized by
\begin{align*}
T_0(1)= \left\{ [\mu]\in T(1) : \beta([\mu])\in
A_2(\Del)\right\}.
\end{align*}
\end{corollary}

Next we prove that the Hilbert manifold $T_{0}(1)$ is
a topological group.
\begin{lemma}
Every $[\mu]\in T_{0}(1)$ has a representative $\mu\in L^2(\Del^*,
\rho(z)d^2z)\cap \mathcal{O}(\Del^*)_1$.
\end{lemma}
\begin{proof}
First observe that $L^2(\Del^*, \rho(z)d^2z)\cap
\mathcal{O}(\Del^*)_1$ is a subgroup of $L^{\infty}(\Del^*)_1$.
Indeed, let $\lambda,\nu\in L^2(\Del^*, \rho(z)d^2z)\cap
\mathcal{O}(\Del^*)_1$. Setting $\mu=\nu^{-1}$ and $\lambda_1
=\lambda, \lambda_2 =\nu$ in Lemma \ref{right1}, we get
\begin{equation*}
\Vert\lambda\ast\nu^{-1}\Vert_2 <C\Vert\lambda-\nu\Vert_2<\infty.
\end{equation*}

Denote by $B_{\alpha}$ the ball of radius $\alpha$ about the
origin in $A_2(\Del)$, where $\alpha < \sqrt{\tfrac{\pi}{3}}\, \delta$ and
$\delta$ is defined in Corollary \ref{righttran}, and let
$W_0=(\Phi\circ \Lambda) (\beta(V_0)\cap B_{\alpha})$. By
definition, $W_0\subset T_0(1)$ and every $[\mu] \in W_0$ has a
representative $ \mu\in L^2(\Del^*, \rho(z)d^2z)\cap
\mathcal{O}(\Del^*)_1$. For every $\nu \in L^{\infty}(\Del^*)_1$
set $W_{\nu} = R_{[\nu]}(W_0)$.

For every $[\mu]\in T_{0}(1)$ let $[0,1]\ni t\mapsto [\mu(t)]\in
T_0(1)$ be a path in $T_{0}(1)$ joining $0$ to $[\mu]$. By
definition of the Hilbert manifold structure on $T(1)$ (see
Section \ref{hilbert-manifold}), there exist $0=t_0< t_1<\ldots
<t_n=1$ and $W_i=W_{\mu(t_{i})}\subset T_0(1)$, which cover the
path and such that $W_{i-1}\cap W_{i}\neq \emptyset$,
$i=1,\dots,n$. We will prove by induction that each $[\mu(t_i)]$
has a representative $\mu_i\in L^2(\Del^*, \rho(z)d^2z)\cap
\mathcal{O}(\Del^*)_1$. For $1\leq i \leq n$, choose $[\nu_i] \in
W_{i-1}\cap W_{i}$. Since $[\nu_1]\in W_0$, it has a
representative $\nu_1\in L^2(\Del^*, \rho(z)d^2z)\cap
\mathcal{O}(\Del^*)_1$. On the other hand, since $[\nu_1]\in W_1$,
$[\nu_1]=[\lambda* \mu(t_1)]$ for some $\lambda\in L^2(\Del^*,
\rho(z)d^2z)\cap \mathcal{O}(\Del^*)_1$. Hence $[\mu(t_1)]$ has a
representative $\mu_1=\lambda^{-1}* \nu_1 \in L^2(\Del^*,
\rho(z)d^2z)\cap \mathcal{O}(\Del^*)_1$. Now suppose that
$[\mu(t_{i-1})]$ has a representative $\mu_{i-1}\in L^2(\Del^*,
\rho(z)d^2z)\cap \mathcal{O}(\Del^*)_1$. Since $[\nu_{i}]\in
W_{i-1}\cap W_i$, there exist $[\lambda_1], [\lambda_2] \in W_0$
such that $[\nu_i] = [\lambda_1* \mu(t_{i-1})]= [\lambda_2 *
\mu(t_i)]$. From the first equality it follows that $[\nu_i]$ has
a representative $\nu_i \in L^2(\Del^*, \rho(z)d^2z)\cap
\mathcal{O}(\Del^*)_1$, and the second equality implies that
$[\mu(t_{i})]$ has a representative $\mu_i \in L^2(\Del^*,
\rho(z)d^2z)\cap \mathcal{O}(\Del^*)_1$. The assertion follows.

\end{proof}

\begin{lemma}\label{group}
$T_{0}(1)$ is a subgroup of $T(1)$.
\end{lemma}
\begin{proof}
Every $[\lambda]\in T_{0}(1)$ has a representative $\lambda \in
L^2(\Del^*, \rho(z)d^2z)\cap \mathcal{O}(\Del^*)_1$. Using Lemma
\ref{rightinv} with $\mu=\lambda^{-1}$ we get by Corollary
\ref{Hilbertcon}
 that $[\lambda^{-1}]\in
T_{0}(1)$. Now let $[\nu]\in T_{0}(1)$ with a representative $\nu
\in L^2(\Del^*, \rho(z)d^2z)\cap \mathcal{O}(\Del^*)_1$. Using
again Lemma \ref{rightinv} with $\mu=\nu$ we get that
$[\lambda\ast\nu]\in T_{0}(1)$.

\end{proof}
The following lemmas and corollary are needed for the proof that
$T_0(1)$ is a topological group. For every $[\mu]\in T(1)$ let
\begin{displaymath}
L_{[\mu]}:T(1)\rightarrow T(1),\;\;\;[\nu]\mapsto[\mu\ast\nu],
\end{displaymath}
be the left translation by $[\mu]$.
\begin{lemma}\label{leftcont}
Left translations are continuous on $T_{0}(1)$.
\end{lemma}
\begin{proof}
Let $[\mu], [\kappa]\in T_{0}(1)$ with representatives $\mu,
\kappa\in L^2(\Del^*, \rho(z)d^2z)\cap \mathcal{O}(\Del^*)_1$, and
let $V_{\kappa}$ and $V_{\mu * \kappa}$ be coordinate charts
introduced in Section \ref{hilbert-manifold}. In terms of
corresponding Bers coordinates the mapping $L_{[\mu]}$ takes the
form
\begin{align*}
\beta\circ \Phi\circ R_{(\mu*\kappa)^{-1}}\circ L_{\mu} \circ
R_{\kappa} \circ \Lambda: \beta(V_0)\subset A_2(\Del) \rightarrow
A_2(\Del).
\end{align*}
Since left translations commute with right translations, it
simplifies to
\begin{align*}
\beta\circ \Phi\circ R_{\mu^{-1}}\circ L_{\mu}\circ \Lambda :
\beta(V_0)\subset A_2(\Del) \rightarrow A_2(\Del),
\end{align*}
which does not depend on $[\kappa]$. Thus to show that $L_{[\mu]}$
is continuous at $[\kappa]\in T_0(1)$ it is sufficient to show
that the above mapping is continuous at $0\in A_2(\Del)$. From
Lemma \ref{right1} it follows that the mapping
\begin{displaymath}
R_{\mu^{-1}}: L^2(\Del^*, \rho(z)d^2z)\cap L^{\infty}(\Del^*)_1
\rightarrow L^2(\Del^*, \rho(z)d^2z)\cap L^{\infty}(\Del^*)_1
\end{displaymath}
is continuous at $\mu$, and from Lemma \ref{rightinv} it follows
that the mapping
\begin{displaymath}
\beta\circ \Phi: L^2(\Del^*, \rho(z)d^2z)\cap L^{\infty}(\Del^*)_1
\rightarrow A_2(\Del)
\end{displaymath}
is continuous at $0$. What remains to show is that the mapping
\begin{gather*}
L_{\mu}: H^{-1,1}(\Del^*)\cap
L^{\infty}(\Del^*)_1\rightarrow L^2(\Del^*,\rho(z)d^2z),\\
\nu\mapsto \mu*\nu =\frac{(\mu \circ w_{\nu})
\frac{\ov{(w_{\nu})_z}}{(w_{\nu})_z}+\nu }{1+ (\mu \circ
w_{\nu})\bar{\nu} \frac{\ov{(w_{\nu})_z}}{(w_{\nu})_z}}
\end{gather*}
is continuous at $0$.

For $\nu\in H^{-1,1}(\Del^{*})\cap L^{\infty}(\Del^{*})_{1}$ we
have
\begin{align*}
L_{\mu}(\nu) - L_{\mu}(0)= \frac{(\mu \circ w_{\nu})
\frac{\ov{(w_{\nu})_z}}{(w_{\nu})_z}  -\mu + \nu -(\mu \circ
w_{\nu}) \mu \bar{\nu} \frac{\ov{(w_{\nu})_z}}{(w_{\nu})_z} }{1+
(\mu \circ w_{\nu})\bar{\nu}\frac{\ov{(w_{\nu})_z}}{(w_{\nu})_z}}.
\end{align*}
Since $\mu \in L^2(\Del^*, \rho(z)d^2z)$, for every $\vep>0$ there
exists $0<r<1$ such that
\begin{align*}
\iint\limits_{\Del \setminus \Del_r} \vert \mu \vert^2 \rho(z)d^2z
<\frac{\vep^2}{16}.
\end{align*}
It follows from Mori's theorem (see, e.g., \cite{Ahl2}) that for
every $0<\delta_{1}<1$ there exists $0<r_{0}<1$ such that for all
$\Vert \nu \Vert_{\infty}< \delta_1$ and $r_{0}\leq r<1$,
\begin{align*}
w_{\nu} (\Del \setminus \Del_{r'}) \subset \Del \setminus \Del_r,
\quad\text{for all}\quad r'\geq\frac{1+r}{2}.
\end{align*}
Hence if $\Vert \nu \Vert_{\infty} < \min \{\delta, \delta_1\}$,
where $\delta$ is as in Corollary \ref{righttran}, we have
\begin{gather*}
\iint\limits_{\Del} \left|\mu \circ w_{\nu}
\frac{\ov{(w_{\nu})_z}}{(w_{\nu})_z}-\mu\right|^2 \rho(z)d^2z \leq
2\iint\limits_{\Del\setminus \Del_{r'}} \left| \mu \circ w_{\nu}
\frac{\ov{(w_{\nu})_z}}{(w_{\nu})_z}\right|^2 \rho(z)d^2z \\+
2\iint\limits_{\Del\setminus \Del_{r'}} \left| \mu
\right|^2\rho(z)d^2z +\iint\limits_{\Del_{r'}} \left|\mu \circ
w_{\nu}
\frac{\ov{(w_{\nu})_z}}{(w_{\nu})_z}-\mu\right|^2 \rho(z)d^2z \\
\leq \frac{3\vep^2}{8} +\iint\limits_{\Del_{r'}} \left|\mu \circ
w_{\nu} \frac{\ov{(w_{\nu})_z}}{(w_{\nu})_z}-\mu\right|^2
\rho(z)d^2z.
\end{gather*}
Since $w_{\nu}(z) \rightarrow z$, $(w_{\nu})_z(z)\rightarrow 1$ as
$\Vert\nu\Vert_{\infty} \rightarrow 0$ uniformly for all $z\in
\Del_{r'}$, we can choose $\delta_2$ such that for all $\Vert \nu
\Vert_{\infty} <\delta_2$, the second integral is less than
$\tfrac{\vep^2}{8}$. Finally, let $\delta_3\leq
\sqrt{\tfrac{4\pi}{3}}\min\{\delta, \delta_1, \delta_2\}$, then
for all $\Vert \nu \Vert_2 <\delta_3$, we have
\begin{align*} \bigl\Vert L_{\mu}(\nu) - L_{\mu}(0) \bigr\Vert_2 <
\frac{1}{1-\sqrt{\tfrac{3}{4\pi}}\, \delta_3} \left(
\frac{\vep}{\sqrt{2}} + 2\delta_3\right),
\end{align*}
 which is less than
$\vep$ for $\delta_3$ sufficiently small.
\end{proof}

Notice that what we actually prove  is the following
\begin{corollary}\label{adjoint}
For every $[\mu]\in T_{0}(1)$ the adjoint mapping
$$T_0(1)\ni [\nu]\mapsto [\mu*\nu*\mu^{-1}]\in T_0(1)$$
is continuous at the origin.
\end{corollary}

\begin{lemma}\label{inverse} The following statements hold.
\begin{itemize}
\item[(i)] The inverse mapping $T_0(1)\ni [\mu]\mapsto
[\mu^{-1}]\in T_0(1)$ is continuous at the origin. \item[(ii)] The
group multiplication $T_{0}(1)\times T_0(1) \ni ([\mu],
[\nu])\mapsto [\mu*\nu]\in T_{0}(1)$ is continuous at the origin.

\end{itemize}
\end{lemma}
\begin{proof} In coordinate chart $V_0$ the inverse mapping takes the
form
\begin{align*}
A_2(\Del) \supset\beta(V_0)\ni\phi &\mapsto \beta\circ \Phi
((\Lambda\phi)^{-1})\in A_2(\Del).
\end{align*}
By Lemma \ref{rightinv}, it is continuous at $0$. This proves (i).
In coordinate chart $V_0$ the group multiplication takes the form
\begin{align*}
\beta(V_0)\times \beta(V_0)\ni(\phi_1, \phi_2)\mapsto \beta\circ
\Phi (\Lambda(\phi_1)*\Lambda(\phi_2))\in A_2(\Del).
\end{align*}
Its continuity at $(0,0)\in A_2(\Del)\times A_2(\Del)$ follows
from Lemma \ref{rightinv} and the proof of Lemma \ref{leftcont}.

\end{proof}

\begin{theorem} The Hilbert manifold $T_0(1)$ is a topological group.
\end{theorem}
\begin{proof}
We need to prove that the map
\begin{align*}
T_0(1)\times T_{0}(1)\ni ([\mu], [\nu])\mapsto [\mu *\nu^{-1}]\in
T_{0}(1)
\end{align*}
is continuous at every point $([\mu], [\nu])\in T_0(1)\times
T_0(1)$ with representatives $\mu, \nu \in L^2(\Del^*,
\rho(z)d^2z)\cap \mathcal{O}(\Del^*)_1$. In coordinate charts
$V_{\mu}$, $V_{\nu}$ and $V_{\kappa}$, where
$\kappa=\mu*\nu^{-1}$,  this map takes the form
\begin{align*}
\beta(V_0)\times \beta(V_0) \ni(\phi_1, \phi_2) \mapsto \beta\circ
\Phi(\Lambda(\phi_1)* \kappa* \Lambda(\phi_2)^{-1}*\kappa^{-1})\in
A_2(\Del).
\end{align*}
Its continuity at $(0,0)\in A_2(\Del)\times A_2(\Del)$ follows
from Corollary \ref{adjoint}, Lemma \ref{inverse} and Lemma
\ref{rightinv}.
\end{proof}

Finally, we show that $\mathcal{T}_0(1)$ is also a topological
group.

\begin{proposition}
The Hilbert manifold $\mathcal{T}_0(1)$ is a topological group.
\end{proposition}
\begin{proof}
We have to prove that the multiplication map
\begin{displaymath}
\mathcal{T}_0(1)\times\mathcal{T}_0(1)\ni (\gamma_{1},\gamma_{2}) \mapsto \gamma_{1}\ast\gamma_{2}^{-1}\in \mathcal{T}_0(1)
\end{displaymath}
is continuous. Setting $\gamma_{1}=([\nu],\zeta)$ and $\gamma_{2}=([\mu],w)$, we get from \eqref{gp2} that  $\gamma_{1}\ast\gamma_{2}^{-1}=([\lambda],z)$, where
$$\lambda=(\sigma_{u}^{-1})^{\ast}(\nu\ast\mu^{-1}),\quad z=(w^{\lambda}\circ\gamma_{1}\circ (w^{\nu})^{-1})(\zeta),$$
and
$$\gamma_{2}=\sigma_{u}\circ w_{\mu}=g_{\mu}^{-1}\circ f^{\mu},\quad u=g_{\mu}^{-1}(w).$$
It easily follows from the proof of Lemma A.5 that the mapping
$$\Mob(S^{1})\times T(1)\ni (\sigma,[\mu])\mapsto [(\sigma^{-1})^{\ast}(\mu)]\in T(1)$$
is continuous. Now the map $\mathcal{T}(1)\ni ([\mu],w)\mapsto g_{\mu}^{-1}(w)\in\Del^{\ast}$ depends continuously on
$[\mu]$ (uniformly when $w$ varies on compact subsets), so using that $T_{0}(1)$ is a topological group, we conclude that $[\lambda]\in T_{0}(1)$ depends continuously on $\gamma_{1}, \gamma_{2}
\in\mathcal{T}_{0}(1)$. Finally, by the result of Bers \cite{Bers1} (which we already have used in the proof of Lemma \ref{right-invariance}) we get that $z\in\Del^{\ast}$ also depends continuously on $\gamma_{1}$ and $\gamma_{2}$.
 \end{proof}

\section{}
Here we prove Lemma \ref{second-in-L^2}.
Using the model $\HH\simeq\UU$ and \eqref{ahlfors-K-bar}, we get for
the second variation of $Q$,
\begin{align*}
&\left.\frac{\pa^2}{\pa \vep\pa
\bar{\vep}}\right\vert_{\vep=0}Q(\mu,
\vep\nu)(z)\\
=&\frac{24}{\pi^3}\iint\limits_{\U}\iint\limits_{\U}\iint\limits_{\U}
\frac{\rho(z)^{-1}\mu(w)\ov{\nu(u)}\nu(v)d^2vd^2ud^2w}
{(w-v)^2(v-\z)^2(w-\bar{u})^2(\bar{u}-\z)^2} \\
+&\frac{24}{\pi^3}\iint\limits_{\U}\iint\limits_{\U}\iint\limits_{\U}
\frac{\rho(z)^{-1}\mu(w)\ov{\nu(u)}\nu(v)d^2vd^2ud^2w}
{(w-\z)^2(w-v)^2(v-\bar{u})^2(\bar{u}-\z)^2}\\
+&\frac{24}{\pi^3}\iint\limits_{\U}\iint\limits_{\U}\iint\limits_{\U}
\frac{\rho(z)^{-1}\mu(w)\ov{\nu(u)}\nu(v)d^2vd^2ud^2w}
{(w-\z)^2(w-\bar{u})^2(v-\bar{u})^2(v-\z)^2} \\
-&\mu(z)\rho^{-1}(z)
\left.\frac{\pa^2}{\pa
\vep\pa\bar{\vep}}\right\vert_{\vep=0}\rho^{\vep\nu}(z)
=I_1(z)+I_2(z)+I_3(z)+I_4(z).
\end{align*}
Using Proposition \ref{secondvarmetric} and property \textbf{RK3},
we immediately get that $I_4\in L^2(\UU,\rho(z)d^2z)$. Now using
the Cauchy-Schwarz inequality, the identity
\begin{displaymath}
\iint\limits_{\U}\frac{d^2w}{|w-\z|^4}=\frac{\pi}{4}\rho(z),
\end{displaymath}
and the property that the Hilbert transform is an isometry on $L^2(\CC,d^2z)$, we get
\begin{align*}
&\Vert I_1\Vert^2_2= \frac{24^2}{\pi^6}
\iint\limits_{\U}\left|\iint\limits_{\U}\iint\limits_{\U}\iint\limits_{\U}
\frac{\mu(w)\ov{\nu(u)}\nu(v)d^2vd^2ud^2w}
{(w-v)^2(v-\z)^2(w-\bar{u})^2(\bar{u}-\z)^2}\right|^2\rho(z)^{-1}d^2z\\
\leq & \frac{24^2}{\pi^6}\iint\limits_{\U} \left(\iint\limits_{\U}
\frac{\rho(z)^{-1}|\nu(v)|^2d^2v}{|v-\z|^4}\right) \\
\times &
\left(\iint\limits_{\U}
\left|\iint\limits_{ \U}\iint\limits_{\U}
\frac{\mu(w)\ov{\nu(u)}d^2ud^2w}
{(w-v)^2(w-\bar{u})^2(\bar{u}-\z)^2}\right|^2d^2v\right)d^2z\\
\leq &\frac{12^2\Vert
\nu\Vert_{\infty}^2}{\pi^5}\iint\limits_{\U}\iint\limits_{\U}\left|
\iint\limits_{\U}\iint\limits_{\U}\frac{\mu(w)\ov{\nu(u)}d^2ud^2w}
{(w-v)^2(w-\bar{u})^2(\bar{u}-\z)^2}\right|^2d^2vd^2z\\
\leq & \frac{12^2\Vert
\nu\Vert_{\infty}^2}{\pi^3}\iint\limits_{\U}\iint\limits_{\U}|\mu(v)|^2
\left|\iint\limits_{\U}\frac{\nu(u)d^2u}{(u-\bar{v})^2(u-z)^2}\right|^2
d^2zd^2v\leq 36\Vert \nu\Vert_{\infty}^4\Vert
\mu\Vert_2^2.
\end{align*}
Similarly, denoting by $\ov{\U}$ the lower half-plane,
\begin{align*}
\Vert  I_3\Vert_2^2 &=\frac{24^2}{\pi^6}
\iint\limits_{\U}\left|\iint\limits_{\U}\iint\limits_{\U}\iint\limits_{\U}
\frac{\mu(w)\ov{\nu(u)}\nu(v)d^2vd^2ud^2w}
{(w-\z)^2(w-\bar{u})^2(v-\bar{u})^2(v-\z)^2}\right|^2\rho(z)^{-1}d^2z\\
&\leq \frac{12^2\Vert
\nu\Vert_{\infty}^2}{\pi^5}\iint\limits_{\U}\iint\limits_{\U}\left|
\iint\limits_{\U}\iint\limits_{\U}\frac{\mu(w)\ov{\nu(u)}d^2ud^2w}
{(w-\z)^2(w-\bar{u})^2(v-\bar{u})^2}\right|^2d^2vd^2z\\
&= \frac{12^2\Vert
\nu\Vert_{\infty}^2}{\pi^5}\iint\limits_{\U}\iint\limits_{\ov{\U}}\left|
\iint\limits_{\U}\iint\limits_{\U}\frac{\mu(w)\ov{\nu(u)}d^2ud^2w}
{(w-\z)^2(w-\bar{u})^2(\bar{v}-\bar{u})^2}\right|^2d^2vd^2z\\
&\leq  \frac{12^2\Vert
\nu\Vert_{\infty}^2}{\pi^3}\iint\limits_{\U}\iint\limits_{\U}
|\nu(v)|^2 \left|\iint\limits_{\U}\frac{\mu(w)d^2w}
{(w-\z)^2(w-\bar{v})^2}\right|^2d^2vd^2z\\
&\leq\frac{12^2\Vert
\nu\Vert_{\infty}^4}{\pi}\iint\limits_{\U}\iint\limits_{\U}\frac{|\mu(v)|^2}{|v-\z|^4}d^2vd^2z=36\Vert
\nu\Vert_{\infty}^4\Vert\mu\Vert^2_2.
\end{align*}
For the term $I_2$ we use the same identity as in the proof of
Proposition \ref{secondvarmetric}\footnote{We could estimate $I_2$
in the same way as $I_1$ If $\nu\in H^{-1,1}(\U)$. However, for
Theorem \ref{WPvariation-Gamma} we only have
$\nu\in\Omega^{-1,1}(\U)$.}
\begin{align*}
\frac{1}{(v-\bar{u})^2(\bar{u}-\z)^2} =&\frac{(v-z)^2}{(v-\z)^2
(z-\bar{u})^2(v-\bar{u})^2}+\frac{2(v-z)^2(z-\z)}{(v-\z)^3(v-\bar{u})(z-\bar{u})^3}\\
+&\frac{(z-\z)^2}{(v-\z)^2(z-\bar{u})^2(\bar{u}-\z)^2}+\frac{2(z-\z)^2(v-z)}{(v-\z)^3(\bar{u}-\z)(\bar{u}-z)^3}.
\end{align*}
As in the proof of Proposition \ref{secondvarmetric}, the last two
terms do not contribute to $I_2$, and we obtain
\begin{align*}
I_2(z)&=\frac{24}{\pi^3}\rho^{-1}(z)
\iint\limits_{\U}\iint\limits_{\U}\iint\limits_{\U}\frac{\mu(w)\nu(v)\ov{\nu(u)}(v-z)^2d^2vd^2ud^2w}
{(w-\z)^2(w-v)^2(v-\z)^2(z-\bar{u})^2(v-\bar{u})^2}\\
&+\frac{48}{\pi^3}\rho^{-1}(z)\iint\limits_{\U}\iint\limits_{\U}\iint\limits_{\U}
\frac{\mu(w)\nu(v)\ov{\nu(u)}(v-z)^2(z-\z)d^2vd^2ud^2w}
{(w-\z)^2(w-v)^2(v-\z)^3(z-\bar{u})^3(v-\bar{u})} \\
&=I_5(z)+I_6(z).
\end{align*}
The $L^2$-norm of $I_5$ is estimated exactly as before and we get
$\Vert I_5\Vert^2_2\leq 36\Vert\nu\Vert_{\infty}^4\Vert \mu\Vert_2^2$.
Finally,
\begin{align*}
\Vert I_6\Vert_2^2 =&\frac{48^2}{\pi^6} \iint\limits_{\U}\left|
\iint\limits_{\U}\iint\limits_{\U}\iint\limits_{\U}\frac{\mu(w)\nu(v)
\ov{\nu(u)}(v-z)^2(z-\z)d^2vd^2ud^2w}
{(w-\z)^2(w-v)^2(v-\z)^3(z-\bar{u})^3(v-\bar{u})}\right|^2\rho(z)^{-1}d^2z\\
\leq & \frac{24^2\Vert
\nu\Vert_{\infty}^2}{\pi^5}\iint\limits_{\U}\iint\limits_{\U}\frac{|z-\z|^2}
{|u-\z|^2} \left|\;\;\iint\limits_{\U}\iint\limits_{\U}
\frac{\mu(w)\nu(v)(v-z)^2d^2vd^2w}{(w-\z)^2(w-v)^2(v-\z)^3(v-\bar{u})}\right|^2d^2ud^2z\\
\leq & \frac{24^2\Vert
\nu\Vert_{\infty}^2}{\pi^5}\iint\limits_{\U}
\left(\iint\limits_{\U}\left|\iint\limits_{\U}
\frac{\mu(w)d^2w}{(w-\z)^2(w-v)^2}\right|^2d^2v\right)\\
& \hspace{2cm}\left(\iint\limits_{\U}\frac{|z-\z|^2}{|u-\z|^2}
\iint\limits_{\U}\frac{|\nu(v)|^2|v-z|^4}{|v-\z|^6|v-\bar{u}|^2}d^2v
d^2u\right)d^2z.
\end{align*}
Making a change of variables $\frac{v-z}{v-\z}\mapsto v$ and
$\frac{u-z}{u-\z}\mapsto u$, we obtain
\begin{align*}
&\iint\limits_{\U}\frac{|z-\z|^2}{|u-\z|^2}
\iint\limits_{\U}\frac{|\nu(v)|^2|v-z|^4}{|v-\z|^6|v-\bar{u}|^2}d^2v
d^2u & \leq \Vert\nu\Vert_{\infty}^2\iint\limits_{\Del}\iint\limits_{\Del}
\frac{|v|^4}{|1-v\bar{u}|^2}d^2vd^2u \\
&=\frac{3\pi^2}{4}\Vert\nu\Vert_{\infty}^2.
\end{align*}
Hence
\begin{align*}
\Vert I_6\Vert_2^2 \leq \frac{3 \cdot 12^2\Vert
\nu\Vert_{\infty}^4}{\pi}\iint\limits_{\U}\iint\limits_{\U}
\frac{|\mu(v)|^2}{|v -\z|^4}d^2vd^2z=
108\Vert\nu\Vert_{\infty}^4 \Vert\mu\Vert_2^2.
\end{align*}

\bibliographystyle{amsalpha}
\bibliography{kahler}

\end{document}